\documentclass[preprint,12pt]{elsarticle}

\usepackage[T1]{fontenc}
\usepackage{lmodern}
\usepackage[utf8]{inputenc}
\usepackage{amsmath,amssymb,amsthm}
\usepackage{mathtools}
\usepackage{enumitem}
\usepackage{graphicx}
\usepackage{booktabs}
\usepackage{array}
\usepackage{longtable}
\usepackage[margin=2.6cm]{geometry}
\usepackage{xr-hyper}
\usepackage[colorlinks=true,linkcolor=blue!50!black,citecolor=blue!50!black,urlcolor=blue!50!black]{hyperref}
\usepackage[capitalise,nameinlink]{cleveref}

\newcommand{\C}{\mathbb{C}}
\newcommand{\R}{\mathbb{R}}
\newcommand{\Z}{\mathbb{Z}}
\newcommand{\N}{\mathbb{N}}
\newcommand{\T}{\mathbb{T}}
\newcommand{\HH}{\mathbb{H}}
\newcommand{\RE}{\operatorname{Re}}
\newcommand{\IM}{\operatorname{Im}}
\newcommand{\Arg}{\operatorname{Arg}}
\newcommand{\Log}{\operatorname{Log}}
\newcommand{\dist}{\operatorname{dist}}
\newcommand{\rank}{\operatorname{rank}}
\newcommand{\Ran}{\operatorname{Ran}}
\newcommand{\Var}{\operatorname{Var}}
\newcommand{\Kr}{\mathcal{K}}
\newcommand{\Sa}{\Sigma_\alpha}
\newcommand{\La}{\Lambda_\alpha}
\newcommand{\B}{\mathcal{B}}
\newcommand{\RR}{\mathcal{R}}
\newcommand{\Jq}{\mathcal{J}}
\newcommand{\Mst}{M^{*}_{\alpha}}

\newcommand{\D}{\mathcal{D}}
\newcommand{\eps}{\varepsilon}

\newtheorem{theorem}{Theorem}[section]
\newtheorem{lemma}[theorem]{Lemma}
\newtheorem{proposition}[theorem]{Proposition}
\newtheorem{corollary}[theorem]{Corollary}

\theoremstyle{definition}
\newtheorem{definition}[theorem]{Definition}

\theoremstyle{remark}
\newtheorem{remark}[theorem]{Remark}
\theoremstyle{plain}
\newtheorem{introthm}{Theorem}

\crefname{introthm}{Theorem}{Theorems}
\crefformat{appendix}{#2#1#3}
\Crefformat{appendix}{#2#1#3}
\crefrangeformat{appendix}{#3#1#4 to #5#2#6}

\usepackage[textsize=tiny]{todonotes}
\usepackage{xcolor}

\journal{Journal of Functional Analysis}
\biboptions{sort&compress}

\begin{document}

\begin{frontmatter}

\title{A Kreiss theorem for fractional Cauchy problems \\ with a finite slow component}

\author[ns]{Vladimir R.\ Kosti\'c}
\ead{vladimir.slk@gmail.com}

\author[impa]{Jo\~ao Pedro Mangi}
\ead{mangi.joao@gmail.com}

\author[mi]{Ilja Uzelac Buji\v{s}i\'c\corref{cor}}
\ead{bujisiciljamat@gmail.com}
\cortext[cor]{Corresponding author.}

\affiliation[ns]{organization={University of Novi Sad},
            addressline={Trg Dositeja Obradovi\'ca 4},
            postcode={21000},
            city={Novi Sad},
            country={Serbia}}
\affiliation[impa]{organization={Instituto de Matem\'atica Pura e Aplicada},
            addressline={Estrada Dona Castorina 110},
            postcode={22460-320},
            city={Rio de Janeiro},
            country={Brazil}}
\affiliation[mi]{organization={Mathematical Institute of the Serbian Academy of Sciences and Arts},
            addressline={Kneza Mihaila 36},
            postcode={11000},
            city={Belgrade},
            country={Serbia}}

\begin{abstract}
{Fractional Cauchy problems with a Caputo derivative of order bounded by one are solved by Mittag--Leffler operator families that lack the semigroup property but remain bounded for classes of sectorial generators, even  with spectrum in the right half-plane. To study their transient growth we introduce the fractional Kreiss constant via the Laplace transform of the propagator family. For a closed generator on a Banach space with a time-scale separation and a finite spectral
component, that is, a finite-rank projector commuting with the generator whose complement carries a bounded
family, we prove that the projected part of the propagator (typically slow spectral part carrying dynamics that is pertinent at mid and long time-scales) is comparable to a deflated Kreiss constant up
to a factor linear in the rank, uniformly in fractional order and independently of the dimension. The
analytic core of this result is a Bernstein-type inequality for a class of non-rational functions on vertical lines proved by bounding the turning of the derivative. This is a fundamentally different technique from the classical one used in Kreiss theorem for Markov semigroups based on counting intersections of algebraic
curves. Importantly, we 
show that the rank dependence cannot be removed while keeping the Kreiss constant upper bound. Finally, we show that in Hilbert spaces, generators whose numerical range avoids the sector of instability obey the sharp bound of normal operators, which no spectral-set estimate improves. To conclude, we discuss how subdiffusion driven by a self-excited oscillator, the time-fractional Schr\"odinger equation, and time-changed Markov dynamics illustrate the theory.}
\end{abstract}


\end{frontmatter}

\section{Introduction}\label{sec:intro}


The stability of a linear evolution equation tells us what happens eventually, while its transient 
tells us what happens on the way, and for non-normal generators the two can be very far apart.  In a standard finite dimensional case, i.e. for a
matrix $A\in\C^{d\times d}$ the transient growth is bounded by the Kreiss constant
$\Kr_1(A):=\sup_{\RE z>0}\RE(z)\,\|(z-A)^{-1}\|$, and the continuous-time Kreiss matrix theorem
\begin{equation}\label{eq:classical}
\Kr_1(A)\;\le\;\sup_{t\ge0}\|e^{tA}\|\;\le\;e\,d\,\Kr_1(A)
\end{equation}
states that the peak of the transient growth is controlled from both sides by a single resolvent quantity,
with a factor linear in the dimension
\citep{Kreiss1959,kreiss1962stabilitatsdefinition,LeVequeTrefethen1984,Spijker1991,DorsselaerKraaijevangerSpijker1993,trefethen2005spectra}.
This two-sided estimate is the backbone of the pseudospectral theory of transient growth
\citep{trefethen2005spectra}, and its infinite-dimensional limits are by now fully understood: on
Banach spaces the Kreiss resolvent condition gives no bound on the semigroup at all, while on Hilbert
spaces it gives sublinear growth in time \citep{EisnerZwart2006,Arnold2022,DellOro2025,Arnold2026}. 

The purpose of this paper is to develop the corresponding theory for fractional Cauchy problems
\begin{equation}\label{eq:caputo}
\partial_t^\alpha u=Lu,\qquad u(0)=u_0,\qquad 0<\alpha\le1,
\end{equation}
with the Caputo derivative $\partial_t^\alpha$ and a closed operator $L$ on a Banach space $X$. Such
problems describe evolutions with memory: subdiffusive transport in heterogeneous or crowded media
arises as the scaling limit of continuous-time random walks with heavy-tailed waiting times
\citep{MetzlerKlafter2000,MeerschaertScheffler2004,MeerschaertSikorskii2019}, the time-fractional
Schr\"odinger equation was proposed as a quantum model with memory \citep{Naber2004}, and fractional
relaxation and oscillation equations are standard models in viscoelasticity and 
control \citep{Podlubny1999,Matignon1996,Gorenflo2020}. The solution of \eqref{eq:caputo} is
$u(t)=E_\alpha(t^\alpha L)u_0$, where $E_\alpha(t^\alpha L)$ is the Mittag--Leffler family generated by
$L$ \citep{Pruss2012,Bazhlekova2001}, given for bounded or sectorial $L$ by the Mittag--Leffler function
$E_\alpha(z)=\sum_{k\ge0}\frac{z^k}{\Gamma(1+\alpha k)}$.

Already at the level of definitions, these families differ from semigroups in ways that call for a
theory of their own. The map $t\mapsto E_\alpha(t^\alpha L)$ has no semigroup property, so arguments
that propagate a bound from short to long times are not available. The stability region is larger: for
a matrix $L$ the family is bounded if and only if the spectrum avoids the open sector
$\Sa:=\{z\ne0:|\Arg z|<\frac{\alpha\pi}{2}\}$ and the eigenvalues on its boundary are semisimple
\citep{Matignon1996}. Since $\Sa$ is a proper subsector of the right half-plane when $\alpha<1$, modes
that grow exponentially in physical time can be perfectly stable under the fractional clock, and the
subordination formula, which writes $E_\alpha(t^\alpha L)$ as an average of $e^{sL}$ over a random time
\citep{Bazhlekova2001}, then gives no information, because it averages an exponentially growing
semigroup. 

The envelope $t\mapsto\|E_\alpha(t^\alpha L)\|$ in general has a short-time, an intermediate and a long-time regime, each governed by a
different resolvent quantity. In particular, in Hilbert spaces for the semigroup case $\alpha=1$, it is well understood that numerical and spectral abscissa govern short-time and long-time behavior, respectively, see e.g. \cite{trefethen2005spectra}. While this remains also the case in the fractional setting (we review these essentially known results in \Cref{app:timescales}), the characterization of intermediate behavior is much more involved problem. Indeed, whenever the operator $L$ is normal, while the envelope of the semigroup $\Vert e^{t L}\Vert$ is sharply bounded by one, in the fractional setting  $|E_\alpha|$ exceeds ${1}/{\alpha}$ at
points of the boundary rays of $\Sa$, implying that even a normal generator with spectrum on these rays amplifies envelope by more than ${1}/{\alpha}$, thus provoking transient growth. Hence, in this paper we focus on answering the question:

\begin{quote}
\textit{Which resolvent quantity controls $\sup_{t\ge0}\|E_\alpha(t^\alpha L)\|$, and with which constants?}    
\end{quote}

In answering it, we particularly focus on the finite-rank operators on infinite-dimensional spaces, which is motivated by recent works on data-driven estimators of
transfer operators and generators of dynamical systems
\citep{KlusEtAl2020,kostic2022koopman,kostic2023sharp,kostic2025laplace}, including estimators of slow modes of molecular systems \citep{devergne2025slow}, and their pseudospectra can be computed at a cost
governed by the rank \citep{kostic2026lowrank}. 

\subsection{Related work}

The continuous-time Kreiss matrix theorem goes back to \citet{Kreiss1959,kreiss1962stabilitatsdefinition};
the linear dependence on the dimension was obtained by \citet{LeVequeTrefethen1984}, the constant $e$
through Spijker's  
\citep{Spijker1991,DorsselaerKraaijevangerSpijker1993,wegert1994buffon}, and
\citet{Kraaijevanger1994} showed that the linear order is sharp, for powers as well as for the
exponential; see also \citet{SpijkerTracognaWelfert2003}, and for general domains \citet{TohTrefethen1999}.
\citet{EisnerZwart2006} showed that the Kreiss condition gives no growth bound on Banach spaces, that
it gives $\|e^{tL}\|=O(t)$ on Hilbert spaces, and that growth of order $t^\gamma$ occurs for every
$\gamma<1$, using Fourier multipliers on $L^2$ spaces with power weights of class $A_2$. The
Hilbert-space bound was improved to $O(\frac{t}{\sqrt{\log t}})$ by \citet{Arnold2022}, extended to general
resolvent growth by \citet{DellOro2025} and \citet{DellOroFarinaPata2026}, and very recently to
$O(t^{1-\eps})$ with $\eps$ depending on the Kreiss constant \citep{Arnold2026}; see also
\citet{RozendaalVeraar2018}. All these arguments use the semigroup law. On the side of spectral sets,
 the numerical range is a $(1+\sqrt2)$-spectral set \citep{CrouzeixPalencia2017}, von Neumann's
inequality makes a half-plane containing the numerical range a spectral set with constant one 
\citep{vonNeumann1951,SzNagyFoias2010}, and \citet{Crouzeix2025} interpolates between the two for
unbounded convex sets.

Mittag--Leffler families, also called $\alpha$-times resolvent families or solution operators, are the
resolvent families of Volterra equations with the kernel $\frac{t^{\alpha-1}}{\Gamma(\alpha)}$
\citep{Pruss2012,Bazhlekova2001}. Their generation theory, subordination and asymptotics are
classical \citep{Pruss2012,Bazhlekova2001,Gorenflo2020}, and resolvent characterisations of growth
bounds and of uniform stability were obtained recently by \citet{GhavidelShahali2023} and
\citet{MeiLi2024}, in the spirit of the Hille--Yosida and Gearhart--Pr\"uss theorems. These results are
qualitative or asymptotic, and we are not aware of quantitative two-sided bounds of Kreiss type for
Mittag--Leffler families beyond the matrix case conjectured in 
\citet{KosticBujisic2025}. For matrices,
the spectral criterion of \citet{Matignon1996} characterises boundedness, and fractional pseudospectra
have been used to localise spectra and to study stability of fractional systems
\citep{SancaKosticCvetkovic2018,DiethelmEtAl2024}.

\subsection{Organisation}

For the reader's convince, we start with the summerize of our main results as four Theorems \ref{thm:A}-\ref{thm:D}  in \Cref{sec:nutshell}. Then, in \Cref{sec:prelim} we review  facts on the Mittag--Leffler functions, Mittag--Leffler operator families, and prove basic properties of the newly introduced the
fractional Kreiss constant. Based on these preliminary results, \Cref{sec:bernstein} proves \Cref{thm:B}, and
\cref{sec:finite-rank} the finite-rank core of \Cref{thm:A}. \Cref{sec:hilbert} treats normal
operators and the numerical range, and \cref{sec:slow} completes \Cref{thm:A} for generators with a
slow component and analyses prototype (a). \Cref{sec:optimal} proves \Cref{thm:C},
\cref{sec:mechanisms} describes the mechanisms and time scales of the transient behavior,
\cref{sec:experiments} revisits the three prototypes introduced in \cref{sec:nutshell}, and \cref{sec:conclusions} discusses open
problems. The appendix contains the table of notation and the proofs of auxiliary results, and the
Supplementary Material the numerical methods, the complete tables, the study of random finite-rank
operators and a description of the code.

\section{Main contributions in a nutshell}\label{sec:nutshell}

Recalling that the Laplace transform of $E_\alpha(t^\alpha L)$ is $\omega^{\alpha-1}(\omega^\alpha-L)^{-1}$, somewhat natural way to define the \emph{fractional Kreiss constant} is
\begin{equation}\label{eq:kreiss-def}
\Kr_\alpha(L):=\sup_{\RE\omega>0}\RE(\omega)\,\big\|\omega^{\alpha-1}(\omega^\alpha-L)^{-1}\big\|
=\sup_{z\in\Sa}a_\alpha(z)\,\|(z-L)^{-1}\|,
\end{equation}
where $a_\alpha(z):=|z|\cos\big(\frac{\Arg z}{\alpha}\big)$; it reduces to $\Kr_1$ for $\alpha=1$. The second
expression replaces the classical weight $\RE z$ on the right half-plane by a weight on $\Sa$ that
vanishes on the boundary rays, so $\Kr_\alpha(L)$ is finite for operators that are unstable as generators
of semigroups. Up to a factor between $1$ and $\frac{1}{\alpha}$, $a_\alpha(z)$ is the distance from $z$ to the boundary
of $\Sa$ (\cref{lem:elementary}(c)), so $\Kr_\alpha(L)$ measures how fast the resolvent may blow up as $z$
approaches the stability region from inside the sector of instability. A constant of this type, for matrices and in pseudospectral form, was proposed in
the precursor \citep{KosticBujisic2025}, 
%
where a fractional Kreiss matrix theorem was conjectured; the
present paper proves this conjecture in a stronger form, for all real $\alpha\in(0,1]$, on Banach spaces and with the rank in place of the dimension.
Our main object is a generator with a finite slow component, i.e. a system whose mid and long-term dynamics is determined by finitely many eigenvalues that are well separated from the fast time-scales. More precisely, a projection $P\in\B(X)$ of finite rank
$m$ is a \emph{slow component} of $L$ if it commutes with $L$, the part of $L$ in $(I-P)X$ generates a bounded Mittag--Leffler family, spectra of $P L$ and $(I-P)L$ are disjoint with the former being closer to stability boundary.

As we first show the transient growth of such a generator is carried out by
the $m$ slow modes and by the geometry of the splitting, and both are seen by the
\emph{deflated Kreiss constant} $\Kr_\alpha(L;P):=\sup_{z\in\Sa}a_\alpha(z)\,\|(z-L)^{-1}P\|$, the Kreiss
quantity of the compressed resolvent.
\begin{introthm}[\Cref{thm:kreiss-lower,thm:finite-rank,thm:slow}]\label{thm:A}
Let $0<\alpha\le1$, let $X\ne\{0\}$ be a Banach space, let $L$ be a closed, densely defined operator
on $X$, and put
\[
\kappa_{\alpha,n}:=\frac e2\big(3+4n(2-\alpha)\big),\qquad
\kappa'_{\alpha,n}:=\frac e2\big(1+4n(2-\alpha)\big).
\]
\begin{enumerate}[label=\textup{(\roman*)},leftmargin=*]
\item If $L$ generates a bounded Mittag--Leffler family, then
$\Kr_\alpha(L)\le\sup_{t\ge0}\|E_\alpha(t^\alpha L)\|$.
\item If $P$ is a slow component of rank $m$, then $L$ generates a Mittag--Leffler family,
\begin{gather*}
\Kr_\alpha(L;P)\le\sup_{t\ge0}\|E_\alpha(t^\alpha L)P\|\le\kappa_{\alpha,m}\,\Kr_\alpha(L;P),\\
\sup_{t\ge0}\|E_\alpha(t^\alpha L)\|\le\kappa_{\alpha,m}\,\Kr_\alpha(L;P)+\sup_{t\ge0}\|E_\alpha(t^\alpha L)(I-P)\|,
\end{gather*}
and $\kappa_{\alpha,m}$ may be replaced by $\kappa'_{\alpha,m}$ when $0$ is not an eigenvalue of the
slow block.
\item If $L$ is bounded with $\rank L\le n$, then $\sup_{t\ge0}\|E_\alpha(t^\alpha L)\|\le\kappa_{\alpha,n}\Kr_\alpha(L)$.
\end{enumerate}
\end{introthm}

The constants depend on $\alpha$ and on the rank, never on the dimension of $X$, and they are uniform
in $\alpha\in(0,1]$. 

Part (iii) is the core case of (ii): a finite-rank operator with
finite Kreiss constant splits as $X=\Ran L\oplus\ker L$, its slow component is the projection onto
$\Ran L$ along $\ker L$, and the stationary contribution of the kernel is absorbed into the estimate
without paying the norm of the oblique complementary projection. Part (ii) applies to unbounded
generators that may have infinite rank, need not generate a $C_0$-semigroup and may have spectrum in
the right half-plane, as in prototype (a), and the fast part may be controlled on Banach spaces by
sectoriality or on Hilbert spaces by the numerical range.

{Second, we note that the } upper bounds rest on a Bernstein-type inequality that can be of independent interest. Inverting the Laplace transform after
one integration by parts reduces the problem to the length of the curve
$y\mapsto h_F(a+iy)$, $h_F(\omega):=\omega^{-1}F(\omega^\alpha)$, where $F$ is a rational function of
degree at most $n$ built from a matrix element of the resolvent. Let $\RR_{\alpha,n}$ denote the
nonzero rational functions $F=\frac{F_1}{F_0}$ in lowest terms with $\deg F_1\le\deg F_0\le n$ and no poles
in $\Sa$.
\begin{introthm}[\Cref{thm:bernstein}]\label{thm:B}
Let $0<\alpha\le1$, $n\in\N_0$, $F\in\RR_{\alpha,n}$ and $a>0$. Then
\[
2\sup_{y\in\R}|h_F(a+iy)|\;\le\;\int_\R|h_F'(a+iy)|\,dy\;\le\;\pi\big(3+4n(2-\alpha)\big)\sup_{y\in\R}|h_F(a+iy)|,
\]
and the constant improves to $\pi(1+4n(2-\alpha))$ when $F(0)=0$.
\end{introthm}

For $\alpha=1$ the function $h_F$ is rational and the classical route is Spijker's lemma, which
counts the intersections of algebraic curves \citep{Spijker1991,wegert1994buffon}. For irrational
$\alpha$ the curve $\{(a+iy)^\alpha:y\in\R\}$ is not algebraic and such counts are not available. We
replace them by a bound on the total turning of $\arg h_F'$, which follows from the convexity of
the image of a vertical line under $\omega\mapsto\omega^\alpha$: its tangent turns monotonically
through an angle $(1-\alpha)\pi$, and this forces the argument of $(a+iy)^\alpha-\zeta$ to vary by at
most $(2-\alpha)\pi$ for every point $\zeta$ off the curve.

 {Next, we obtain that} rank dependence in \Cref{thm:A} cannot be removed, and  {show that} one family of weighted spaces on the
circle suffices to see it. For $0<\gamma<1$ let $L^2_\gamma$ be the space $L^2(\T,|\tan(\frac{x}{2})|^\gamma\,dx)$, and let
$\D^\alpha$ be the Weyl fractional derivative, the Fourier multiplier with symbol $(ik)^\alpha$.
\begin{introthm}[\Cref{thm:optimal,cor:optimal-n,cor:infinite,prop:weyl}]\label{thm:C}
Let $0<\alpha\le1$ and $0<\gamma<1$. The rank-$2N$ truncations $B_N$ of $N^{-\alpha}\D^\alpha$ to the
trigonometric polynomials of degree at most $N$ satisfy $\sup_N\Kr_\alpha(B_N)<\infty$ and
$\sup_{t\ge0}\|E_\alpha(t^\alpha B_N)\|\ge cN^\gamma$ on $L^2_\gamma$. Consequently the best constant in
\Cref{thm:A}\textup{(iii)} grows at least like $n^\gamma$ for every $\gamma<1$; there is a bounded
operator on a Hilbert space with $\Kr_\alpha(L)<\infty$ and $\limsup_t t^{-\gamma}\|E_\alpha(t^\alpha L)\|>0$;
and $\D^\alpha$ itself satisfies $\Kr_\alpha(\D^\alpha)<\infty$ on $L^2_\gamma$ but generates no
Mittag--Leffler family.
\end{introthm}

The growth in this construction is oscillatory: on the boundary rays
$E_\alpha((is)^\alpha)=\alpha^{-1}e^{is}$ up to an error of bounded variation, so the Mittag--Leffler family
of $\D^\alpha$ contains a translation, and translation by $\pi$ is unbounded on $L^2_\gamma$ because it
moves the zero of the weight onto its pole. For $\alpha=1$ this recovers the weighted Fourier
multipliers of \citet{EisnerZwart2006}.

 {Then, we focus on } Hilbert spaces  {where} the fast part of the dynamics is often controlled by its numerical range, and  {in this setting give a sharp answer with the fractional theory}. Let $\Mst:=\sup_{z\notin\Sa}|E_\alpha(z)|$, which is the
best uniform bound for normal generators with spectrum outside $\Sa$ and satisfies $\Mst>\frac{1}{\alpha}$ for
$\alpha<1$.

\begin{introthm}[\Cref{thm:numrange}]\label{thm:D}
Let $L$ be a closed, densely defined operator on a Hilbert space whose numerical range does not meet
$\Sa$ and whose resolvent set meets $\Sa$. Then $L$ generates a Mittag--Leffler family and
$\sup_{t\ge0}\|E_\alpha(t^\alpha L)\|\le\Mst$, with equality for suitable normal operators. No spectral-set
estimate on a convex set containing the numerical range yields a better uniform bound, whatever its
constant.
\end{introthm}

The generation statement covers operators that generate no $C_0$-semigroup, such as the generator of
prototype (b), and the proof is by separation and Sz.-Nagy's dilation theorem. The bound is consistent
with, but not covered by, Crouzeix's recent estimates for unbounded convex sets \citep{Crouzeix2025}:
its constant is that of a half-plane, where his spectral constant equals one, but in the extremal case
the functions $E_\alpha(t^\alpha\,\cdot)$ have no limit at infinity and fall outside his function class.
The last statement of \Cref{thm:D} is elementary once the right geometric lemma is in place, and it
explains why the smaller constants that \citet{Crouzeix2025} obtains for sectors cannot help for the
uniform bound: the union over $t$ of the dilated sets generates a cone on which $|E_\alpha|$ has the same
supremum as on a half-plane. Combined with \Cref{thm:A}, \Cref{thm:D} gives explicit certificates for
generators whose fast part has its numerical range outside $\Sa$.


Finally, to show that the fractional transient has a structure of its own, note that normal operators amplify by $\Mst$, a factor larger than
$\frac{1}{\alpha}$ that tends to one as $\alpha\uparrow1$. On the other hand, as  we
describe in \cref{sec:mechanisms}, strongly non-normal Jordan blocks, which produce the
dimension-dependent amplification $e\sqrt{\frac{n}{2\pi}}$ in \eqref{eq:classical}, produce for each fixed
$\alpha<1$ only a bounded amplification ratio. This is because the Mittag--Leffler distribution spreads the peak
of $\frac{s^{n-1}e^{-s}}{(n-1)!}$ over a range of times, and the bound deteriorates as $\alpha\uparrow1$, as it must, suggesting that the growth exhibited in \Cref{thm:C} is oscillatory rather than nilpotent in
origin. 
In addition to this informative Jordan block example, for the three prototypes bellow, in \cref{sec:experiments} we compare both sides of \Cref{thm:A} with the bounds that are otherwise available: the semigroup route through subordination, the
numerical-range bounds of \citet{CrouzeixPalencia2017} and \citet{Crouzeix2025}, deflation bounds that
estimate the slow block and the splitting separately, and, where a complete eigenbasis exists, the
eigenvector bound.

\begin{enumerate}[label=(\alph*),leftmargin=*]
\item \emph{Subdiffusion driven by a self-excited oscillator.} A two-dimensional oscillatory mode
with eigenvalues $e^{\pm i\theta}$, $\frac{\alpha\pi}{2}<\theta<\frac{\pi}{2}$, grows exponentially as an ordinary
differential equation but is stable under the fractional clock; when it feeds a source into a
subdiffusive medium, the generator is unbounded, it has a two-dimensional slow component, its
semigroup grows exponentially, and its numerical range meets $\Sa$.
\item \emph{The time-fractional Schr\"odinger equation.} The generator $e^{-i\frac{\alpha\pi}{2}}A_V$, with $A_V$
a self-adjoint Schr\"odinger operator, is normal and generates no $C_0$-semigroup because $A_V$ is
unbounded above, yet its Mittag--Leffler family is uniformly bounded and amplifies by exactly the
sharp constant of normal operators.
\item \emph{Time-changed Markov dynamics.} Expectations of a Markov process run on the clock of an
inverse stable subordinator evolve by $E_\alpha(t^\alpha L)$, with $L$ the generator of the
process. Linear chains of compartments are strongly non-normal, and in data-driven practice the
generator is replaced by a finite-rank estimate on an infinite-dimensional function space
\citep{KlusEtAl2020,kostic2022koopman,kostic2024learning,kostic2025laplace}, whose number of retained
modes, and not the dimension of the ambient space, is the natural parameter.
\end{enumerate}

In prototype (a) the semigroup and numerical-range
bounds are all infinite, while the deflated Kreiss constant follows the propagator within $13\,\%$ for
every coupling strength in all computed cases; in prototype (b) the amplification equals $\Mst$ exactly,
although no semigroup exists; and in prototype (c) the classical amplification ratio of linear
compartment chains grows like the square root of their length while the fractional one stays small,
and for low-rank generators learned from data the slow-component quantities are computed at a cost set
by the rank of the estimator, while a spurious eigenvalue in $\Sa$ is flagged by an infinite deflated
Kreiss constant. The numerical methods, the complete data and a statistical study of random
finite-rank operators are collected in a separate Supplementary Material, since none of the results of
the paper depends on them.

\section{Preliminary results}\label{sec:prelim}

\subsection{Notation}

All powers are principal, $\omega^\beta:=e^{\beta\Log\omega}$ with $\Arg\in(-\pi,\pi]$, and
$0^\beta:=0$ for $\beta>0$. We write $\HH:=\{\omega\in\C:\RE\omega>0\}$ for the open right half-plane
and, for $0<\theta\le\pi$, $\Sigma(\theta):=\{z\in\C\setminus\{0\}:|\Arg z|<\theta\}$ for the open
sector of half-opening $\theta$. Two sets appear throughout,
\[
\Sa:=\Sigma(\tfrac{\alpha\pi}{2})\qquad\text{and}\qquad\La:=\C\setminus\Sa ,
\]
namely the open sector of instability and the closed fractional stability region; $\La$ contains $0$, the
left half-plane and the two boundary rays $\{re^{\pm i\frac{\alpha\pi}{2}}:r\ge0\}$. For $0<\alpha<2$ the
map $\omega\mapsto\omega^\alpha$ is a biholomorphism of $\HH$ onto $\Sa$, and the weight
\[
a_\alpha(z):=|z|\cos\Big(\frac{\Arg z}{\alpha}\Big)\quad(z\in\Sa)\qquad\text{satisfies}\qquad
a_\alpha(\omega^\alpha)=\RE(\omega)\,|\omega|^{\alpha-1}\quad(\omega\in\HH).
\]
It is positive on $\Sa$, extends continuously by $0$ to $\partial\Sa$, and coincides with $\RE z$ when
$\alpha=1$. Throughout, $X$ is a complex Banach space with dual $X^*$ and pairing
$\langle x,x^*\rangle:=x^*(x)$, $H$ is a complex Hilbert space, $\B(X)$ is the algebra of bounded
operators, and $L$ denotes a closed, densely defined operator with domain $D(L)$, resolvent set
$\rho(L)$ and spectrum $\sigma(L)$; bounded operators are included. An operator $L$ on $H$ is
dissipative if $\RE\langle Lv,v\rangle\le0$ for $v\in D(L)$, and m-dissipative if moreover
$(0,\infty)\subset\rho(L)$. A table of notation is given in \cref{app:notation}.

\subsection{The Mittag--Leffler function}

For $\alpha>0$ the Mittag--Leffler function $E_\alpha(z)=\sum_{k\ge0}\frac{z^k}{\Gamma(1+\alpha k)}$ is entire
and real on the real axis, so $E_\alpha(\bar z)=\overline{E_\alpha(z)}$; moreover 
$E_1(z)=e^z$ and
$E_{\frac{1}{2}}(z)=e^{z^2}\operatorname{erfc}(-z)$ \citep{Gorenflo2020}. All later estimates rest on its behaviour at infinity: it grows like $\alpha^{-1}\exp(z^{\frac{1}{\alpha}})$ inside $\Sa$, where
$z^{\frac{1}{\alpha}}$ has positive real part, it oscillates with modulus close to $\frac{1}{\alpha}$ on the boundary
rays, where $z^{\frac{1}{\alpha}}$ is purely imaginary, and it decays algebraically in the interior of $\La$.

For $0<\alpha<1$ the function $x\mapsto E_\alpha(-x)$ is moreover the Laplace transform of a probability
distribution, so that the fractional dynamics is an average of exponential ones. Let $\Theta_\alpha$ be a
positive random variable with $\mathbb E e^{-\lambda\Theta_\alpha}=e^{-\lambda^\alpha}$ (the one-sided $\alpha$-stable
law), $g_\alpha$ its density, and $m_\alpha(u)=\alpha^{-1}u^{-1-\frac{1}{\alpha}}g_\alpha(u^{-\frac{1}{\alpha}})$ the density of
$Y_\alpha:=\Theta_\alpha^{-\alpha}$, the Mittag--Leffler distribution. 



Since the supremum of $|E_\alpha|$ over the stability region is the constant that governs normal
operators, recall \cref{thm:D},  we record its basic properties here. Let
\[
\Mst:=\sup_{z\in\La}|E_\alpha(z)|,\qquad M_\alpha(\theta):=\sup_{r\ge0}|E_\alpha(re^{i\theta})|\quad
(\tfrac{\alpha\pi}{2}\le\theta\le\pi).
\]
For $\alpha=1$ we have $M_1(\theta)=1$ for $\frac{\pi}{2}\le\theta\le\pi$ and $M^*_1=1$.

The next lemma collects the facts on $E_\alpha$ that we use, like the fact that $|E_\alpha|$ is largest on the boundary rays of
$\Sigma_\alpha$, where it exceeds $1/\alpha$; this gives the sharp constant for
normal and numerical-range generators in Section~\ref{sec:hilbert} and the transient
growth of normal operators.
\begin{lemma}\label{lem:mstar}
Let $0<\alpha<1$. Then $\Mst=M_\alpha(\frac{\alpha\pi}{2})<\infty$, the function $\theta\mapsto M_\alpha(\theta)$ is
nonincreasing, $M_\alpha(\theta)=1$ for $\frac{\pi}{2}\le\theta\le\pi$,
$\lim_{\theta\downarrow\frac{\alpha\pi}{2}}M_\alpha(\theta)=\Mst$, and $\Mst>\frac{1}{\alpha}$. More generally, for every
closed sector $\overline{\Sigma}$ with vertex $0$ and opening less than $2\pi$ that does not meet $\Sa$,
the supremum of $|E_\alpha|$ over $\overline\Sigma$ is the larger of the values $M_\alpha(|\theta|)$ over
the arguments $\theta\in(-\pi,\pi]$ of its two boundary rays.
\end{lemma}

\begin{proof}
By \cref{lem:ML}(c), $E_\alpha$ is bounded on $\La$. A closed sector of opening less than
$2\pi$ that does not meet $\Sa$ is contained in $\La$, so $E_\alpha$ is bounded and continuous on it and
holomorphic in its interior. By the Phragm\'en--Lindel\"of principle for bounded holomorphic
functions in sectors \citep[Chapter~VI]{Conway1978}, its supremum equals the supremum over the two
boundary rays, which is the last claim; by conjugation symmetry the value on the ray of argument
$\theta$ is $M_\alpha(|\theta|)$. Applied to the sectors $\{\theta\le|\Arg z|\le\pi\}\cup\{0\}$,
$\frac{\alpha\pi}{2}\le\theta<\pi$, it shows that $M_\alpha$ is nonincreasing and that $\Mst=M_\alpha(\frac{\alpha\pi}{2})$.
For $\RE z\le0$ the integral $\int_0^\infty m_\alpha(u)e^{zu}\,du$ converges absolutely, it is holomorphic
in $\{\RE z<0\}$ and continuous on $\{\RE z\le0\}$, and it equals $E_\alpha$ on $(-\infty,0]$ by
\cref{lem:ML}(d); hence it equals $E_\alpha$ on $\{\RE z\le0\}$, so $|E_\alpha(z)|\le1=E_\alpha(0)$
there and $M_\alpha(\theta)=1$ for $\theta\ge\frac{\pi}{2}$. As a supremum of continuous functions, $M_\alpha$ is
lower semicontinuous, which together with monotonicity gives the limit. Finally, on the ray
$z=re^{i\frac{\alpha\pi}{2}}$ we have $z^{\frac{1}{\alpha}}=ir^{\frac{1}{\alpha}}$, and \cref{lem:ML}(a) gives
$E_\alpha(z)=\alpha^{-1}e^{ir^{\frac{1}{\alpha}}}-\frac{r^{-1}e^{-i\frac{\alpha\pi}{2}}}{\Gamma(1-\alpha)}+O(r^{-2})$. Along a
sequence $r_k\to\infty$ with $e^{ir_k^{\frac{1}{\alpha}}}=-e^{-i\frac{\alpha\pi}{2}}$ this yields
$|E_\alpha(z)|=\alpha^{-1}+\frac{r_k^{-1}}{\Gamma(1-\alpha)}+O(r_k^{-2})$, which exceeds $\alpha^{-1}$ for large $k$.
\end{proof}

Numerically, $\alpha\Mst=1.326$, $1.245$, $1.171$, $1.068$ and $1.034$ for $\alpha=0.1$, $0.3$, $0.5$,
$0.8$ and $0.9$; the maximum is attained at moderate arguments, $r_*^{\frac{1}{\alpha}}\approx2$ to $2.6$, and
not in the asymptotic regime (\cref{tab:constants}).

\subsection{Mittag--Leffler families}\label{sec:families}

We write $k_\beta(t):=\frac{t^{\beta-1}}{\Gamma(\beta)}$ for $\beta>0$ and $t>0$. The natural notion of a solution
operator for \eqref{eq:caputo} is defined through the Laplace transform, which for bounded $L$ is
computed termwise from the power series and for unbounded $L$ becomes the definition.

\begin{definition}\label{def:family}
Let $0<\alpha\le1$. The operator $L$ \emph{generates a Mittag--Leffler family} if there are a
strongly continuous function $S:[0,\infty)\to\B(X)$ and $w\ge0$ such that
$\int_0^\infty e^{-wt}\|S(t)\|\,dt<\infty$ and, for every real $\lambda>w$, $\lambda^\alpha\in\rho(L)$ and
\begin{equation}\label{eq:laplace}
  \int_0^\infty e^{-\lambda t}S(t)x\,dt=\lambda^{\alpha-1}(\lambda^\alpha-L)^{-1}x
  \qquad (x\in X).
\end{equation}By the uniqueness theorem for the Laplace transform, $S$ is determined by $L$, and we write
$E_\alpha(t^\alpha L):=S(t)$.
\end{definition}

For $\alpha=1$ these are the generators of $C_0$-semigroups and $E_1(tL)=e^{tL}$. For $0<\alpha<1$,
$E_\alpha(t^\alpha L)$ is the resolvent family of the Volterra equation with kernel $k_\alpha$, also called
the $\alpha$-times resolvent family or solution operator of $L$ 
\citep{GhavidelShahali2023}. If
$x\in D(L)$, then
$\lambda^{\alpha-1}(\lambda^\alpha-L)^{-1}x=\lambda^{-1}x+\lambda^{-\alpha}\cdot\lambda^{\alpha-1}(\lambda^\alpha-L)^{-1}Lx$,
and uniqueness of the Laplace transform turns \eqref{eq:laplace} into
\begin{equation}\label{eq:integrated}
E_\alpha(t^\alpha L)x=x+\int_0^tk_\alpha(t-s)\,E_\alpha(s^\alpha L)Lx\,ds\qquad(t\ge0,\ x\in D(L)),
\end{equation}
the integrated form of \eqref{eq:caputo}; in particular $E_\alpha(0\cdot L)=I$ by density.

In the following two lemmas, whose proofs we provide in \cref{app:aux} for completeness, we summarize main facts on Mittag--Leffler families generated by bounded and sectorial operators. 

\begin{lemma}\label{lem:families}
Let $0<\alpha\le1$.
\begin{enumerate}[label=\textup{(\alph*)}]
\item If $L\in\B(X)$, then $L$ generates a Mittag--Leffler family,
$E_\alpha(t^\alpha L)=\sum_{k\ge0}\frac{t^{\alpha k}L^k}{\Gamma(1+\alpha k)}$, the map
$t\mapsto E_\alpha(t^\alpha L)$ is norm continuous, $\|E_\alpha(t^\alpha L)\|\le E_\alpha(t^\alpha\|L\|)$, and
\eqref{eq:laplace} holds for every $\omega\in\C$ with $\RE\omega>\|L\|^{\frac{1}{\alpha}}$ in place of $\lambda$,
the integral converging absolutely in norm.
\item \textup{(Subordination.)} Let $0<\alpha<1$, and let $L$ generate a $C_0$-semigroup with
$\|e^{sL}\|\le Me^{\varpi s}$ for some $M\ge1$ and $\varpi\ge0$. Then $L$ generates a Mittag--Leffler family,
\[
E_\alpha(t^\alpha L)x=\int_0^\infty m_\alpha(u)\,e^{t^\alpha uL}x\,du\qquad(t\ge0,\ x\in X),
\]
and $\sup_{t\ge0}\|E_\alpha(t^\alpha L)\|\le\sup_{s\ge0}\|e^{sL}\|$.
\item If $X=H$ and $L$ is m-dissipative, then $\|E_\alpha(t^\alpha L)\|\le1$ for all $t\ge0$.
\end{enumerate}
\end{lemma}

Part (b) is the subordination principle \citep{Bazhlekova2001}: $E_\alpha(t^\alpha L)$ averages the semigroup
over the random time $t^\alpha Y_\alpha$. It yields no bound when $\sigma(L)$ meets the right half-plane,
and it does not apply when $L$ generates no $C_0$-semigroup. Such operators are the main object of
this paper, and the next lemma covers the most important class of them. 

\begin{lemma}[Sectorial generators]\label{lem:sectorial}
Let $0<\alpha\le1$ and $\frac{\alpha\pi}{2}<\psi\le\pi$. Suppose that $\Sigma(\psi)\subset\rho(L)$ and
\[
s_{\psi'}(L):=\sup_{z\in\Sigma(\psi')}\|z(z-L)^{-1}\|<\infty\qquad\text{for every }\psi'\in(\tfrac{\alpha\pi}{2},\psi).
\]
Fix $\psi'\in(\frac{\alpha\pi}{2},\min\{\psi,\alpha\pi\})$, put $\varphi:=\frac{\psi'}{\alpha}\in(\frac{\pi}{2},\pi)$, and let
$\mathrm{Ha}_\varphi$ be the Hankel contour consisting of the ray $\{\varrho e^{-i\varphi}:\varrho\ge1\}$
traversed towards $e^{-i\varphi}$, the arc $\{e^{i\vartheta}:|\vartheta|\le\varphi\}$ traversed
counterclockwise, and the ray $\{\varrho e^{i\varphi}:\varrho\ge1\}$ traversed away from $e^{i\varphi}$. Then
$L$ generates a Mittag--Leffler family,
\begin{equation}\label{eq:hankel}
E_\alpha(t^\alpha L)=\frac1{2\pi i}\int_{\mathrm{Ha}_\varphi}e^{\omega t}\,\omega^{\alpha-1}(\omega^\alpha-L)^{-1}\,d\omega\qquad(t>0),
\end{equation}
with absolute convergence in norm, and
\[
\sup_{t\ge0}\|E_\alpha(t^\alpha L)\|\le c_\varphi\,s_{\psi'}(L),\qquad
c_\varphi:=\frac1{2\pi}\int_{\mathrm{Ha}_\varphi}|e^\omega|\,|\omega|^{-1}\,|d\omega| .
\]
\end{lemma}

For $\alpha=1$ this is the classical construction of bounded analytic semigroups. For $\alpha<1$ and
$\psi\le\frac{\pi}{2}$ the operator $L$ may have spectrum in the right half-plane, and then it need not generate
a $C_0$-semigroup. 

\subsection{The novel fractional Kreiss constant}\label{sec:kreiss-constant}

\begin{definition}\label{def:kreiss}
For a closed operator $L$ on $X$ with $\Sa\subset\rho(L)$, define the  \emph{fractional Kreiss constant} to be
\[
\Kr_\alpha(L):=\sup_{\omega\in\HH}\RE(\omega)\,\big\|\omega^{\alpha-1}(\omega^\alpha-L)^{-1}\big\|
=\sup_{z\in\Sa}a_\alpha(z)\,\|(z-L)^{-1}\|.
\]
If $\Sa\not\subset\rho(L)$ we put $\Kr_\alpha(L):=\infty$.
\end{definition}

The two expressions coincide by the substitution $z=\omega^\alpha$ and the identity for
$a_\alpha(\omega^\alpha)$. The first one shows that $\Kr_\alpha(L)$ is the classical Kreiss quantity
$\sup_{\RE\omega>0}\RE(\omega)\|\widehat S(\omega)\|$ of the Laplace transform $\widehat S$ of the family, and the
second one shows how the resolvent of $L$ enters: near the boundary rays of $\Sa$ the weight vanishes
linearly, so eigenvalues on these rays are allowed as long as they are semisimple. First, we will prove some of the properties of the fractional Kreiss constant.

\begin{lemma}\label{lem:elementary}
Let $X\ne\{0\}$ and let $L$ be closed and densely defined.
\begin{enumerate}[label=\textup{(\alph*)}]
\item $\Kr_\alpha(L)\ge1$.
\item $\Kr_\alpha(L)\le\Kr_1(L)$, and the inequality can be strict with $\Kr_1(L)=\infty$.
\item For $z\in\Sa$, $\dist(z,\partial\Sa)\le a_\alpha(z)\le\alpha^{-1}\dist(z,\partial\Sa)$.
\item \textup{(Pseudospectral form.)} Let $\sigma_\eps(L):=\{z:\|(z-L)^{-1}\|>\eps^{-1}\}$, with
$\|(z-L)^{-1}\|:=\infty$ on $\sigma(L)$, and
$a_{\alpha,\eps}(L):=\sup\{a_\alpha(z):z\in\sigma_\eps(L)\cap\Sa\}$, with $\sup\emptyset:=0$. Then
$\Kr_\alpha(L)=\sup_{\eps>0}\frac{a_{\alpha,\eps}(L)}{\eps}$.
\item \textup{(Scaling.)} $\Kr_\alpha(cL)=\Kr_\alpha(L)$ for every $c>0$.
\end{enumerate}
\end{lemma}

\begin{proof}
(a) We may assume $\Kr:=\Kr_\alpha(L)<\infty$. Since $a_\alpha(\lambda)=\lambda$ for $\lambda>0$, we have
$\|\lambda(\lambda-L)^{-1}\|\le\Kr$. For $x\in D(L)$,
$\|\lambda(\lambda-L)^{-1}x-x\|=\|(\lambda-L)^{-1}Lx\|\le\frac{\Kr\|Lx\|}{\lambda}\to0$ as $\lambda\to\infty$. By density and the
uniform bound, $\lambda(\lambda-L)^{-1}x\to x$ for every $x$, whence
$\Kr\ge\liminf_{\lambda\to\infty}\|\lambda(\lambda-L)^{-1}\|\ge1$.

(b) If $z=|z|e^{i\vartheta}$ with $|\vartheta|<\frac{\alpha\pi}{2}$, then $\frac{|\vartheta|}{\alpha}\ge|\vartheta|$, so
$a_\alpha(z)\le\RE z$, and $\Sa\subset\HH$. For $L=e^{i\theta}$ on $X=\C$ with $\frac{\alpha\pi}{2}<\theta<\frac{\pi}{2}$ we
have $\Kr_1(L)=\infty$, while $|z-e^{i\theta}|\ge\dist(z,\partial\Sa)$ for $z\in\Sa$, so $\Kr_\alpha(L)\le\frac{1}{\alpha}$
by (c).

(c) Let $\delta:=\frac{\alpha\pi}{2}-|\Arg z|\in(0,\frac{\alpha\pi}{2}]$. Then $\dist(z,\partial\Sa)=|z|\sin\delta$ and
$a_\alpha(z)=|z|\sin(\frac{\delta}{\alpha})$. Since $\delta\le\frac{\delta}{\alpha}\le\frac{\pi}{2}$, we have
$\sin\delta\le\sin(\frac{\delta}{\alpha})$, and the function $g(\delta):=\alpha^{-1}\sin\delta-\sin(\frac{\delta}{\alpha})$ satisfies
$g(0)=0$ and $g'(\delta)=\alpha^{-1}(\cos\delta-\cos(\frac{\delta}{\alpha}))\ge0$.

(d) If $z\in\sigma_\eps(L)\cap\Sa$, then $\frac{a_\alpha(z)}{\eps}<a_\alpha(z)\|(z-L)^{-1}\|$. Conversely, if
$\|(z-L)^{-1}\|=\eps_0^{-1}$ with $z\in\Sa$, then $z\in\sigma_\eps(L)$ for every $\eps>\eps_0$, and
$\frac{a_\alpha(z)}{\eps}\to a_\alpha(z)\|(z-L)^{-1}\|$ as $\eps\downarrow\eps_0$. If $z_0\in\sigma(L)\cap\Sa$, then both
sides are infinite.

(e) The map $z\mapsto cz$ preserves $\Sa$, $a_\alpha(cz)=c\,a_\alpha(z)$ and
$\|(cz-cL)^{-1}\|=c^{-1}\|(z-L)^{-1}\|$.
\end{proof}

Part (d) relates $\Kr_\alpha$ to the fractional pseudospectra of \citet{SancaKosticCvetkovic2018}: it
measures how far the $\eps$-pseudospectrum reaches into $\Sa$, in the units of $a_\alpha$, relative to
$\eps$. Part (e) shows that $\Kr_\alpha$ is a dimensionless quantity, which is why a uniform bound on
$\Kr_\alpha$ carries no information about the size of $L$; this observation will be used in
\cref{sec:optimal}.

\begin{remark}[Other normalisations]\label{rem:Kplus}
Since $\RE\omega\le|\omega|$, we have that the quantity
$
\Kr^+_\alpha(L):=\sup_{\omega\in\HH}(\RE\omega)^\alpha\|(\omega^\alpha-L)^{-1}\|
$ 
dominates $\Kr_\alpha(L)$, but it is not
equivalent to boundedness of $E_\alpha(t^\alpha L)$. For $0<\alpha<1$ and the scalar operator
$L=e^{i\frac{\alpha\pi}{2}}$, $\sup_t|E_\alpha(t^\alpha L)|=\Mst$ by \cref{lem:mstar} and $\Kr_\alpha(L)\le\frac{1}{\alpha}$ by
\cref{lem:elementary}(c), while the points $z=e^{i(\frac{\alpha\pi}{2}-\delta)}=\omega^\alpha$,
$\omega=e^{i(\frac{\pi}{2}-\frac{\delta}{\alpha})}$, give $\Kr^+_\alpha(L)\ge\frac{\sin^\alpha(\frac{\delta}{\alpha})}{\delta}\to\infty$ as
$\delta\downarrow0$. For the weight $\RE\omega$ of \cref{def:kreiss}, by contrast, finiteness is necessary for
boundedness of the family (\cref{thm:kreiss-lower}).
\end{remark}

\subsection{Slow components}\label{sec:slow-def}

We say that a projection $P\in\B(X)$ \emph{commutes with} $L$ if $PD(L)\subset D(L)$ and $PLx=LPx$ for
$x\in D(L)$. The part of $L$ in $PX$ is the restriction of $L$ to $D(L)\cap PX$; it is closed and densely
defined in $PX$, because $PD(L)\subset D(L)\cap PX$ and $Lx=LPx=PLx$ for $x\in D(L)\cap PX$.

The idea behind the following definition is that the transient growth of many generators is produced
by finitely many modes, while the remaining dynamics is bounded for a structural reason. The word slow
refers to the role of these modes in the transient rather than to the size of their eigenvalues: a slow
component collects the finitely many modes that can produce growth, and its complement carries a bounded
family.

\begin{definition}\label{def:slow}
Let $L$ be closed and densely defined on $X$. A \emph{slow component of rank $m\ge1$} is a projection
$P\in\B(X)$ of rank $m$ that commutes with $L$ and such that the part $L_Q$ of $L$ in $QX$, $Q:=I-P$,
generates a Mittag--Leffler family on $QX$ with $\sup_{t\ge0}\|E_\alpha(t^\alpha L_Q)\|<\infty$. We call
$L_P$, the part of $L$ in $PX$, the \emph{slow block}, and
\[
\Kr_\alpha(L;P):=\sup_{z\in\Sa}a_\alpha(z)\,\|(z-L_P)^{-1}P\|
\]
the \emph{deflated Kreiss constant}, with $\Kr_\alpha(L;P):=\infty$ if $\sigma(L_P)\cap\Sa\ne\emptyset$; the
norm is that of $\B(X)$. By \cref{lem:slow-basic}(a) below, $(z-L_P)^{-1}P=(z-L)^{-1}P$ for $z\in\rho(L)$.
\end{definition}

The typical candidate for a slow component is the Riesz projection
$P=(2\pi i)^{-1}\oint_\Gamma(z-L)^{-1}dz$ onto finitely many isolated eigenvalues of finite algebraic
multiplicity, enclosed by a positively oriented contour $\Gamma$ that encloses no other point of
$\sigma(L)$; then $P$ has finite rank, $PX\subset D(L)$, $P$ commutes with $L$, $\sigma(L_P)$ consists of the
enclosed eigenvalues and $\sigma(L_Q)$ is the rest of $\sigma(L)$ \citep[Theorem~III.6.17]{Kato1995}. It is a slow
component as soon as $L_Q$ generates a bounded family, for instance by \cref{lem:sectorial} or, on Hilbert
spaces, by the numerical-range criterion of \cref{sec:hilbert}.

\begin{lemma}\label{lem:slow-basic}
Let $P$ be a projection of finite rank $m\ge1$ commuting with $L$, and let $Q:=I-P$.
\begin{enumerate}[label=\textup{(\alph*)}]
\item $PX\subset D(L)$, $L_P\in\B(PX)$, and $\sigma(L)=\sigma(L_P)\cup\sigma(L_Q)$; for $z\in\rho(L)$,
$(z-L)^{-1}P=(z-L_P)^{-1}P$ and $(z-L)^{-1}Q=(z-L_Q)^{-1}Q$.
\item If $P=UV$ with $U:\C^m\to X$ injective, $V:X\to\C^m$ surjective and $VU=I_m$, and if
$\Lambda\in\C^{m\times m}$ is the matrix of $L_P$ in this factorisation, then
$(z-L_P)^{-1}P=U(z-\Lambda)^{-1}V$ for $z\in\rho(L_P)$.
\item $\max\{\|P\|,\Kr_\alpha(L_P)\}\le\Kr_\alpha(L;P)\le\|P\|\,\Kr_\alpha(L_P)$, and if $\Sa\subset\rho(L)$, then
$\Kr_\alpha(L_P)\le\Kr_\alpha(L)$. If $X=H$ and $P$ is an orthogonal projection, then
$\Kr_\alpha(L;P)=\Kr_\alpha(L_P)$.
\end{enumerate}
\end{lemma}

\begin{proof}
(a) Since $D(L)$ is dense, $PD(L)$ is a dense subspace of the finite-dimensional space $PX$, hence
$PX=PD(L)\subset D(L)$, and $Ly=LPy=PLy\in PX$ for $y\in PX$. If $z\in\rho(L_P)\cap\rho(L_Q)$, then
$R:=(z-L_P)^{-1}P+(z-L_Q)^{-1}Q$ maps $X$ into $D(L)$ and satisfies $(z-L)R=I$ and $R(z-L)x=x$ for $x\in D(L)$,
so $z\in\rho(L)$. Conversely, applying $(z-L)^{-1}$ to $P(z-L)x=(z-L)Px$ shows that $(z-L)^{-1}$ commutes with
$P$ and $Q$, so its restrictions invert $z-L_P$ and $z-L_Q$.

(b) For $y\in PX$ we have $L_Py=U\Lambda Vy$, and $U(z-\Lambda)^{-1}V$ inverts $z-L_P$ on $PX$.

(c) The lower bound $\|P\|$ holds because $\lambda(\lambda-L_P)^{-1}P\to P$ as $\lambda\to\infty$ and
$a_\alpha(\lambda)=\lambda$; the lower bound $\Kr_\alpha(L_P)$ because $y=Py$ for $y\in PX$; the upper bound because
$\|(z-L_P)^{-1}P\|\le\|(z-L_P)^{-1}\|\,\|P\|$; and $\Kr_\alpha(L_P)\le\Kr_\alpha(L)$ by (a). In the orthogonal case
$\|(z-L_P)^{-1}P\|=\|(z-L_P)^{-1}\|$.
\end{proof}

By (b), $\Kr_\alpha(L;P)$ is a Kreiss quantity of an operator of rank at most $m$ that involves only
the slow block and the two factors of the projection, so it can be evaluated at a cost governed by
$m$ and by the cost of applying $U$ and $V$. By (c) it lies between
the Kreiss constant of the slow block and the same constant multiplied by the obliqueness $\|P\|$ of
the splitting; in prototype (a), analysed in \cref{sec:slow,sec:experiments}, the propagator follows
$\Kr_\alpha(L;P)$ closely although $\|P\|$ is large.

\section{A Bernstein-type inequality on vertical lines}\label{sec:bernstein}


Before stating the main inequality that is the foundation of the theorem for the fractional Kreiss case we explain where it comes from.
Let $L$ be bounded with
$\sigma(L)\cap\Sa=\emptyset$, and let $x\in X$, $x^*\in X^*$ be fixed. Since the main aim is to find a suitable upper-bound for the norm of the function $E_\alpha(t^\alpha L)$, it is sensible to consider the scalar function and its Laplace transform
\[
g(t):=\langle E_\alpha(t^\alpha L)x,x^*\rangle, \;\text{ and }
\widehat g(\omega)=\omega^{\alpha-1}\langle(\omega^\alpha-L)^{-1}x,x^*\rangle,
\] 
respectively, which will, by the choice of the fractional Kreiss constant, control
$|\widehat g|$ on vertical lines: 
\[
|\widehat g(a+iy)|\le\frac{\Kr_\alpha(L)\|x\|\|x^*\|}{a}
\]
To recover $g(t)$ we will invert the Laplace
transform along the line $\RE\omega=a$. The direct inversion integral does not converge absolutely,
because $\widehat g$ decays only like $|\omega|^{-1}$, but after one integration by parts it does
(see \eqref{eq:inversion} in \cref{sec:finite-rank}), and it gives
\[
|g(t)|\le\frac{e^{at}}{2\pi t}\int_\R|\widehat g\,'(a+iy)|\,dy .
\]
The integral on the right is the length of the curve $y\mapsto\widehat g(a+iy)$, while the Kreiss constant
controls its distance from the origin. The main task is to compare the two, that is, to show
that a curve of this special form cannot be much longer than its maximal modulus. If $L$ has rank $n$,
the function $z\mapsto z\langle(z-L)^{-1}x,x^*\rangle$ is rational of degree at most $n$
(\cref{lem:struct}), so $\widehat g=h_F$ with
\[
h_F(\omega):=\omega^{-1}F(\omega^\alpha)\qquad(\omega\in\HH)
\]
and $F$ in the class $\RR_{\alpha,n}$ of nonzero rational functions $F=\frac{F_1}{F_0}$, with coprime
polynomials $F_1,F_0$, such that
\[
\deg F_1\le\deg F_0\le n\qquad\text{and}\qquad F_0(z)\ne0\ \text{ for all } z\in\Sa .
\]
Zeros of $F_0$ at $0$ or on $\partial\Sa$ are allowed. For $F\in\RR_{\alpha,n}$ and $a>0$ we consider the
quotient of length and maximal modulus,
\[
\Jq_F(a):=\frac{\int_\R|h_F'(a+iy)|\,dy}{\sup_{y\in\R}|h_F(a+iy)|}.
\]
The main aim will be to bound this quotient and to do so, we will prove the following Bernstein-type theorem : 
\begin{theorem}[Bernstein-type inequality]\label{thm:bernstein}
Let $0<\alpha\le1$, $n\in\N_0$ and $F\in\RR_{\alpha,n}$. For every $a>0$ the quotient $\Jq_F(a)$ is well
defined and
\[
2\le\Jq_F(a)\le\pi\big(3+4n(2-\alpha)\big).
\]
If $F$ is regular at $0$ with $F(0)=0$, then $\Jq_F(a)\le\pi\big(1+4n(2-\alpha)\big)$.
\end{theorem}
We call it Bernstein-type because it bounds an $L^1$ norm of the derivative by the supremum of the
function, with a constant that depends only on a degree. For $\alpha=1$, $h_F$ is rational of degree at
most $n+1$, and Spijker's lemma \citep{Spijker1991,wegert1994buffon}, transported from the circle to the
line by a M\"obius map, which preserves the degree, gives the better constant $2\pi(n+1)$. That proof counts the intersections of
the curve with straight lines by B\'ezout's theorem. However for $\alpha<1$ the curve $\{(a+iy)^\alpha:y\in\R\}$
is not algebraic when $\alpha$ is irrational, and when $\alpha=\frac{p}{q}$ is rational its degree grows with
$q$, so such count will not be uniform in $\alpha$. The proof that we will present avoids intersection counts altogether: it
bounds the total turning of the tangent of the curve $y\mapsto h_F(a+iy)$, and this bound is uniform in
$\alpha\in(0,1]$.

Our proof technique has two ingredients. The first is classical: the length of a plane curve is the average, over
all directions, of the total variation of its projections onto lines, and a projection changes monotonicity
only when the tangent of the curve is orthogonal to the direction of projection, so the length is controlled
by the maximal modulus times the total turning of the tangent. The second ingredient, which is where the fractional order enters, is a geometric bound on that turning, uniform in $\alpha$, which we discuss in detail first.


Throughout, for a continuous function $u$ on an open interval $I\subseteq\R$ we write:
\begin{itemize}
  \item $\Var_I(u)$ for its total variation;
  \item $N_u(p):=\#\{y\in I:u(y)=p\}\in\N_0\cup\{\infty\}$ for its counting function;
  \item $N_u^{[c]}(\vartheta):=\#\{y\in I:u(y)\in\vartheta+c\Z\}$, for $c>0$, for its counting function modulo $c$.
\end{itemize} There is a link between the stated quantities which comes from the Banach's indicatrix theorem
\citep{Banach1925}, extended from compact intervals to $I$ by monotone convergence, and by Tonelli's theorem applied to
$N_u^{[c]}(\vartheta)=\sum_{k\in\Z}N_u(\vartheta+kc)$,
\begin{equation}\label{eq:banach}
\Var_I(u)=\int_\R N_u(p)\,dp=\int_0^cN^{[c]}_u(\vartheta)\,d\vartheta\qquad(c>0),
\end{equation}
and $\Var_I(u)=\int_I|u'|$ if $u\in C^1(I)$. We also use the elementary fact that a $C^1$ function whose derivative has
$k<\infty$ zeros takes every value at most $k+1$ times, since it is strictly monotone between consecutive critical points.

For $\alpha=1$ the image $\{1+iy:y\in\R\}$ of a vertical line is a straight line, and the argument of $1+iy-\zeta$ seen
+from any point $\zeta$ off the line changes by exactly $\pi$. For $\alpha<1$ the image
\[
\mathcal C_\alpha:=\{(1+iy)^\alpha:y\in\R\}\subset\Sa
\]
bends, and the following lemma shows that a curve whose tangent turns monotonically through an angle $\Theta$ is seen from
any external point under an angle of variation at most $\pi+\Theta$. This geometric fact replaces the algebraic
intersection counts of the case $\alpha=1$.

\begin{lemma}[Turning lemma]\label{lem:turning}
Let $\gamma\in C^1(I;\C)$ with $\gamma'\ne0$ on an open interval $I$, and suppose that a continuous lift $T$ of
$\arg\gamma'$ is either constant or strictly monotone, with range contained in an open interval $J$ of
length $\Theta\le\pi$; when $T$ is constant, $J$ is replaced by its single value and $\Theta:=0$. Let $\zeta\in\C\setminus\gamma(I)$ and let $\beta_\zeta$ be a
continuous lift of $s\mapsto\arg(\gamma(s)-\zeta)$. Then $\Var_I(\beta_\zeta)\le\pi+\Theta$.
\end{lemma}

\begin{proof}
For $\tau\in[0,\pi)$ let $\ell_\tau:=\zeta+\R e^{i\tau}$ be the line through $\zeta$ in direction $e^{i\tau}$, and
consider $f_\tau(s):=\IM\big(e^{-i\tau}(\gamma(s)-\zeta)\big)$, whose zeros are the parameters with
$\gamma(s)\in\ell_\tau$. Since $f_\tau'(s)=|\gamma'(s)|\sin(T(s)-\tau)$, the critical points of $f_\tau$ are the
parameters with $T(s)\in\tau+\pi\Z$. If $T$ is strictly monotone, there is at most one such parameter for each
element of $(\tau+\pi\Z)\cap J$, and this set has at most one element because $J$ is open of length at most $\pi$;
if $T$ is constant, $\gamma$ parametrises a segment, a half-line or a line, which meets $\ell_\tau$ at most once
except for the single direction parallel to it. Hence, for almost every $\tau$, $\gamma$ meets $\ell_\tau$ at most twice,
and at most once if $(\tau+\pi\Z)\cap J=\emptyset$; the set of $\tau\in[0,\pi)$ with $(\tau+\pi\Z)\cap J\ne\emptyset$ has
measure at most $\Theta$. Now $\arg(\gamma(s)-\zeta)\in\tau+2\pi\Z$ exactly when $\gamma(s)$ lies on the open half-line
$\zeta+(0,\infty)e^{i\tau}$. If $K(\tau)$ denotes the number of such parameters, then $K(\tau)+K(\tau+\pi)$ is the number of
intersections of $\gamma$ with $\ell_\tau$, because $\gamma(s)\ne\zeta$. By \eqref{eq:banach} with $c=2\pi$,
\[
\Var_I(\beta_\zeta)=\int_0^{2\pi}K(\tau)\,d\tau=\int_0^\pi\big(K(\tau)+K(\tau+\pi)\big)\,d\tau\le2\Theta+(\pi-\Theta)=\pi+\Theta.\qedhere
\]
\end{proof}

We apply the lemma to $\mathcal C_\alpha$. With $y=\tan\varphi$, $|\varphi|<\frac{\pi}{2}$, the curve is parametrised by
$w(\varphi)=(1+i\tan\varphi)^\alpha=\sec^\alpha\!\varphi\,e^{i\alpha\varphi}$, and after some algebra
$w'(\varphi)=\alpha\sec^{\alpha+1}\!\varphi\,e^{iT(\varphi)}$ with $T(\varphi)=\frac\pi2-(1-\alpha)\varphi$. This lift is constant for
$\alpha=1$ and strictly decreasing for $\alpha<1$, with range $(\frac{\alpha\pi}2,\pi-\frac{\alpha\pi}2)$ of length $(1-\alpha)\pi$.
Consequently, for every $\zeta\notin\mathcal C_\alpha$ and every continuous lift $\beta_\zeta$ of $y\mapsto\arg((1+iy)^\alpha-\zeta)$,
\begin{equation}\label{eq:factor}
\Var_\R(\beta_\zeta)\le(2-\alpha)\pi .
\end{equation}
The bound is attained: seen from a point $\zeta>1$ inside the convex region bounded by $\mathcal C_\alpha$, for instance
$\zeta=2$, the argument of $(1+iy)^\alpha-\zeta$ decreases monotonically from $-\frac{\alpha\pi}{2}$ through $\pm\pi$ to $\frac{\alpha\pi}{2}$,
whereas seen from the origin it only sweeps the angle $\alpha\pi$ between the asymptotic directions of the curve.


\begin{proof}[Proof of \Cref{thm:bernstein}]
Fix $F=\frac{F_1}{F_0}\in\RR_{\alpha,n}$, write $q:=\deg F_0\le n$ and $h:=h_F$, and put $D_F(z):=-F(z)+\alpha zF'(z)$. The proof is organized in four components focus on: (i) structure of $h$, (ii) length controlled by turning, (iii) turning of the derivative $h'$ and (iv) lines with critical points.

(i) Since $F$ is holomorphic on $\Sa=\{\omega^\alpha:\omega\in\HH\}$ and
$\omega^{-1}\omega^{\alpha-1}=\omega^{-2}\omega^\alpha$, the product rule gives
\begin{equation}\label{eq:hprime}
h'(\omega)=\omega^{-2}D_F(\omega^\alpha)\qquad(\omega\in\HH).
\end{equation}
Moreover $D_F=\frac{G}{F_0^2}$ with $G:=-F_1F_0+\alpha z(F_1'F_0-F_1F_0')$ of degree at most $2q$, because $\deg F_1\le q$; after
cancelling common factors, $D_F=\frac{G_N}{G_D}$ with $d_N:=\deg G_N\le2q$, $d_D:=\deg G_D\le2q$ and $G_D$ dividing $F_0^2$, so $D_F$ has
no poles in $\Sa$. The function $h$ does not vanish identically on any vertical line, since $F\not\equiv0$ on the open set $\Sa$,
and $D_F\not\equiv0$, since otherwise $h$ would be constant and $h(t)=t^{-1}F(t^\alpha)\to0$ as $t\to+\infty$. As $F$ and $D_F$ are
bounded near infinity, there are $C,y_*>0$ with $|h(a+iy)|\le\frac{C}{|y|}$ and $|h'(a+iy)|\le\frac{C}{|y|^2}$ for all $a>0$ and $|y|\ge y_*$.
Hence the length $\ell_a:=\int_\R|h'(a+iy)|\,dy$ and the maximal modulus $h^*_a:=\sup_y|h(a+iy)|$ are finite and positive, and
$\ell_a\ge2h^*_a$ because the curve $y\mapsto h(a+iy)$ starts and ends at $0$; this is the lower bound of the theorem. Finally,
with $F_a(z):=F(a^\alpha z)\in\RR_{\alpha,n}$ we have $h_F(a\omega)=a^{-1}h_{F_a}(\omega)$, so $\Jq_F(a)=\Jq_{F_a}(1)$ and $F_a(0)=F(0)$, and it
suffices to consider $a=1$.

(ii) Suppose that $h'$ has no zero on $\{\RE\omega=1\}$, and let $\eta$ be a continuous
lift of $y\mapsto\arg h'(1+iy)$. We claim that
\begin{equation}\label{eq:reduction}
\ell_1\le h^*_1\big(\pi+\Var_\R(\eta)\big).
\end{equation}
Put $\gamma(y):=h(1+iy)$ and $u_\vartheta:=\RE(e^{-i\vartheta}\gamma)$ for $\vartheta\in[0,\pi)$. Since
$\int_0^\pi|\cos(\varphi-\vartheta)|\,d\vartheta=2$ for every $\varphi$, Tonelli's theorem gives the Cauchy--Crofton formula
$2\ell_1=\int_0^\pi\Var_\R(u_\vartheta)\,d\vartheta$. A direct computation gives
$u_\vartheta'(y)=-|h'(1+iy)|\sin(\eta(y)-\vartheta)$, so the critical points of $u_\vartheta$ are the solutions of $\eta(y)\in\vartheta+\pi\Z$,
and there are $N^{[\pi]}_\eta(\vartheta)$ of them. Since $|u_\vartheta|\le h^*_1$, \eqref{eq:banach} and the elementary counting fact give
$\Var_\R(u_\vartheta)=\int_{-h^*_1}^{h^*_1}N_{u_\vartheta}(p)\,dp\le2h^*_1\big(N^{[\pi]}_\eta(\vartheta)+1\big)$. Integrating over $\vartheta$ and using
\eqref{eq:banach} with $c=\pi$ gives $2\ell_1\le2h^*_1(\pi+\Var_\R(\eta))$.

(iii) Under the assumption of Step 2, write $D_F(z)=c_0\prod_i(z-a_i)^{k_i}\prod_j(z-b_j)^{-l_j}$ with $c_0\ne0$,
distinct zeros $a_i$ of multiplicities $k_i$ and distinct poles $b_j$ of orders $l_j$, so that $\sum_ik_i=d_N$ and $\sum_jl_j=d_D$.
The poles lie outside $\Sa\supset\mathcal C_\alpha$ by Step 1, and a zero on $\mathcal C_\alpha$ would produce a zero of $h'$ on the line by
\eqref{eq:hprime}; so no factor vanishes on $\mathcal C_\alpha$. By \eqref{eq:hprime},
\[
\eta(y):=\Arg c_0-2\arctan y+\sum_ik_i\beta_{a_i}(y)-\sum_jl_j\beta_{b_j}(y)
\]
is a continuous lift of $\arg h'(1+iy)$, with $\beta_\zeta$ as in \eqref{eq:factor}. Subadditivity of the variation,
$\Var_\R(\arctan)=\pi$ and \eqref{eq:factor} give $\Var_\R(\eta)\le2\pi+(d_N+d_D)(2-\alpha)\pi\le2\pi+4q(2-\alpha)\pi$, and
\eqref{eq:reduction} gives $\Jq_F(1)\le\pi(3+4q(2-\alpha))$.

If $F(0)=0$, then $F_0(0)\ne0$ and $F_1(0)=0$, so $F=\frac{z\tilde F_1}{F_0}$ with $\deg\tilde F_1\le q-1$, and after some algebra
$D_F(z)=z\tilde D(z)$ with $\tilde D=\frac{\tilde G}{F_0^2}$, $\tilde G:=-(1-\alpha)\tilde F_1F_0+\alpha z(\tilde F_1'F_0-\tilde F_1F_0')$ of degree at most
$2q-1$. Then $h'(\omega)=\omega^{\alpha-2}\tilde D(\omega^\alpha)$, the prefactor contributes $(2-\alpha)\pi$ instead of $2\pi$, the factors
contribute at most $(4q-1)(2-\alpha)\pi$, and $\Var_\R(\eta)\le4q(2-\alpha)\pi$, so that $\Jq_F(1)\le\pi(1+4q(2-\alpha))$.

(iv) By \eqref{eq:hprime}, $h'(\omega)=0$ if and only if $G_N(\omega^\alpha)=0$, and since
$\omega\mapsto\omega^\alpha$ is injective on $\HH$, $h'$ has at most $2n$ zeros in $\HH$. By Step 1, $\Jq_F(a)=\Jq_{F_a}(1)$, so the bounds of
Step 3 hold for all but at most $2n$ values of $a$. It remains to see that $a\mapsto\Jq_F(a)$ is continuous. On
$K:=\{a+iy:a\in[a_1,a_2],\,|y|\le y_*\}$, $0<a_1<a_2$, the functions $h$ and $h'$ are jointly continuous and bounded, because the
image of $K$ under $\omega\mapsto\omega^\alpha$ is a compact subset of $\Sa$ and poles of $F$ on $\partial\Sa$ stay at positive distance from
it. Dominated convergence with the majorant $C'\mathbf 1_{\{|y|\le y_*\}}+C|y|^{-2}\mathbf 1_{\{|y|>y_*\}}$ gives continuity of $a\mapsto\ell_a$,
and the uniform tail bound on $h$ gives continuity of $a\mapsto h^*_a>0$. This completes the proof. 
\end{proof}

\begin{remark}[The range $1<\alpha<2$]\label{rem:alpha-large}
For $1<\alpha<2$ the map $\omega\mapsto\omega^\alpha$ is still injective on $\HH$, and the tangent angle
$T(\varphi)=\frac\pi2+(\alpha-1)\varphi$ of $\mathcal C_\alpha$ is increasing with range of length $(\alpha-1)\pi$. The
turning lemma then gives $\Var_\R(\beta_\zeta)\le\alpha\pi$ for each factor, the prefactor $(1+iy)^{-2}$ still contributes
$2\pi$, and the same proof yields $\Jq_F(a)\le\pi(3+4n\alpha)$. We do not pursue this range, because for
$1<\alpha<2$ the family $E_\alpha(t^\alpha L)$ solves \eqref{eq:caputo} only with $u'(0)=0$, and several results below, among
them the subordination formula, do not extend.
\end{remark}

\section{The finite-rank core}\label{sec:finite-rank}
This section proves \Cref{thm:A} for bounded operators of finite rank, which is the heart of the
matter, while later in \cref{sec:slow} we prove the result for generators with a slow component by reduction to this case. The strategy proceeds in
three steps: (i) a lower bound by the Kreiss constant holds for every bounded family; (ii) for a finite-rank operator,
every matrix element of the resolvent is a rational function of degree at most the rank, and the
operator splits into a slow part acting through a small matrix and a stationary kernel; and finally (iii) the
inversion of the Laplace transform along a vertical line, after one integration by parts, turns the
upper bound into the Bernstein-type inequality proven in \cref{sec:bernstein}.

\subsection{The lower bound}

The lower bound in \Cref{thm:A} holds for every generator of a bounded family and requires no
finiteness assumption. Its proof shows that the Laplace transform of a bounded family is the
resolvent, transported to $\HH$, on the whole half-plane and not only on the real half-line where
\cref{def:family} prescribes it.

\begin{theorem}\label{thm:kreiss-lower}
Let $0<\alpha\le1$, let $X\ne\{0\}$, and let $L$ be closed, densely defined and generate a Mittag--Leffler
family with $C_L:=\sup_{t\ge0}\|E_\alpha(t^\alpha L)\|<\infty$. Then $\Sa\subset\rho(L)$,
\[
\omega^{\alpha-1}(\omega^\alpha-L)^{-1}x=\int_0^\infty e^{-\omega t}E_\alpha(t^\alpha L)x\,dt\qquad(\omega\in\HH,\ x\in X),
\]
and $\Kr_\alpha(L)\le C_L$.
\end{theorem}

\begin{proof}
Write $S(t):=E_\alpha(t^\alpha L)$. For $\omega\in\mathbb{H}$ and $x\in X$, the proof of the lower bound basically follows from the Laplace theorem : 
\begin{equation}\label{eq:lower-bound}
\operatorname{Re}(\omega)\,\bigl\|\omega^{\alpha-1}(\omega^\alpha-L)^{-1}x\bigr\|
=\operatorname{Re}(\omega)\,\Bigl\|\int_0^\infty e^{-\omega t}S(t)x\,dt\Bigr\|
\le \operatorname{Re}(\omega)\int_0^\infty e^{-t\operatorname{Re}\omega}\,dt\; C_L\|x\|
= C_L\|x\|.
\end{equation}
Taking the supremum over $\omega\in\mathbb{H}$ and $\|x\|\le 1$ gives
$K_\alpha(L)\le C_L$ by Definition~\ref{def:kreiss}. It remains to justify the
equality in \eqref{eq:lower-bound}, which is the Laplace formula of the theorem
and includes the claim $\omega^\alpha\in\rho(L)$.
\newline
\newline
Since $\|S(t)\|\le C_L$, the integral
$\widehat S(\omega):=\int_0^\infty e^{-\omega t}S(t)\,dt$ converges absolutely for
every $\omega\in\mathbb{H}$. It defines a holomorphic function
$\widehat S\colon\mathbb{H}\to B(X)$, being locally bounded and strongly holomorphic.
The map $\omega\mapsto\omega^\alpha$ is a biholomorphism of $\mathbb{H}$ onto
$\Sigma_\alpha$, so
\[
R(z):=z^{(1-\alpha)/\alpha}\,\widehat S\bigl(z^{1/\alpha}\bigr)
\]
is holomorphic on $\Sigma_\alpha$. By \eqref{eq:laplace}, $R(z)=(z-L)^{-1}$
for real $z>w^\alpha$. 
\newline
\newline
Let $U$ be the connected component of the open set $\rho(L)\cap\Sigma_\alpha$
that contains $(w^\alpha,\infty)$. By the identity theorem, $R(z)=(z-L)^{-1}$ on
$U$. Suppose $z_*\in\Sigma_\alpha$ is a boundary point of $U$. Then
$z_*\in\sigma(L)$, and as $z\to z_*$ within $U$,
\[
\|R(z)\|=\|(z-L)^{-1}\|\ge \operatorname{dist}\bigl(z,\sigma(L)\bigr)^{-1}\to\infty .
\]
This contradicts the continuity of $R$ at $z_*$. Hence $U=\Sigma_\alpha$, that is,
$\Sigma_\alpha\subset\rho(L)$. Moreover,
$\omega^{\alpha-1}(\omega^\alpha-L)^{-1}=\omega^{\alpha-1}R(\omega^\alpha)=\widehat S(\omega)$
for every $\omega\in\mathbb{H}$, which is the equality in \eqref{eq:lower-bound}.
\end{proof}

\subsection{Finite-rank operators and their slow projection}
A finite-rank operator $L=UV$, with $U:\C^n\to X$ and $V:X\to\C^n$, acts nontrivially only through the
small matrix $\widehat L:=VU$, and the identity $(z-UV)^{-1}=z^{-1}(I+U(z-VU)^{-1}V)$ reduces its resolvent
to that of $\widehat L$. The next lemma collects what we need from this reduction, and it identifies the
slow component of a finite-rank operator.

\begin{lemma}\label{lem:struct}
Let $L\in\B(X)$ have rank $n\ge1$. Choose a basis $x_1,\dots,x_n$ of $\Ran L$, write
$Lx=\sum_{j=1}^n\varphi_j(x)x_j$ with $\varphi_j\in X^*$, and let $U:\C^n\to X$, $Uc:=\sum_jc_jx_j$, and
$V:X\to\C^n$, $Vx:=(\varphi_j(x))_j$, so that $L=UV$ with $U$ injective and $V$ surjective. Put
$\widehat L:=VU\in\C^{n\times n}$.
\begin{enumerate}[label=\textup{(\alph*)}]
\item $\sigma(L)\setminus\{0\}=\sigma(\widehat L)\setminus\{0\}$, and for $z\notin\sigma(\widehat L)\cup\{0\}$,
$(z-L)^{-1}=z^{-1}\big(I+U(z-\widehat L)^{-1}V\big)$.
\item If $\Kr_\alpha(L)<\infty$, then $\widehat L$ is invertible, $\Pi:=U\widehat L^{-1}V$ is a projection with
$\Ran\Pi=\Ran L$ and $\ker\Pi=\ker L$, so that $X=\Ran L\oplus\ker L$, and $L\Pi=\Pi L=L$. Moreover
\[
E_\alpha(t^\alpha L)=U\,E_\alpha(t^\alpha\widehat L)\,\widehat L^{-1}V+(I-\Pi)\qquad(t\ge0).
\]
\item For $x\in X$ and $x^*\in X^*$ the function $F(z):=z\langle(z-L)^{-1}x,x^*\rangle$ extends to a proper rational
function whose reduced denominator has degree at most $n$ and whose poles lie in $\sigma(\widehat L)$. If
$\Kr_\alpha(L)<\infty$, then $F\in\RR_{\alpha,n}\cup\{0\}$, $F$ is regular at $0$, and
$F(0)=\langle(I-\Pi)x,x^*\rangle$.
\end{enumerate}
\end{lemma}


\begin{proof}
(a) The identity is verified by multiplying out, and it shows that every
$z\notin\sigma(\widehat L)\cup\{0\}$ lies in $\rho(L)$; conversely, if $\widehat Lc=zc$ with $z\ne0$ and $c\ne0$, then
$L(Uc)=zUc$ and $Uc\ne0$.
\newline
\newline
(b) Suppose $\widehat L$ is singular, so that $(s-\widehat L)^{-1}$ has a pole at $s=0$ and
$\|(s-\widehat L)^{-1}\|\ge c_0s^{-1}$ for small $s>0$. Since $U$ is injective and $V$ is surjective, there is
$c_1>0$ with $\|UMV\|\ge c_1\|M\|$ for all $M\in\C^{n\times n}$, and (a) gives
$s\|(s-L)^{-1}\|\ge c_1\|(s-\widehat L)^{-1}\|-1\to\infty$ as $s\downarrow0$.
\newline
\newline
Since $a_\alpha(s)=s$, this
contradicts $\Kr_\alpha(L)<\infty$. Hence $\widehat L$ is invertible, and 
\[
\Pi^2=U\widehat L^{-1}\widehat L\widehat L^{-1}V=\Pi,
L\Pi=UVU\widehat L^{-1}V=L \\
\Pi L=U\widehat L^{-1}VUV=L
\]
The range of $\Pi$ is $\Ran U=\Ran L$, and
$\Pi x=0$ if and only if $Vx=0$, that is, $Lx=0$. Finally, $L^k=U\widehat L^{k-1}V$ for $k\ge1$, so
\[
E_\alpha(t^\alpha L)=I+U\Big(\sum_{k\ge1}\frac{t^{\alpha k}\widehat L^{k-1}}{\Gamma(1+\alpha k)}\Big)V
=I+U\big(E_\alpha(t^\alpha\widehat L)-I\big)\widehat L^{-1}V,
\]
which is the stated formula.

(c) By (a), $F(z)=\langle x,x^*\rangle+\ell^\top(z-\widehat L)^{-1}k$ with $k:=Vx$ and $\ell_i:=\langle x_i,x^*\rangle$. By
Cramer's rule,
\[
F(z)=\frac{\langle x,x^*\rangle\det(z-\widehat L)+\ell^\top\operatorname{adj}(z-\widehat L)\,k}{\det(z-\widehat L)},
\]
which is proper with denominator of degree $n$; cancelling common factors preserves properness, and the
poles lie in $\sigma(\widehat L)$. If $\Kr_\alpha(L)<\infty$, then $\sigma(L)\cap\Sa=\emptyset$, so
$\sigma(\widehat L)\cap\Sa=\emptyset$ by (a) and $0\notin\Sa$; and $\widehat L$ is invertible by (b), so $F$ is regular
at $0$ with $F(0)=\langle x,x^*\rangle-\ell^\top\widehat L^{-1}Vx=\langle(I-\Pi)x,x^*\rangle$.
\end{proof}

By (b), the family acts on $\Ran L$ through the $n\times n$ matrix $\widehat L$ and leaves the kernel of $L$
fixed. The two subspaces need not be orthogonal, so the propagator of the slow part is
$U E_\alpha(t^\alpha\widehat L)\widehat L^{-1}V=E_\alpha(t^\alpha L)\Pi$, an operator on all of $X$ whose norm feels the angle
between $\Ran L$ and $\ker L$. This obliqueness of the splitting already shows up on both sides of \Cref{thm:A}.

\begin{remark}[Obliqueness]\label{rem:oblique}
If $\Kr_\alpha(L)<\infty$, then $\|I-\Pi\|\le\Kr_\alpha(L)$. Indeed, by \cref{lem:struct}(a),
$s(s-L)^{-1}=I+U(s-\widehat L)^{-1}V\to I-U\widehat L^{-1}V=I-\Pi$ as $s\downarrow0$, and $a_\alpha(s)=s$. If moreover every
eigenvalue of $\widehat L$ lies in the interior of $\La$, then $E_\alpha(t^\alpha\widehat L)\to0$, because in the Jordan form
of $\widehat L$ the entries are multiples of $t^{\alpha j}E_\alpha^{(j)}(t^\alpha\lambda)=O(t^{-\alpha})$ by
\cref{lem:ML}(c) (trivially for $\alpha=1$), and \cref{lem:struct}(b) gives $E_\alpha(t^\alpha L)\to I-\Pi$ in norm. A strongly oblique
splitting therefore forces a large Kreiss constant and a large stationary state, whatever the spectrum of $\widehat L$.
\end{remark}

\subsection{The upper bound}

\begin{theorem}\label{thm:finite-rank}
Let $0<\alpha\le1$, let $X\ne\{0\}$ be a Banach space, and let $L\in\B(X)$ have rank at most $n$. Then
\begin{equation}\label{eq:kreiss}
\Kr_\alpha(L)\;\le\;\sup_{t\ge0}\|E_\alpha(t^\alpha L)\|\;\le\;\kappa_{\alpha,n}\,\Kr_\alpha(L),\qquad
\kappa_{\alpha,n}=\frac e2\big(3+4n(2-\alpha)\big).
\end{equation}
In particular, $\sup_t\|E_\alpha(t^\alpha L)\|<\infty$ if and only if $\Kr_\alpha(L)<\infty$.
\end{theorem}

\begin{proof}
The lower bound is \cref{thm:kreiss-lower}. For the upper bound we may assume $\Kr:=\Kr_\alpha(L)<\infty$, so
that $\sigma(L)\cap\Sa=\emptyset$. Since $t=0$ the claim reads $1\le\kappa_{\alpha,n}\Kr$, which holds by
\cref{lem:elementary}(a), and if $L=0$, then $E_\alpha(t^\alpha L)=I$, we might assume that $1\le\rank L\le n$.
Then fix unit vectors
$x\in X$ and $x^*\in X^*$. 
\newline
\newline
Define the function $g : \mathbb{R}_{t \geq 0} \to \mathbb{C}$ to be 
\[
g(t):=\langle E_\alpha(t^\alpha L)x,x^*\rangle
\]
and the function $\widehat g : \HH \to \mathbb{C}$ as 
\[
\widehat g(\omega) := \omega^{\alpha-1}\langle(\omega^\alpha-L)^{-1}x,x^*\rangle
\quad\text{for } \omega\in\HH.
\]
Since $\{\omega^\alpha:\omega\in\HH\}=\Sa\subset\rho(L)$,
$\widehat g$ is holomorphic on $\HH$, and by \cref{lem:struct}(c) it equals $h_F$, with the corresponding  $F \in \RR_{\alpha,n}\cup\{0\}$ being $F(z)=z\langle(z-L)^{-1}x,x^*\rangle$. Let $R:=(\omega^\alpha-L)^{-1}$, then the following holds : 
\[
\widehat g\,'(\omega)=(\alpha-1)\omega^{\alpha-2}\langle Rx,x^*\rangle-\alpha\omega^{2\alpha-2}\langle R^2x,x^*\rangle ,
\]
and $\|R\|\le2|\omega|^{-\alpha}$ for $|\omega|^\alpha\ge2\|L\|$, so that $|\widehat g\,'(\omega)|\le4\|x\|\|x^*\||\omega|^{-2}$
\newline
\newline
Let $c:=\|L\|^{\frac{1}{\alpha}}$. By \cref{lem:families}(a),
\[
\int_0^\infty e^{-\omega t}g(t)\,dt=\widehat g(\omega)
\]
for $\RE\omega>c$, with absolute convergence, and differentiation under
the integral gives
\[
\int_0^\infty e^{-\omega t}t\,g(t)\,dt=-\widehat g\,'(\omega).
\]
For $b>c$ the function
$t\mapsto e^{-bt}t\,g(t)\mathbf 1_{[0,\infty)}(t)$ is continuous and integrable on $\R$, its Fourier transform
$y\mapsto-\widehat g\,'(b+iy)$ is continuous and
integrable. The Fourier inversion theorem gives
the formula below on the line $\RE\omega=b$. For $0<a<b$, Cauchy's theorem on the rectangles with vertices $a\pm iY$, $b\pm iY$,
$Y\ge(2\|L\|)^{\frac{1}{\alpha}}$, whose horizontal sides contribute at most $(b-a)e^{bt}\,4\|x\|\|x^*\|Y^{-2}$, moves it to every line $\RE\omega=a>0$:
\begin{equation}\label{eq:inversion}
\langle E_\alpha(t^\alpha L)x,x^*\rangle=-\frac{1}{2\pi i\,t}\int_{a-i\infty}^{a+i\infty}e^{\omega t}\,\widehat g\,'(\omega)\,d\omega
\qquad(t>0,\ a>0),
\end{equation}
with an absolutely convergent integral. The integration by parts hidden in the factor $t$ is what makes the
integrand integrable, and it is the reason why the variation of $h_F$, and not $h_F$ itself, enters the bound. Apart from
the identification $\widehat g=h_F$, Steps~1 and~2 use only that $L\in\B(X)$ and $\Sa\subset\rho(L)$, so \eqref{eq:inversion} holds
for every such operator and all $x\in X$, $x^*\in X^*$.
\newline
\newline
Otherwise \eqref{eq:inversion} and
\cref{thm:bernstein} give, for every $a>0$ and $t>0$,
\[
|g(t)|\le\frac{e^{at}}{2\pi t}\int_\R|h_F'(a+iy)|\,dy
\le\frac{e^{at}}{2\pi t}\,\pi\big(3+4n(2-\alpha)\big)\sup_{y\in\R}|h_F(a+iy)|.
\]
Since $|h_F(\omega)|\le\|\omega^{\alpha-1}(\omega^\alpha-L)^{-1}\|\le\frac{\Kr}{\RE\omega}$, we obtain
\[
|g(t)|\le\frac{3+4n(2-\alpha)}{2}\cdot\frac{e^{at}}{at}\cdot\Kr .
\]
The choice $a=\frac{1}{t}$ minimises $\frac{e^{at}}{at}$, with value $e$, and the supremum over $x$ and $x^*$ gives
\eqref{eq:kreiss}.
\end{proof}

\begin{remark}\label{rem:finite-rank}
(i) \Cref{thm:finite-rank} holds on every Banach space and uses no Hilbert-space structure, and only the rank
of $L$ enters, not the dimension of $X$: for $L\in\C^{d\times d}$ with $\rank L=n<d$ the pole of the resolvent at
$0$ is absorbed by the factor $z$ in $F(z)=z\langle(z-L)^{-1}x,x^*\rangle$. For $\alpha=1$, Spijker's lemma replaces
$\pi(3+4n)$ by $2\pi(n+1)$ and the same argument gives $\sup_t\|e^{tL}\|\le e(n+1)\Kr_1(L)$ for operators of rank $n$
on Banach spaces. \Cref{thm:bernstein} replaces Spijker's lemma for every $\alpha<1$, where no algebraic
argument is available.
(ii) The constant is $\kappa_{\alpha,n}$ and not $\kappa'_{\alpha,n}$ because $F(0)=\langle(I-\Pi)x,x^*\rangle$ need not
vanish: the fixed kernel contributes a stationary term to every matrix element. If one estimates the slow
part $E_\alpha(t^\alpha L)\Pi$ separately, \cref{thm:slow} below gives the refined constant $\kappa'_{\alpha,n}$ for it,
at the price of the term $\|I-\Pi\|$ for the kernel; by \cref{rem:oblique} this term never exceeds the left-hand
side of \eqref{eq:kreiss}.
(iii) By \cref{lem:elementary}(b), $\Kr_\alpha(L)$ may be finite while $\Kr_1(L)=\infty$. Then $e^{tL}$ is unbounded,
and it grows exponentially when $\sigma(L)$ meets $\HH$, so neither subordination nor \eqref{eq:classical} gives
any information about $E_\alpha(t^\alpha L)$.
\end{remark}

\section{Hilbert spaces: normal operators and the numerical range}\label{sec:hilbert}

On a Hilbert space the fast part of a generator is often controlled by its numerical range, and the
question is what the numerical range can say about the Mittag--Leffler family. For semigroups the answer
is the Lumer--Phillips theorem: a dissipative generator has a contractive semigroup. For $\alpha<1$ the
family of a dissipative generator is still contractive (\cref{lem:families}(c)), but the stability region
$\La$ is larger than the left half-plane, and the relevant question becomes what happens when the numerical
range avoids only the sector $\Sa$. We first settle the normal case, which identifies the sharp constant.

\subsection{Normal operators}

\begin{proposition}[Normal operators]\label{prop:normal}
Let $0<\alpha\le1$ and let $L$ be a normal operator on $H$, possibly unbounded, with spectral measure $P_L$ and
$\sigma(L)\subset\La$. Then $L$ generates a Mittag--Leffler family,
\[
E_\alpha(t^\alpha L)=\int_{\sigma(L)}E_\alpha(t^\alpha z)\,dP_L(z),\qquad\|E_\alpha(t^\alpha L)\|=\sup_{z\in\sigma(L)}|E_\alpha(t^\alpha z)|,
\]
\eqref{eq:laplace} holds for every $\lambda>0$, and
$1\le\Kr_\alpha(L)\le1/\alpha$,$\sup_{t\ge0}\|E_\alpha(t^\alpha L)\|\le\Mst$,
independently of $\dim H$. Both bounds are sharp: if $\sigma(L)$ contains a nonzero point of $\partial\Sa$, then
$\Kr_\alpha(L)=\frac{1}{\alpha}$ and $\sup_t\|E_\alpha(t^\alpha L)\|=\Mst$. In particular, for $0<\alpha<1$ normal operators exhibit
transient growth, and the ratio of the two sides of \Cref{thm:A} equals $\alpha\Mst>1$ for them.
\end{proposition}

\begin{proof}
Since $t^\alpha\La=\La$ and $E_\alpha$ is bounded and continuous on $\La$, the operators
$S(t):=\int E_\alpha(t^\alpha z)\,dP_L(z)$ are bounded, with norm equal to the supremum of $|E_\alpha(t^\alpha\cdot)|$ over the
support $\sigma(L)$ of $P_L$, and $S$ is strongly continuous by dominated convergence. For $z\in\La$ the functions
$\lambda\mapsto\int_0^\infty e^{-\lambda t}E_\alpha(t^\alpha z)\,dt$ and $\lambda\mapsto\lambda^{\alpha-1}(\lambda^\alpha-z)^{-1}$ are holomorphic on $\HH$,
since the integrand is bounded and $\lambda^\alpha\in\Sa\not\ni z$, and they agree for $\lambda>|z|^{\frac{1}{\alpha}}$ by
\cref{lem:families}(a), hence for all $\lambda>0$. Fubini's theorem with respect to the complex measures
$\langle P_L(\cdot)x,y\rangle$ gives \eqref{eq:laplace}, and $\|S(t)\|\le\Mst$.

For $z\in\Sa$ we have $\|(z-L)^{-1}\|=\dist(z,\sigma(L))^{-1}\le\dist(z,\partial\Sa)^{-1}$, so $\Kr_\alpha(L)\le\frac{1}{\alpha}$ by
\cref{lem:elementary}(c), and $\Kr_\alpha(L)\ge1$ by \cref{lem:elementary}(a). If $\lambda=re^{i\frac{\alpha\pi}{2}}\in\sigma(L)$ with
$r>0$ (the lower ray is treated in the same way, by conjugation symmetry), the points $z=\lambda+s\mathbf n$, with $\mathbf n$ the inner unit normal of the boundary ray at $\lambda$, satisfy
$\|(z-L)^{-1}\|\ge s^{-1}$ and $\frac{a_\alpha(z)}{s}\to\frac{1}{\alpha}$ as $s\downarrow0$, by the computation in the proof of
\cref{lem:elementary}(c). Moreover $t^\alpha\lambda$ runs over the whole boundary ray as $t$ varies, so
$\sup_t\|E_\alpha(t^\alpha L)\|\ge M_\alpha(\frac{\alpha\pi}{2})=\Mst$ by \cref{lem:mstar}.
\end{proof}

\begin{remark}[Generators without a semigroup]\label{rem:normal-ray}
Let $0<\alpha<1$, let $A\ne0$ be a nonnegative self-adjoint operator and $\frac{\alpha\pi}{2}\le\theta<\frac{\pi}{2}$. Then
$L:=e^{\pm i\theta}A$ is normal with spectrum on a ray in $\La$, so it generates a Mittag--Leffler family with
$\sup_t\|E_\alpha(t^\alpha L)\|=M_\alpha(\theta)$ and $\Kr_\alpha(L)\le\frac{1}{\alpha}$; if $A$ is unbounded, then $\RE\sigma(L)$ is unbounded
above and $L$ generates no $C_0$-semigroup \citep{EngelNagel2000}. The time-fractional Schr\"odinger equation
$i^\alpha\partial_t^\alpha u=A_Vu$ of \citet{Naber2004}, with a Schr\"odinger operator $A_V$, corresponds to
$L=i^{-\alpha}A_V=e^{-i\frac{\alpha\pi}{2}}A_V$. If $A_V$ is self-adjoint,
bounded below and unbounded above, the positive part of its spectrum is mapped onto the lower boundary ray of
$\Sa$ and its negative eigenvalues onto the ray of argument $\pi-\frac{\alpha\pi}{2}$, so the solution operators are
uniformly bounded with $\sup_t\|E_\alpha(t^\alpha L)\|=\Mst$ exactly, although $e^{tL}$ does not exist. This is prototype
(b) of the introduction, to which \cref{sec:exp-schroedinger} returns.
\end{remark}

\subsection{The numerical range}

For a closed, densely defined operator $L$ on $H$ let $W(L):=\{\langle Lv,v\rangle:v\in D(L),\ \|v\|=1\}$ be its numerical
range, which is convex, and let $C_W$ be the closed convex cone generated by it,i.e. the closure of
$\{sw:s\ge0,\ w\in W(L)\}$. If $W(L)$ does not meet the open cone $\Sa$, neither does $C_W$.  The following
geometric lemma shows that a half-plane, with spectral-set constant being equal to one, 
is as good as any convex set for the functions $E_\alpha(t^\alpha\,\cdot)$.

\begin{lemma}\label{lem:cone}
Let $0<\alpha\le1$ and let $C\ne\C$ be a closed convex cone with vertex $0$ that does not meet $\Sa$. There is a
closed half-plane $\Omega_0$ with $0\in\partial\Omega_0$, $C\subset\Omega_0$ and $\Omega_0\cap\Sa=\emptyset$ such that
$\sup_{z\in\Omega_0}|E_\alpha(z)|=\sup_{z\in C}|E_\alpha(z)|$. This common value is the larger of $M_\alpha(|\theta|)$ over the
arguments $\theta\in(-\pi,\pi]$ of the boundary rays of $C$, and it is at most $\Mst$.
\end{lemma}

\begin{proof}
For $\alpha=1$ every such cone lies in the closed left half-plane, on which $|e^z|\le1$, so we let $\alpha<1$; we
may also assume $C\ne\{0\}$. If $C$ is a line, it bounds a closed half-plane disjoint from $\Sa$, on which the
supremum is attained on the line by \cref{lem:mstar}. Otherwise $C=\{re^{i\theta}:r\ge0,\ \theta_1\le\theta\le\theta_2\}$
with $\frac{\alpha\pi}{2}\le\theta_1  -\theta_2\le2\pi-\frac{\alpha\pi}{2}$ and $\theta_2-\theta_1\le\pi$. Let $d(\theta):=\min\{\theta,2\pi-\theta\}$
and $\delta_0:=\min\{d(\theta_1),d(\theta_2)\}\ge\frac{\alpha\pi}{2}$. By \cref{lem:mstar}, $\sup_C|E_\alpha|=M_\alpha(\delta_0)$, and
for $\phi\in[\frac{\alpha\pi}{2},\pi-\frac{\alpha\pi}{2}]$ the half-plane $\Omega_\phi:=\{re^{i\theta}:r\ge0,\ \phi\le\theta\le\phi+\pi\}$
does not meet $\Sa$ and satisfies $\sup_{\Omega_\phi}|E_\alpha|=M_\alpha(\min\{\phi,\pi-\phi\})$. If $\delta_0\ge\frac{\pi}{2}$, then
$C\subset\Omega_{\frac{\pi}{2}}$ and both suprema equal $1$. If $\delta_0<\frac{\pi}{2}$, choose
\[
\phi\in\big[\max\{\theta_2-\pi,\delta_0\},\ \min\{\theta_1,\pi-\delta_0\}\big],
\]
an interval that is nonempty because $\theta_2-\pi\le\theta_1$, $\theta_2\le2\pi-d(\theta_2)\le2\pi-\delta_0$,
$\delta_0\le d(\theta_1)\le\theta_1$ and $\delta_0<\frac{\pi}{2}$. Then $\phi\le\theta_1$ and $\phi+\pi\ge\theta_2$, so $C\subset\Omega_\phi$, and
$\min\{\phi,\pi-\phi\}\ge\delta_0$, so $\sup_{\Omega_\phi}|E_\alpha|\le M_\alpha(\delta_0)=\sup_C|E_\alpha|$ by monotonicity of
$M_\alpha$. The reverse inequality holds because $C\subset\Omega_\phi$.
\end{proof}

\begin{theorem}[Numerical range]\label{thm:numrange}
Let $0<\alpha\le1$ and let $L$ be a closed, densely defined operator on $H$ with $W(L)\cap\Sa=\emptyset$ and
$\rho(L)\cap\Sa\ne\emptyset$.
\begin{enumerate}[label=\textup{(\alph*)}]
\item $\Sa\subset\rho(L)$, $\|(z-L)^{-1}\|\le\dist(z,W(L))^{-1}$ for $z\in\Sa$, and $1\le\Kr_\alpha(L)\le\frac{1}{\alpha}$.
\item Let $\Omega$ be a closed half-plane with $W(L)\subset\Omega$ and $\Omega\cap\Sa=\emptyset$. Then $L$ generates a
Mittag--Leffler family, and
\begin{equation}\label{eq:vN}
\|E_\alpha(t^\alpha L)\|\le\sup_{z\in\partial\Omega}|E_\alpha(t^\alpha z)|\qquad(t\ge0).
\end{equation}
\item Consequently
\[
\sup_{t\ge0}\|E_\alpha(t^\alpha L)\|\le\sup_{z\in C_W}|E_\alpha(z)|\le\Mst,
\]
and for every closed convex cone $C\ne\C$ disjoint from $\Sa$ the first bound is attained by a normal operator
whose numerical range generates $C$. If $W(L)\subset\{\RE z\le-\delta\}$ with $\delta\ge0$, then $\|E_\alpha(t^\alpha L)\|\le G_\alpha(\delta t^\alpha)$,
where $G_\alpha(s):=\sup_{y\in\R}|E_\alpha(-s+iy)|$. In particular, if $L$ is dissipative, then
$\Kr_\alpha(L)=\sup_t\|E_\alpha(t^\alpha L)\|=1$.
\end{enumerate}
\end{theorem}

For $0<\alpha<1$ we have $E_\alpha(-s)\le G_\alpha(s)\le\min\{1,\frac{C_\alpha}{s}\}$ for $s>0$, by \cref{lem:mstar} and
\cref{lem:ML}(b), so a uniformly dissipative generator satisfies
$\|E_\alpha(t^\alpha L)\|\le\min\{1,C_\alpha\delta^{-1}t^{-\alpha}\}$, and the order $t^{-\alpha}$ is attained by $L=-\delta$. For
$\alpha=1$, $G_1(s)=e^{-s}$ and \eqref{eq:vN} reduces to the Lumer--Phillips bound.

\begin{proof}
(a) For $z\notin\overline{W(L)}$ and unit vectors $v\in D(L)$, $\|(z-L)v\|\ge|\langle(z-L)v,v\rangle|\ge\dist(z,W(L))$, so
$z-L$ is injective with closed range, and by \citet[Theorem~V.3.2]{Kato1995} the codimension of its range is
constant on each connected component of $\C\setminus\overline{W(L)}$. The open connected set $\Sa$ does not meet
$\overline{W(L)}$ and meets $\rho(L)$, so it lies in a component on which $z-L$ is bijective. This gives
$\Sa\subset\rho(L)$ and the resolvent bound. Since $\dist(z,W(L))\ge\dist(z,\La)=\dist(z,\partial\Sa)$ for $z\in\Sa$, the
bound $\Kr_\alpha(L)\le\frac{1}{\alpha}$ follows from \cref{lem:elementary}(c).

(b) Write $\Omega=\{z:\RE(e^{-i\chi}z)\le-\delta\}$. Since $\Omega\cap\Sa=\emptyset$ and $\Sa$ is a cone, replacing $z$ by
$rz$ and letting $r\to0$ and $r\to\infty$ shows that $\delta\ge0$ and $|\chi|\le\frac{(1-\alpha)\pi}{2}$. Put
$T:=e^{-i\chi}L+\delta$ with $D(T)=D(L)$, so that $\RE\langle Tv,v\rangle\le0$ for $v\in D(L)$. The open half-plane
$V:=\{z:\RE(e^{-i\chi}z)>-\delta\}$ is connected, does not meet $\overline{W(L)}$, and contains
$1\in\Sa\subset\rho(L)$ because $\cos\chi>0$; by the argument in (a), $V\subset\rho(L)$. For $\lambda>0$ the point
$e^{i\chi}(\lambda-\delta)$ lies in $V$, so $\lambda-T=e^{-i\chi}\big(e^{i\chi}(\lambda-\delta)-L\big)$ is bijective. Thus $T$ is
m-dissipative and generates a $C_0$-semigroup of contractions. By the dilation theorem of Sz.-Nagy 
\citep{SzNagyFoias2010} there are a Hilbert space $K\supset H$ and a strongly continuous unitary
group $(U(s))_{s\in\R}$ on $K$ such that $e^{sT}=P_HU(s)|_H$ for $s\ge0$, where $P_H$ is the orthogonal projection
onto $H$. Let $\mathcal G$ be the skew-adjoint generator of $(U(s))$. Taking Laplace transforms,
\begin{equation}\label{eq:dilation-res}
(\varsigma-T)^{-1}h=P_H(\varsigma-\mathcal G)^{-1}h\qquad(h\in H,\ \RE\varsigma>0).
\end{equation}
The operator $N:=e^{i\chi}(\mathcal G-\delta)$ is normal on $K$ with $\sigma(N)\subset e^{i\chi}(i\R-\delta)=\partial\Omega\subset\La$. By
\cref{prop:normal}, $N$ generates a Mittag--Leffler family with
$\|E_\alpha(t^\alpha N)\|\le\sup_{z\in\partial\Omega}|E_\alpha(t^\alpha z)|$. Put $S(t):=P_HE_\alpha(t^\alpha N)|_H$, which is strongly
continuous and satisfies \eqref{eq:vN}. Let $\lambda>0$ and $\varsigma:=e^{-i\chi}\lambda^\alpha+\delta$, so that
$\RE\varsigma=\lambda^\alpha\cos\chi+\delta>0$, $\lambda^\alpha-N=e^{i\chi}(\varsigma-\mathcal G)$ and $\varsigma-T=e^{-i\chi}(\lambda^\alpha-L)$. By \cref{prop:normal}
and \eqref{eq:dilation-res}, for $h\in H$,
\[
\int_0^\infty e^{-\lambda t}S(t)h\,dt=\lambda^{\alpha-1}e^{-i\chi}P_H(\varsigma-\mathcal G)^{-1}h
=\lambda^{\alpha-1}e^{-i\chi}(\varsigma-T)^{-1}h=\lambda^{\alpha-1}(\lambda^\alpha-L)^{-1}h,
\]
so $L$ generates the Mittag--Leffler family $S$.

(c) By \cref{lem:cone} there is a closed half-plane $\Omega_0\supset C_W\supset W(L)$ with $0\in\partial\Omega_0$, $\Omega_0\cap\Sa=\emptyset$
and $\sup_{\Omega_0}|E_\alpha|=\sup_{C_W}|E_\alpha|$. Since $t^\alpha\partial\Omega_0=\partial\Omega_0$, \eqref{eq:vN} with $\Omega=\Omega_0$ gives the first bound,
and the second is \cref{lem:cone}. For sharpness, let $\Gamma$ be the union of the boundary rays of a given cone $C$ and let $L$ be multiplication by
$z$ on $L^2(\Gamma,|dz|)$. It is normal with $\sigma(L)=\Gamma$, the cone generated by $W(L)$ is $C$, and by
\cref{prop:normal} and \cref{lem:mstar} its family has supremum $\sup_{z\in\Gamma}|E_\alpha(z)|=\sup_C|E_\alpha|$. For the second claim take $\Omega=\{\RE z\le-\delta\}$ and note
$t^\alpha\partial\Omega=\{\RE z=-\delta t^\alpha\}$. If $L$ is dissipative, the case $\delta=0$ and $M_\alpha(\frac{\pi}{2})=1$ give
$\|E_\alpha(t^\alpha L)\|\le1$, and $1\le\Kr_\alpha(L)\le\sup_t\|E_\alpha(t^\alpha L)\|$ by (a) and \cref{thm:kreiss-lower}.
\end{proof}


Estimate \eqref{eq:vN} says that a closed half-plane containing $W(L)$ and disjoint from $\Sa$ is a spectral set
with constant one for the functions $E_\alpha(t^\alpha\,\cdot)$, which is von Neumann's inequality for a half-plane 
\citep{vonNeumann1951}. In addition, we note that \citet{CrouzeixPalencia2017} proved that the closure of the numerical range of a bounded
operator is a 
$(1+\sqrt2)$-spectral set, and \citet[Theorem~3]{Crouzeix2025} showed that an unbounded closed
convex set $\Omega\ne\C$ containing the numerical range, with aperture $2\beta$, the maximal opening of a sector
contained in $\Omega$, is a $K(\beta)$-spectral set for the rational functions bounded on $\Omega$ and for the
holomorphic functions on $\Omega$ that are continuous on $\overline\Omega$ and have a limit at infinity, with
$K(\beta)=1-\beta/\pi+\sqrt{2-{4\beta}/\pi+{\beta^2}/{\pi^2}}$
his Section~3 extends this to closed, densely defined operators. 
\section{Generators with a finite slow component}\label{sec:slow}

We now remove the assumption that the whole operator has finite rank, and allow the remaining dynamics to
be unbounded and of infinite rank, provided it is bounded for a structural reason such as dissipativity or a
numerical range outside $\Sa$. A slow component (\cref{def:slow}) isolates finitely many modes, and the rank
dependence of \cref{thm:finite-rank} is confined to them. The geometry of the splitting, which may be very oblique, is carried by the deflated Kreiss constant
$\Kr_\alpha(L;P)$ and by the fast propagator $E_\alpha(t^\alpha L)(I-P)$, both of which can be computed.

\begin{theorem}[Kreiss theorem with a finite slow component]\label{thm:slow}
Let $0<\alpha\le1$, let $X\ne\{0\}$, let $L$ be closed and densely 
defined on $X$, and let $P\in\B(X)$ be a projection
of finite rank $m\ge1$ commuting with $L$, with slow block $L_P$ and $Q:=I-P$.
\begin{enumerate}[label=\textup{(\roman*)}]
\item The slow part of the propagator satisfies
\begin{equation}\label{eq:compressed}
\Kr_\alpha(L;P)\;\le\;\sup_{t\ge0}\big\|E_\alpha(t^\alpha L_P)P\big\|\;\le\;\kappa_{\alpha,m}\,\Kr_\alpha(L;P),
\end{equation}
and $\kappa_{\alpha,m}$ may be replaced by $\kappa'_{\alpha,m}$ if $0\notin\sigma(L_P)$. If $L$ generates a Mittag--Leffler
family, then $E_\alpha(t^\alpha L)P=E_\alpha(t^\alpha L_P)P$.
\item If $P$ is a slow component, then $L$ generates a Mittag--Leffler family,
$E_\alpha(t^\alpha L)=E_\alpha(t^\alpha L_P)P+E_\alpha(t^\alpha L_Q)Q$, and with $C_Q:=\sup_{t\ge0}\|E_\alpha(t^\alpha L)Q\|$,
\begin{equation}\label{eq:deflation}
\Kr_\alpha(L)\;\le\;\sup_{t\ge0}\|E_\alpha(t^\alpha L)\|\;\le\;\kappa_{\alpha,m}\,\Kr_\alpha(L;P)+C_Q,
\end{equation}
with $\kappa'_{\alpha,m}$ in place of $\kappa_{\alpha,m}$ if $0\notin\sigma(L_P)$.
\item If $X=H$ and $P$ is an orthogonal slow component, then $\Kr_\alpha(L;P)=\Kr_\alpha(L_P)$ and
\[
\sup_{t\ge0}\|E_\alpha(t^\alpha L)\|=\max\Big\{\sup_{t\ge0}\|E_\alpha(t^\alpha L_P)\|,\ \sup_{t\ge0}\|E_\alpha(t^\alpha L_Q)\|\Big\}.
\]
\end{enumerate}
\end{theorem}

\begin{proof}
(i) Write $S_P(t):=E_\alpha(t^\alpha L_P)P$ and $C:=\sup_t\|S_P(t)\|$. By \cref{lem:slow-basic}(a), $L_P\in\B(PX)$, and
$\|E_\alpha(t^\alpha L_P)\|_{\B(PX)}\le\|S_P(t)\|$ because $y=Py$ for $y\in PX$. By \cref{lem:families}(a), $S_P$ is norm
continuous and
\begin{equation}\label{eq:compressed-laplace}
\int_0^\infty e^{-\omega t}S_P(t)\,dt=\omega^{\alpha-1}(\omega^\alpha-L_P)^{-1}P\qquad(\RE\omega>\|L_P\|^{\frac{1}{\alpha}}).
\end{equation}

\emph{Lower bound.} We may assume $C<\infty$. Then $\Sa\subset\rho(L_P)$ by \cref{thm:kreiss-lower} applied on $PX$, both
sides of \eqref{eq:compressed-laplace} are holomorphic on $\HH$ and agree there, and the estimate
$\RE(\omega)\|\int_0^\infty e^{-\omega t}S_P(t)\,dt\|\le C$ gives $\Kr_\alpha(L;P)\le C$.

\emph{Upper bound.} We may assume $\Kr:=\Kr_\alpha(L;P)<\infty$, so that $\sigma(L_P)\cap\Sa=\emptyset$. For $t=0$ the claim
reads $\|P\|\le\kappa_{\alpha,m}\Kr$, which holds by \cref{lem:slow-basic}(c). Let $t>0$ and let $x\in X$ and $x^*\in X^*$ be
unit vectors. Since $L_P\in\B(PX)$ and $\Sa\subset\rho(L_P)$, \eqref{eq:inversion} applies on the Banach space $PX$, with $L_P$, $Px$
and the restriction of $x^*$ in place of $L$, $x$ and $x^*$, and gives
\[
\langle S_P(t)x,x^*\rangle=-\frac1{2\pi i\,t}\int_{a-i\infty}^{a+i\infty}e^{\omega t}\widehat g\,'(\omega)\,d\omega,\qquad
\widehat g(\omega):=\omega^{\alpha-1}\langle(\omega^\alpha-L_P)^{-1}Px,x^*\rangle,
\]
for every $a>0$. Put $F(z):=z\langle(z-L_P)^{-1}Px,x^*\rangle$, so that $\widehat g=h_F$. With the factorisation of
\cref{lem:slow-basic}(b), $F(z)=z\,\ell^\top(z-\Lambda)^{-1}k$ with $k:=Vx$ and $\ell_i:=\langle Ue_i,x^*\rangle$, and Cramer's
rule shows that $F$ is a proper rational function with denominator $\det(z-\Lambda)$ of degree $m$, whose poles lie in
$\sigma(L_P)$ and hence outside $\Sa$. Thus $F\in\RR_{\alpha,m}\cup\{0\}$. If $0\notin\sigma(L_P)$, then $F$ is regular at $0$
and $F(0)=0$. If $F\equiv0$, then $\langle S_P(t)x,x^*\rangle=0$. Otherwise \cref{thm:bernstein} gives
\[
|\langle S_P(t)x,x^*\rangle|\le\frac{e^{at}}{2\pi t}\int_\R|h_F'(a+iy)|\,dy\le\frac{c}{2\pi}\cdot\frac{e^{at}}{t}\sup_{y\in\R}|h_F(a+iy)|,
\]
with $c=\pi(3+4m(2-\alpha))$, or $c=\pi(1+4m(2-\alpha))$ if $0\notin\sigma(L_P)$, and
$|h_F(\omega)|\le\|\omega^{\alpha-1}(\omega^\alpha-L_P)^{-1}P\|\le\frac{\Kr}{\RE\omega}$. The choice $a=\frac{1}{t}$ and the supremum over $x$ and
$x^*$ give $\|S_P(t)\|\le\frac{ce}{2\pi}\Kr$, which is the claim.

Finally, if $L$ generates a Mittag--Leffler family, then
\[
\int_0^\infty e^{-\lambda t}E_\alpha(t^\alpha L)Px\,dt=\lambda^{\alpha-1}(\lambda^\alpha-L)^{-1}Px=\lambda^{\alpha-1}(\lambda^\alpha-L_P)^{-1}Px
\]
for large $\lambda$ by \cref{lem:slow-basic}(a), and uniqueness of the Laplace transform together with
\eqref{eq:compressed-laplace} gives $E_\alpha(t^\alpha L)P=S_P(t)$.

(ii) Let $S(t):=E_\alpha(t^\alpha L_P)P+E_\alpha(t^\alpha L_Q)Q$. This family is strongly continuous, and
$\int_0^\infty e^{-wt}\|S(t)\|\,dt<\infty$ for large $w$. For large $\lambda>0$ the Laplace transforms of the two terms are
$\lambda^{\alpha-1}(\lambda^\alpha-L_P)^{-1}Px$ and $\lambda^{\alpha-1}(\lambda^\alpha-L_Q)^{-1}Qx$, and by \cref{lem:slow-basic}(a) their sum is
$\lambda^{\alpha-1}(\lambda^\alpha-L)^{-1}x$. Hence $L$ generates $S$, and since $PQ=QP=0$,
\[
S(t)P=E_\alpha(t^\alpha L_P)P,\qquad S(t)Q=E_\alpha(t^\alpha L_Q)Q .
\]
The upper bound in \eqref{eq:deflation} follows from $\|S(t)\|\le\|S(t)P\|+\|S(t)Q\|$ and (i), and the
lower bound is trivial if the supremum is infinite and otherwise follows from \cref{thm:kreiss-lower}.

(iii) In the orthogonal case $S(t)$ is the orthogonal direct sum of $E_\alpha(t^\alpha L_P)$ and $E_\alpha(t^\alpha L_Q)$, and
$\Kr_\alpha(L;P)=\Kr_\alpha(L_P)$ by \cref{lem:slow-basic}(c).
\end{proof}

\Cref{thm:finite-rank} should be compared with the special case of \cref{thm:slow} in which the fast part is trivial: for a bounded
$L$ of rank $n$ with $\Kr_\alpha(L)<\infty$, the projection $\Pi$ of \cref{lem:struct} is a slow component with
$L_Q=0$ on $\ker L$, so \eqref{eq:deflation} gives $\sup_t\|E_\alpha(t^\alpha L)\|\le\kappa'_{\alpha,n}\Kr_\alpha(L;\Pi)+\|I-\Pi\|$.
\Cref{thm:finite-rank} avoids the term $\|I-\Pi\|$ altogether, at the price of the constant $\kappa_{\alpha,n}$ instead of
$\kappa'_{\alpha,n}$, and in this sense it absorbs the stationary kernel into the Kreiss constant.

The fast part can be controlled on every Banach space by sectoriality, and on Hilbert spaces by its numerical
range; in both cases \cref{thm:slow} becomes fully explicit.

\begin{corollary}\label{cor:sectorial-fast}
Let $P$ be a projection of finite rank $m\ge1$ commuting with $L$, and suppose that the fast part $L_Q$ satisfies the
hypothesis of \cref{lem:sectorial} on $QX$ for some $\psi\in(\frac{\alpha\pi}{2},\pi]$, with $\psi'$ and $\varphi$ as there. Then $P$ is a slow
component, $L$ generates a Mittag--Leffler family, and
\[
\Kr_\alpha(L)\le\sup_{t\ge0}\|E_\alpha(t^\alpha L)\|\le\kappa_{\alpha,m}\,\Kr_\alpha(L;P)+\|Q\|\,c_\varphi\,s_{\psi'}(L_Q),
\]
with $\kappa'_{\alpha,m}$ in place of $\kappa_{\alpha,m}$ if $0\notin\sigma(L_P)$. Neither $L$ nor $L_Q$ needs to generate a
$C_0$-semigroup.
\end{corollary}

\begin{proof}
By \cref{lem:sectorial}, $L_Q$ generates a Mittag--Leffler family bounded by $c_\varphi s_{\psi'}(L_Q)$, so $P$ is a slow
component, and \cref{thm:slow}(ii) applies with $C_Q\le\|Q\|c_\varphi s_{\psi'}(L_Q)$.
\end{proof}

\begin{corollary}\label{cor:deflation}
Let $X=H$ and let $P$ be a projection of finite rank $m\ge1$ commuting with $L$, with $Q=I-P$. Suppose that the
numerical range $W(L_Q)$, taken in the Hilbert space $QH$, does not meet $\Sa$, and that $\rho(L_Q)\cap\Sa\ne\emptyset$.
Let $C_Q^W$ be the closed convex cone generated by $W(L_Q)$. Then $P$ is a slow component, $L$ generates a
Mittag--Leffler family, and
\[
\Kr_\alpha(L)\;\le\;\sup_{t\ge0}\|E_\alpha(t^\alpha L)\|\;\le\;\kappa_{\alpha,m}\,\Kr_\alpha(L;P)+\|Q\|\sup_{z\in C_Q^W}|E_\alpha(z)|
\;\le\;\kappa_{\alpha,m}\,\Kr_\alpha(L;P)+\|Q\|\,\Mst ,
\]
with $\kappa'_{\alpha,m}$ in place of $\kappa_{\alpha,m}$ if $0\notin\sigma(L_P)$; if $W(L_Q)$ lies in the closed left half-plane,
the last term is $\|Q\|$. If $P$ is orthogonal, then
$\sup_t\|E_\alpha(t^\alpha L)\|\le\max\{\kappa_{\alpha,m}\Kr_\alpha(L_P),\sup_{C_Q^W}|E_\alpha|\}$.
\end{corollary}

\begin{proof}
By \cref{thm:numrange} applied to $L_Q$ on $QH$, $L_Q$ generates a Mittag--Leffler family with
$\sup_t\|E_\alpha(t^\alpha L_Q)\|\le\sup_{C_Q^W}|E_\alpha|\le\Mst$, and the bound is $1$ if $W(L_Q)\subset\{\RE z\le0\}$. The claims
follow from \cref{thm:slow}, since $C_Q\le\|Q\|\sup_t\|E_\alpha(t^\alpha L_Q)\|$.
\end{proof}

In \eqref{eq:deflation} both quantities on the right-hand side refer to the splitting rather than to the two parts alone.
Estimating them separately gives back the cruder bound
\begin{equation}\label{eq:crude}
\sup_{t\ge0}\|E_\alpha(t^\alpha L)\|\le\|P\|\,\kappa_{\alpha,m}\Kr_\alpha(L_P)+\|Q\|\sup_{t\ge0}\|E_\alpha(t^\alpha L_Q)\|,
\end{equation}
in which the non-normality between the slow modes and the remainder appears as a separate factor $\|P\|$. This factor
reflects a property of $L$ and not of a choice. If $P$ is the Riesz projection onto an isolated part of the spectrum and
$P'$ is another projection commuting with $L$ with $P'X=PX$, then $P'$ commutes with every resolvent and hence with $P$, so
$P'=PP'=P'P=P$: among projections with the given range, the splitting is unique and so is its obliqueness. 

We conclude this section by a
showing that no estimate in terms of the slow block and the fast part alone can hold true, since the obliqueness might the dominating part of the transient growth. To that end, we consider Recalling Prototype (a):
a two-dimensional oscillatory mode $a(t)\in\C^2$, a linearised self-excited oscillator with growth rate
$\cos\theta>0$ and frequency $\sin\theta$, which feeds a source into a subdiffusive medium on $(0,\pi)$ with absorbing
boundaries, both running on the same fractional clock:
\begin{equation}\label{eq:cascade}
\partial_t^\alpha a=L_1a,\qquad\partial_t^\alpha v=\partial_{xx}v+\mu\,(a_1b_1+a_2b_2),\quad
L_1=\begin{pmatrix}\cos\theta&-\sin\theta\\ \sin\theta&\cos\theta\end{pmatrix},\quad\frac{\alpha\pi}2<\theta<\frac\pi2,
\end{equation}
with linearly independent source profiles $b_1,b_2\in L^2(0,\pi)$ and a coupling strength $\mu\ge0$. On
$H:=\C^2\oplus L^2(0,\pi)$ the generator is $L(a,v)=(L_1a,Av+\mu \mathcal Sa)$, where $A$ is the Dirichlet Laplacian with
$D(A)=H^2(0,\pi)\cap H^1_0(0,\pi)$, $\mathcal Sa:=a_1b_1+a_2b_2$, and $D(L)=\C^2\oplus D(A)$. Let $\lambda_\pm=e^{\pm i\theta}$ be the
eigenvalues of $L_1$, with eigenvectors $e_\pm$, and let $Y_0:\C^2\to L^2(0,\pi)$ be the solution of the Sylvester equation
$AY_0-Y_0L_1=-\mathcal S$, that is, $Y_0e_\pm=(\lambda_\pm-A)^{-1}\mathcal Se_\pm$; it exists because $\sigma(A)\subset(-\infty,-1]$ does not meet
$\{\lambda_\pm\}$.

\begin{proposition}\label{prop:cascade}
Let $0<\alpha<1$ and suppose that $Y_0$ is injective, with smallest singular value $s_0>0$ and norm $y_0$. Then
$L$ is closed and densely defined, $\|e^{tL}\|\ge e^{t\cos\theta}$, and $W(L)$ meets $\Sa$ for every $\mu\ge0$. The Riesz
projection $P(a,v)=(a,\mu Y_0a)$ onto $\{\lambda_\pm\}$ is a slow component of rank $2$ with fast part $L_Q=A$, and for every
$\mu\ge0$:
\begin{enumerate}[label=\textup{(\alph*)}]
\item $\|P\|=\|Q\|=(1+\mu^2y_0^2)^{\frac{1}{2}}$, $C_Q=\|Q\|$ and $\sup_t\|E_\alpha(t^\alpha L_Q)\|=1$;
\item $\|P\|\le\Kr_\alpha(L;P)\le\|P\|\,\Kr_\alpha(L_1)$, and $\Kr_\alpha(L_P)\le\frac{y_0}{s_0}\,\Kr_\alpha(L_1)$, a bound independent of $\mu$;
\item $c_1\mu\le\Kr_\alpha(L)\le\sup_t\|E_\alpha(t^\alpha L)\|\le\kappa'_{\alpha,2}\,\Kr_\alpha(L_1)\|P\|+\|P\|$, with
$c_1:=\|(1-A)^{-1}\mathcal S(1-L_1)^{-1}\|>0$.
\end{enumerate}
In particular all of $\Kr_\alpha(L)$, $\Kr_\alpha(L;P)$, $C_Q$ and $\sup_t\|E_\alpha(t^\alpha L)\|$ grow linearly in $\mu$, while
$\Kr_\alpha(L_P)$ and $\sup_t\|E_\alpha(t^\alpha L_Q)\|$ stay bounded, so no estimate of the form
$\sup_t\|E_\alpha(t^\alpha L)\|\le C(\alpha,m)\big(\Kr_\alpha(L_P)+\sup_t\|E_\alpha(t^\alpha L_Q)\|\big)$ can hold.
\end{proposition}

\begin{proof}
$L$ is a bounded perturbation of $L_1\oplus A$, hence closed, densely defined and the generator of a $C_0$-semigroup, and the
eigenvector $(e_+,\mu Y_0e_+)$ for $\lambda_+$ gives $\|e^{tL}\|\ge e^{t\cos\theta}$. For unit $a\in\C^2$,
$\langle L(a,0),(a,0)\rangle=\langle L_1a,a\rangle$, and $W(L_1)$ is the segment $[\lambda_-,\lambda_+]$, which contains
$\cos\theta\in\Sa$. With $Y:=\mu Y_0$, one checks directly that $T:a\mapsto(a,Ya)$ maps $\C^2$ onto the invariant subspace
for $\lambda_\pm$, that $\{0\}\oplus L^2(0,\pi)$ is invariant, and that $P=TV$ with $V(a,v):=a$, so that $P$ is the Riesz
projection, $L_PT=TL_1$ and $L_Q=A$ on $QH=\{0\}\oplus L^2(0,\pi)$.

(a) In the decomposition $H=\C^2\oplus L^2$, $P=\begin{psmallmatrix}I&0\\ Y&0\end{psmallmatrix}$ and
$Q=\begin{psmallmatrix}0&0\\ -Y&I\end{psmallmatrix}$, whose norms both equal $(1+\|Y\|^2)^{\frac{1}{2}}$. Since $A$ is self-adjoint and
negative, $\|E_\alpha(t^\alpha A)\|\le1$ with equality at $t=0$, and $E_\alpha(t^\alpha L)Q=\begin{psmallmatrix}0&0\\ -E_\alpha(t^\alpha A)Y&E_\alpha(t^\alpha A)\end{psmallmatrix}$
has norm at most $\|Q\|$, with equality at $t=0$.

(b) The lower bound is \cref{lem:slow-basic}(c). Since $(z-L)^{-1}P=T(z-L_1)^{-1}V$, $\|V\|=1$ and $\|T\|=\|P\|$, the upper
bound follows. On $PH$, with the norm of $H$, $L_P=TL_1T^{-1}$, where $\|T\|=(1+\mu^2y_0^2)^{\frac{1}{2}}$ and
$\|T^{-1}\|=(1+\mu^2s_0^2)^{-\frac{1}{2}}$; hence $\Kr_\alpha(L_P)\le\|T\|\,\|T^{-1}\|\,\Kr_\alpha(L_1)$, and
$\frac{1+\mu^2y_0^2}{1+\mu^2s_0^2}\le\frac{y_0^2}{s_0^2}$ because $y_0\ge s_0$.

(c) Since $a_\alpha(1)=1$ and $(1-L)^{-1}(a,0)=\big((1-L_1)^{-1}a,\ \mu(1-A)^{-1}\mathcal S(1-L_1)^{-1}a\big)$, we get
$\Kr_\alpha(L)\ge\|(1-L)^{-1}\|\ge c_1\mu$, and $c_1>0$ because $(1-A)^{-1}$ is injective and $\mathcal S\ne0$. The middle inequality is
\cref{thm:kreiss-lower}, and the upper bound is \cref{thm:slow}(ii) with $\kappa'_{\alpha,2}$, since $0\notin\sigma(L_P)$, together with
(a) and (b).
\end{proof}

The injectivity of $Y_0$ is a nondegeneracy condition on the source profiles. For the profiles
$b_1=\pi^{-\frac{1}{2}}\mathbf 1_{(0,\pi)}$ and $b_2$ the normalised parabola $x(\pi-x)$ used in \cref{sec:exp-cascade}, a computation
with $2\cdot10^4$ sine modes gives $y_0\approx0.78$ and $s_0\approx0.030$ for $\theta=\frac{\pi}{3}$, so the bound of (b) on $\Kr_\alpha(L_P)$ is
about $26\,\Kr_\alpha(L_1)$, uniformly in $\mu$. For eigenfunction profiles $b_1=\sqrt{\frac{2}{\pi}}\sin x$ and
$b_2=\sqrt{\frac{2}{\pi}}\sin2x$ injectivity follows from an explicit $2\times2$ determinant.

\section{Optimality of the rank dependence}\label{sec:optimal}

We now show that the rank cannot be removed from \Cref{thm:A}, even on Hilbert spaces, and that in infinite
dimensions a finite fractional Kreiss constant implies neither boundedness of the family nor even its
existence. All examples live on one fixed Hilbert space of functions on the circle, and they are truncations of
one natural operator, the Weyl fractional derivative.

On the boundary rays of $\Sa$ the Mittag--Leffler function behaves like an exponential of an imaginary argument:
$E_\alpha((is)^\alpha)=\alpha^{-1}e^{is}+b_\alpha(s)$, with an error $b_\alpha$ that is bounded and of bounded variation
(see \eqref{eq:bstar} below). Consequently, if $A$ is a Fourier multiplier with symbol $ik$, so that $e^{tA}$ is the translation
$f\mapsto f(\cdot+t)$, and $B$ is the multiplier with symbol $(ik)^\alpha$, whose spectrum lies on $\partial\Sa$, then
\[
E_\alpha(t^\alpha B)=\alpha^{-1}e^{tA}+\big(\text{multiplier with symbol }b_\alpha(tk)\big).
\]
On $L^2(\T)$ this is harmless, because translations are unitary, but on a weighted space translation may amplify: if
the weight is small near $x=0$ and large near $x=\pm\pi$, a function concentrated near $0$, such as the Dirichlet kernel,
has a small norm, while its translate by $\pi$ is concentrated where the weight is large and has a large norm.
Translation by $\pi$ thus amplifies, and this amplification survives in $E_\alpha(t^\alpha B)$ the other hand, the
resolvent of $B$ is a multiplier whose symbol $(z-(ik)^\alpha)^{-1}$ varies slowly along the boundary rays, and on
weighted spaces with weights of Muckenhoupt class $A_2$ such multipliers are controlled by the Riesz projections,
uniformly in the number of frequencies. The Kreiss constant stays bounded while the propagator grows, and the
growth is limited only by the weight, which must remain in $A_2$. For $\alpha=1$ this is the mechanism of 
\citet{EisnerZwart2006}.


We first identify $\T=\R/2\pi\Z$ with $(-\pi,\pi]$, write $e_k(x):=e^{ikx}$ for $k\in\Z$, and let $\mathcal P$ be the space of
trigonometric polynomials. For $0<\gamma<1$ let
\[
w_\gamma(x):=|\tan(\tfrac{x}{2})|^\gamma,\qquad L^2_\gamma:=L^2(\T,w_\gamma\,dx),\qquad\|f\|_\gamma^2:=\int_{-\pi}^{\pi}|f(x)|^2w_\gamma(x)\,dx .
\]
The weight behaves like $(\frac{|x|}{2})^\gamma$ near $0$ and like $(\frac{2}{\pi-|x|})^\gamma$ near $\pm\pi$. Since $\gamma<1$, both $w_\gamma$
and $w_\gamma^{-1}$ are integrable, so $\mathcal P\subset L^2_\gamma\subset L^1(\T)$ by the Cauchy--Schwarz inequality, Fourier
coefficients $\hat f(k)$ of $f\in L^2_\gamma$ are defined and depend continuously on $f$, and $\mathcal P$ is dense in
$L^2_\gamma$. For $j\in\Z$ let $R_j$ be the Riesz projection $\sum_kc_ke_k\mapsto\sum_{k\ge j}c_ke_k$ on $\mathcal P$, and put
\[
r_\gamma:=\sup_{j\in\Z}\ \sup\big\{\|R_jf\|_\gamma:\ f\in\mathcal P,\ \|f\|_\gamma\le1\big\}.
\]

For a sequence $m=(m_k)_{k\in\Z}$ of bounded variation, $\Var_\Z(m):=\sum_{k\in\Z}|m_k-m_{k-1}|<\infty$, the limits
$m_{\pm\infty}$ exist, and we let $T_m$ denote the Fourier multiplier $T_me_k:=m_ke_k$ on $\mathcal P$. The following lemma
contains everything we need about multipliers on $L^2_\gamma$.

\begin{lemma}[Weighted multipliers]\label{lem:multiplier}
Let $0<\gamma<1$.
\begin{enumerate}[label=\textup{(\alph*)}]
\item $w_\gamma\in A_2(\T)$ and $r_\gamma<\infty$. Consequently every $R_j$ extends to a bounded projection on $L^2_\gamma$, and so
does the Fourier partial sum projection $S_N:=R_{-N}-R_{N+1}$ onto $V_N:=\operatorname{span}\{e_k:|k|\le N\}$, with
$\|S_N\|\le2r_\gamma$.
\item If $\Var_\Z(m)<\infty$, then $\|T_mf\|_\gamma\le\big(|m_{-\infty}|+r_\gamma\Var_\Z(m)\big)\|f\|_\gamma$ for $f\in\mathcal P$, so $T_m$ extends to a
bounded operator on $L^2_\gamma$ with $\widehat{T_mf}(k)=m_k\hat f(k)$.
\item If $m_k=g(\zeta_k)$, where the points $\zeta_k$ lie in the order of $k$ along a locally rectifiable path $\Gamma$, possibly
unbounded, and $g$ is holomorphic in a neighbourhood of $\Gamma$, then $\Var_\Z(m)\le\int_\Gamma|g'(\zeta)|\,|d\zeta|$.
\end{enumerate}
\end{lemma}

\begin{proof}
(a) Near $0$ and near $\pm\pi$ the weight is comparable to the power weights $|x|^\gamma$ and $|x\mp\pi|^{-\gamma}$, which belong
to $A_2$ since $\gamma<1$, and the integrability of $w_\gamma^{\pm1}$ controls long arcs, so $w_\gamma\in A_2(\T)$. The theorem of
\citet{HuntMuckenhouptWheeden1973} gives the boundedness of the conjugate function $f\mapsto\tilde f$ on $L^2_\gamma$, hence of
$R_0f=\frac12(f+i\tilde f)+\frac12\hat f(0)e_0$; since multiplication by $e_j$ is an isometry of $L^2_\gamma$ and $R_j=e_jR_0e_{-j}$, the bound is uniform
in $j$. The details are routine and given in \cref{app:optimal}.

(b) Summation by parts gives $T_mf=m_{-\infty}f+\sum_{j\in\Z}(m_j-m_{j-1})R_jf$ for $f\in\mathcal P$, a series that converges absolutely because $\Var_\Z(m)<\infty$ and $\|R_jf\|_\gamma\le r_\gamma\|f\|_\gamma$. To verify the
identity, note that $R_jf=f$ for $j\le\min\{k:\hat f(k)\ne0\}$ and $R_jf=0$ for $j>\max\{k:\hat f(k)\ne0\}$, so the terms with small $j$
telescope. This gives the bound; the formula for the Fourier coefficients follows by continuity.

(c) $|g(\zeta_j)-g(\zeta_{j-1})|\le\int_{\Gamma_j}|g'|\,|d\zeta|$, where $\Gamma_j$ is the part of $\Gamma$ between $\zeta_{j-1}$ and $\zeta_j$.
\end{proof}

The symbols we need are evaluations of resolvents along the boundary of $\Sa$, and their variation is controlled by an
elementary estimate (\cref{app:optimal}): for a ray $\mathsf R_\theta:=\{se^{i\theta}:s\ge0\}$ and $z\notin\mathsf R_\theta$,
\begin{equation}\label{eq:ray}
\int_0^\infty\frac{ds}{|z-se^{i\theta}|^2}\le\frac{\pi}{\dist(z,\mathsf R_\theta)} .
\end{equation}


Next, the construction needs two quantitative facts, one about the Mittag--Leffler function and one about the weight. Both
are direct computations, and we state them here and give the details in \cref{app:optimal}.

The first one makes the heuristic of the previous subsection precise. Let $\Phi_\alpha(s):=E_\alpha((is)^\alpha)$ for $s\in\R$ and
$b_\alpha(s):=\Phi_\alpha(s)-\alpha^{-1}e^{is}$. Then $b_1\equiv0$, so that $b^*_1=0$ in the notation below, and for $0<\alpha<1$ the
function $b_\alpha$ is continuous with
\begin{equation}\label{eq:bstar}
b^*_\alpha:=\sup_{\R}|b_\alpha|+\Var_\R(b_\alpha)<\infty .
\end{equation}
For large $|s|$ this follows from \cref{lem:ML}(c), which gives $|b_\alpha(s)|=O(|s|^{-\alpha})$ and $|b_\alpha'(s)|=O(|s|^{-1-\alpha})$ on
the boundary ray, and for small $|s|$ from the power series.

The second one quantifies the amplification produced by translation. We use the Dirichlet kernel
\[
D_N:=\sum_{|k|\le N}e_k,\qquad D_N(x)=\frac{\sin((N+\frac12)x)}{\sin(\frac{x}{2})} .
\] Since $|D_N|$ is of size $N$ on an interval of length about $\frac{1}{N}$ around $0$ and bounded by $|\sin(\frac{x}{2})|^{-1}$ elsewhere, one
finds $\|D_N\|_\gamma^2\asymp N^{1-\gamma}$, whereas the translate places this peak near $\pm\pi$, where the weight is of size $N^\gamma$, so
that $\|D_N(\cdot+\pi)\|_\gamma^2\gtrsim N^{1+\gamma}$. Estimating the two weighted integrals
gives, for all $N\in\N$,
\begin{equation}\label{eq:dirichlet}
\|D_N(\cdot+\pi)\|_\gamma\ge c_\gamma N^\gamma\|D_N\|_\gamma,\qquad
c_\gamma:=\Big(\frac{2\cdot3^{2\gamma-2}}{(1-\gamma)Q_\gamma}\Big)^{\frac{1}{2}},
\end{equation}
where $Q_\gamma:=\frac{2}{1+\gamma}+\frac{2\pi^2}{1-\gamma}+4\int_{\frac{\pi}{2}}^{\pi}\tan^\gamma(\frac{y}{2})\,dy$.


\medskip

Finally, we proceed to the construction. Fix $\alpha\in(0,1]$ and $\gamma\in(0,1)$. The \emph{Weyl fractional derivative} of order $\alpha$ on $L^2_\gamma$ is the Fourier
multiplier with symbol $(ik)^\alpha=|k|^\alpha e^{i\operatorname{sgn}(k)\frac{\alpha\pi}{2}}$,
\[
\D^\alpha f:=g\quad\text{whenever } f,g\in L^2_\gamma \text{ and }\hat g(k)=(ik)^\alpha\hat f(k)\text{ for all }k\in\Z,
\]
a closed, densely defined operator, since the Fourier coefficients are continuous on $L^2_\gamma$ and $\mathcal P\subset D(\D^\alpha)$.
Its eigenvalues $(ik)^\alpha$ lie on the two boundary rays of $\Sa$, and $\Sa\subset\rho(\D^\alpha)$ by \cref{prop:weyl} below. For $N\in\N$
we let
\[
B_N:=N^{-\alpha}\D^\alpha S_N,\qquad
\ell_N(k):=\begin{cases}(\frac{ik}{N})^\alpha,&|k|\le N,\\ 0,&|k|>N,\end{cases}
\]
so that $B_N$ is the multiplier with symbol $\ell_N$.
It is a bounded operator on $L^2_\gamma$ of rank $2N$, since $\ell_N(0)=0$. 
\begin{theorem}\label{thm:optimal}
Let $0<\alpha\le1$ and $0<\gamma<1$. For every $N\in\N$, the operator $B_N\in\B(L^2_\gamma)$ has rank $2N$ and satisfies
\begin{enumerate}[label=\textup{(\roman*)}]
\item $\sigma(B_N)=\{\ell_N(k):|k|\le N\}\subset\partial\Sa$ and $\|B_N\|\le4r_\gamma$;
\item $\Kr_\alpha(B_N)\le\frac{1+4\pi r_\gamma}{\alpha}$;
\item $\|E_\alpha((\pi N)^\alpha B_N)\|\ge\alpha^{-1}c_\gamma N^\gamma-(1+r_\gamma)b^*_\alpha$.
\end{enumerate}
\end{theorem}

\begin{proof}
(i) $B_N$ is diagonal in the basis $(e_k)$ of $\mathcal P$ with eigenvalues $\ell_N(k)$, and $\ell_N(k)=0$ for $k=0$ and for
$|k|>N$. As $k$ increases from $-\infty$ to $\infty$, the points $\ell_N(k)$ stay at $0$, jump to $e^{-i\frac{\alpha\pi}{2}}$ at
$k=-N$, move along $\partial\Sa$ to $0$ at $k=0$ and on to $e^{i\frac{\alpha\pi}{2}}$ at $k=N$, and jump back to $0$. They lie in this
order on the closed path $\Gamma$ that runs along the segment $[0,e^{-i\frac{\alpha\pi}{2}}]$ and the segment $[0,e^{i\frac{\alpha\pi}{2}}]$
twice each, which has length $4$. \Cref{lem:multiplier}(b),(c) with $g(\zeta)=\zeta$ and $m_{-\infty}=0$ gives $\|B_N\|\le4r_\gamma$. If
$z\notin\{\ell_N(k):|k|\le N\}$, the symbol $(z-\ell_N(k))^{-1}$ takes finitely many values and is constant for $|k|>N$, so it has finite
variation, and by \cref{lem:multiplier}(b) the corresponding multiplier inverts $z-B_N$; this gives the spectrum.

(ii) Let $z\in\Sa$. Then $z-B_N$ is the multiplier with symbol $z-\ell_N(k)\ne0$, and its inverse is $T_m$ with
$m_k=g(\ell_N(k))$, $g(\zeta):=(z-\zeta)^{-1}$, which is holomorphic near $\Gamma\subset\partial\Sa$. Here
$|m_{-\infty}|=|z|^{-1}\le\dist(z,\partial\Sa)^{-1}$, and by \cref{lem:multiplier}(c) and \eqref{eq:ray}, each segment of $\Gamma$ being
traversed twice,
\[
\Var_\Z(m)\le\int_\Gamma\frac{|d\zeta|}{|z-\zeta|^2}\le\frac{2\pi}{\dist(z,\mathsf R_{\frac{\alpha\pi}{2}})}+\frac{2\pi}{\dist(z,\mathsf R_{-\frac{\alpha\pi}{2}})}
\le\frac{4\pi}{\dist(z,\partial\Sa)} .
\]
Hence $\|(z-B_N)^{-1}\|\le\frac{1+4\pi r_\gamma}{\dist(z,\partial\Sa)}$, and \cref{lem:elementary}(c) gives (ii).

(iii) Since $B_N$ is diagonal, $E_\alpha(t^\alpha B_N)$ is the multiplier with symbol $E_\alpha(t^\alpha\ell_N(k))$, and
$t^\alpha\ell_N(k)=(\frac{itk}{N})^\alpha$ for $|k|\le N$. On $V_N$ it therefore acts as the multiplier with symbol
$\Phi_\alpha(\frac{tk}{N})=\alpha^{-1}e^{\frac{itk}{N}}+b_\alpha(\frac{tk}{N})$. The first part maps $D_N$ to $\alpha^{-1}D_N(\cdot+\frac{t}{N})$. For the second, the
points $\frac{tk}{N}$, $|k|\le N$, increase with $k$, so the restriction of the symbol to $|k|\le N$, extended by constants outside,
has modulus at most $\sup|b_\alpha|$ and variation at most $\Var_\R(b_\alpha)$; by \cref{lem:multiplier}(a),(b) and \eqref{eq:bstar} it maps $D_N$ to a function of
norm at most $(1+r_\gamma)b^*_\alpha\|D_N\|_\gamma$. Taking $t=\pi N$ and using \eqref{eq:dirichlet},
$\|E_\alpha((\pi N)^\alpha B_N)D_N\|_\gamma\ge\big(\alpha^{-1}c_\gamma N^\gamma-(1+r_\gamma)b^*_\alpha\big)\|D_N\|_\gamma$.
\end{proof}

For $n\in\N$ let $\kappa^{\rm opt}_\alpha(n)$ denote the least constant $\kappa$ such that
$\sup_t\|E_\alpha(t^\alpha L)\|\le\kappa\,\Kr_\alpha(L)$ for all Hilbert spaces $H\ne\{0\}$ and all $L\in\B(H)$ with $\rank L\le n$. By
\cref{thm:finite-rank}, $\kappa^{\rm opt}_\alpha(n)\le\kappa_{\alpha,n}$.

\begin{corollary}\label{cor:optimal-n}
Let $\alpha\in(0,1]$ and $\gamma\in(0,1)$. For every $n\ge2$,
\[
\kappa^{\rm opt}_\alpha(n)\ge\frac{c_\gamma\lfloor \frac{n}{2}\rfloor^\gamma-\alpha(1+r_\gamma)b^*_\alpha}{1+4\pi r_\gamma},\qquad\text{in particular}\qquad
\liminf_{n\to\infty}\frac{\kappa^{\rm opt}_\alpha(n)}{n^\gamma}\ge\frac{2^{-\gamma}c_\gamma}{1+4\pi r_\gamma}>0.
\]
Hence $n^\gamma\lesssim\kappa^{\rm opt}_\alpha(n)\lesssim n$ for every $\gamma<1$.
\end{corollary}

\begin{proof}
Take $N:=\lfloor \frac{n}{2}\rfloor$, so that $\rank B_N\!=\!2N\!\le\! n$, and combine \cref{thm:optimal}(ii) and (iii).
\end{proof}

\begin{corollary}\label{cor:infinite}
Let $H:=\bigoplus_{N\ge1}L^2_\gamma$ (Hilbert direct sum) and $L:=\bigoplus_{N\ge1}B_N$. Then $L\in\B(H)$,
$\sigma(L)\cap\Sa=\emptyset$, $\Kr_\alpha(L)\le\frac{1+4\pi r_\gamma}{\alpha}$, and
\[
\limsup_{t\to\infty}t^{-\gamma}\|E_\alpha(t^\alpha L)\|\ge\alpha^{-1}\pi^{-\gamma}c_\gamma>0.
\]
In particular $\Kr_\alpha(L)<\infty$ does not imply $\sup_t\|E_\alpha(t^\alpha L)\|<\infty$ on Hilbert spaces.
\end{corollary}

\begin{proof}
By \cref{thm:optimal}(i), $\|L\|\le4r_\gamma$. For $z\in\Sa$ the inverses $(z-B_N)^{-1}$ are uniformly bounded by
$\frac{1+4\pi r_\gamma}{\dist(z,\partial\Sa)}$, as shown in the proof of \cref{thm:optimal}(ii), so their direct sum inverts $z-L$ and the bound
on $\Kr_\alpha(L)$ follows. Since $E_\alpha(t^\alpha L)=\bigoplus_NE_\alpha(t^\alpha B_N)$, the choice $t=\pi N$ in \cref{thm:optimal}(iii) gives
the claim.
\end{proof}

The untruncated operator shows that in infinite dimensions the fractional Kreiss condition does not even guarantee
the existence of the family.

\begin{proposition}\label{prop:weyl}
Let $0<\alpha\le1$ and $0<\gamma<1$. Then $\Sa\subset\rho(\D^\alpha)$ and $\Kr_\alpha(\D^\alpha)\le\frac{2\pi r_\gamma}{\alpha}$ on $L^2_\gamma$, but $\D^\alpha$
generates no Mittag--Leffler family on $L^2_\gamma$. For $\alpha=1$, the derivative $\D=\frac{d}{dx}$ satisfies the Kreiss resolvent
condition on $L^2_\gamma$ but generates no $C_0$-semigroup.
\end{proposition}

\begin{proof}
Let $z\in\Sa$ and $m_k:=(z-(ik)^\alpha)^{-1}$. The points $(ik)^\alpha$ run in the order of $k$ along the lower boundary ray from
infinity to $0$ and then along the upper one to infinity, so $m_{-\infty}=0$ and, by \cref{lem:multiplier}(c) and
\eqref{eq:ray}, $\Var_\Z(m)\le\frac{2\pi}{\dist(z,\partial\Sa)}$. The bounded operator $T_m$ maps $L^2_\gamma$ into $D(\D^\alpha)$, and comparing
Fourier coefficients gives $(z-\D^\alpha)T_m=I$ and $T_m(z-\D^\alpha)f=f$ for $f\in D(\D^\alpha)$. Hence $z\in\rho(\D^\alpha)$,
$\|(z-\D^\alpha)^{-1}\|\le\frac{2\pi r_\gamma}{\dist(z,\partial\Sa)}$, and \cref{lem:elementary}(c) gives the bound on $\Kr_\alpha(\D^\alpha)$.

Suppose that $\D^\alpha$ generates a Mittag--Leffler family $S$. For $f\in\mathcal P$ and $\lambda$ large,
\[
\int_0^\infty e^{-\lambda t}S(t)f\,dt=\sum_k\frac{\lambda^{\alpha-1}}{\lambda^\alpha-(ik)^\alpha}\,\hat f(k)e_k,
\]
and the scalar identity $\int_0^\infty e^{-\lambda t}E_\alpha(t^\alpha\zeta)\,dt=\lambda^{\alpha-1}(\lambda^\alpha-\zeta)^{-1}$ for $\zeta\in\La$ (proof of
\cref{prop:normal}) and uniqueness of the Laplace transform give
\[
S(t)f=\sum_k\Phi_\alpha(tk)\hat f(k)e_k=\alpha^{-1}f(\cdot+t)+T_{b(t)}f,\qquad b(t)_k:=b_\alpha(tk).
\]
By \cref{lem:multiplier}(b) and \eqref{eq:bstar}, $\|T_{b(t)}f\|_\gamma\le(1+r_\gamma)b^*_\alpha\|f\|_\gamma$, so
$\|f(\cdot+\pi)\|_\gamma\le\alpha\big(\|S(\pi)\|+(1+r_\gamma)b^*_\alpha\big)\|f\|_\gamma$ for all $f\in\mathcal P$, which contradicts \eqref{eq:dirichlet} with
$f=D_N$ and $N\to\infty$. For $\alpha=1$, $b_1\equiv0$ and $S(t)$ would be the translation by $t$.
\end{proof}

By \cref{lem:elementary}(e), $\Kr_\alpha(cL)=\Kr_\alpha(L)$ for $c>0$, so the absence of generation could also be produced by
rescaling the blocks of \cref{cor:infinite}; \cref{prop:weyl} shows that it happens for a natural differential
operator.

\begin{figure}[t]
\centering
\includegraphics[width=\textwidth]{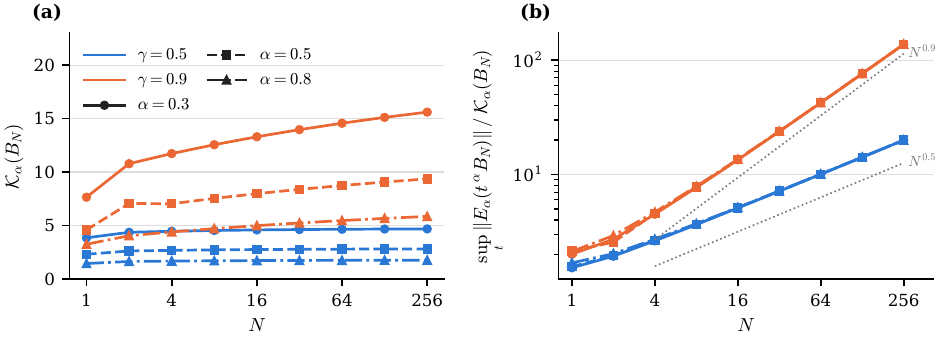}
\caption{The truncated fractional derivatives $B_N$ on $L^2_\gamma$ (\cref{thm:optimal}). (a) Computed $\Kr_\alpha(B_N)$; it is bounded by
\cref{thm:optimal}(ii), and for $\gamma=0.9$ the computed values are still increasing slowly at $N=256$. (b) The ratio
$\frac{\sup_t\|E_\alpha(t^\alpha B_N)\|}{\Kr_\alpha(B_N)}$ grows at a rate close to $N^\gamma$ (dotted guides); the curves for the three values of $\alpha$
nearly coincide.}\label{fig:optimal}
\end{figure}

The constants of \cref{thm:optimal} are far from sharp, and it is instructive to see how the construction behaves in
practice. \Cref{fig:optimal} shows the Kreiss constants $\Kr_\alpha(B_N)$ and the ratios of the two sides of
\Cref{thm:A}, computed on $V_N$ through the Gram matrix of the weight, for $N$ up to $256$, that is, for ranks up to $512$. The Kreiss constants stay between $1.4$
and $16$, while the ratios grow at a rate close to $N^\gamma$, with fitted exponents $0.49$ for $\gamma=0.5$ and $0.84$ for $\gamma=0.9$, and
reach $138$ at $N=256$. The curves for different $\alpha$ nearly coincide, because the translation term
$\alpha^{-1}D_N(\cdot+\frac{t}{N})$ in the proof of \cref{thm:optimal} dominates. The amplification is therefore a property of the weighted geometry, which the fractional clock transmits
essentially unchanged. The complete data are in Section S4 of the Supplementary Material.

Together with \cref{thm:finite-rank}, \cref{cor:optimal-n} leaves open only whether the exponent one is attained, which for
$\alpha=1$ is the case by the examples of \citet{Kraaijevanger1994}.

\section{Mechanisms of fractional transient growth}\label{sec:mechanisms}

The two-sided estimates of \Cref{thm:A} say how large the transient peak is; this section explains where it comes
from and when it occurs. Several of its features have no counterpart for semigroups: normal operators amplify, modes
that are unstable for the semigroup are stable for the family, and nilpotent non-normality, the classical source of
dimension-dependent growth, is tamed by the Mittag--Leffler distribution for each fixed $\alpha<1$.

\subsection{Normal amplification and fractionally stable modes}

\begin{figure}[t]
\centering
\includegraphics[width=\textwidth]{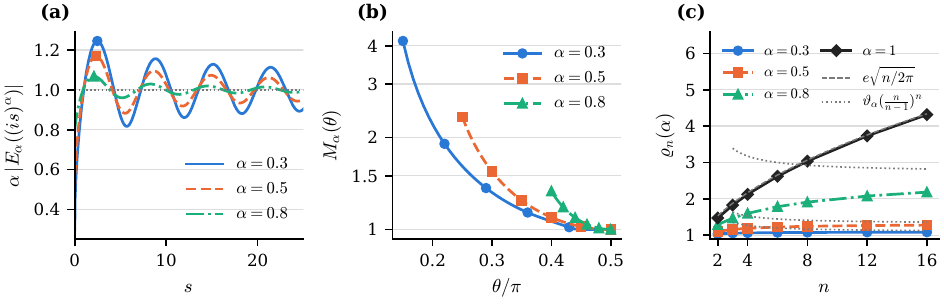}
\caption{(a) The function $s\mapsto\alpha|E_\alpha((is)^\alpha)|$ on the boundary ray of $\Sa$ overshoots $1$ before it settles into the
oscillation $\alpha^{-1}e^{is}$ (\cref{lem:mstar}). (b) The amplification $M_\alpha(\theta)$ of normal operators with spectrum on the ray of
argument $\theta$ decreases from $\Mst$ at $\theta=\frac{\alpha\pi}{2}$ to $1$ at $\theta=\frac{\pi}{2}$. (c) The Jordan-block limits $\varrho_n(\alpha)$ of
\cref{prop:jordan}: they grow like $e\sqrt{\frac{n}{2\pi}}$ for $\alpha=1$ and stay below $\vartheta_\alpha(\frac{n}{n-1})^n$ for $\alpha<1$.}\label{fig:constants}
\end{figure}

By \cref{prop:normal}, a normal generator with spectrum on a boundary ray of $\Sa$ amplifies by exactly $\Mst$, while its
Kreiss constant equals $\frac{1}{\alpha}$, so the ratio of the two sides of \Cref{thm:A} equals $\alpha\Mst$. This ratio decreases from
$1.33$ at $\alpha=0.1$ to $1$ at $\alpha=1$ (\cref{tab:constants}). The amplification originates on the boundary ray, where
$E_\alpha(re^{i\frac{\alpha\pi}{2}})$ is the sum of the unimodular oscillation $\alpha^{-1}e^{ir^{\frac{1}{\alpha}}}$ and a slowly decaying correction of
order $r^{-1}$, and the two add constructively at moderate $r$ (\cref{fig:constants}(a)). If the spectrum moves away from
the boundary ray into the interior of $\La$, the amplification decreases monotonically to $1$, which is reached at the
imaginary axis (\cref{lem:mstar} and \cref{fig:constants}(b)).

\begin{table}[t]
\centering\small
\caption{Constants of the theory for selected orders. $\Mst$ is the sharp bound for normal operators, attained at $s_*=r_*^{\frac{1}{\alpha}}$ on the
boundary ray; $\alpha\Mst$ is the ratio of the two sides of \Cref{thm:A} for normal operators with spectrum on $\partial\Sa$; $\varrho_2(\alpha)$ is the
Jordan-block limit of \cref{prop:jordan}; $\kappa_{\alpha,1}$ and $\kappa'_{\alpha,1}$ are the constants of \Cref{thm:A} for rank one. Further orders,
including irrational ones, are listed in Table~S1 of the Supplementary Material.}\label{tab:constants}
\begin{tabular}{@{}lcccccc@{}}
\toprule
$\alpha$ & $\Mst$ & $s_*=r_*^{1/\alpha}$ & $\alpha\Mst$ & $\varrho_2(\alpha)$ & $\kappa_{\alpha,1}$ & $\kappa'_{\alpha,1}$\\
\midrule
0.10 & 13.2608 & 2.5966 & 1.3261 & 1.0041 & 14.41 & 11.69\\
0.30 & 4.1496 & 2.4670 & 1.2449 & 1.0376 & 13.32 & 10.60\\
0.50 & 2.3413 & 2.3272 & 1.1707 & 1.1069 & 12.23 & 9.51\\
0.80 & 1.3347 & 2.0956 & 1.0677 & 1.2885 & 10.60 & 7.88\\
0.90 & 1.1491 & 2.0122 & 1.0342 & 1.3731 & 10.06 & 7.34\\
\bottomrule
\end{tabular}

\end{table}

The stability region of the family contains the sector $\{\frac{\alpha\pi}{2}\le|\Arg z|<\frac{\pi}{2}\}$ of the right half-plane, where the
semigroup grows exponentially. The simplest example is the oscillator $L=\begin{psmallmatrix}\cos\theta&-\sin\theta\\ \sin\theta&\cos\theta\end{psmallmatrix}$ with
$\frac{\alpha\pi}{2}<\theta<\frac{\pi}{2}$, whose eigenvalues $e^{\pm i\theta}$ have positive real part. Then $\|e^{tL}\|=e^{t\cos\theta}$ and
$\Kr_1(L)=\infty$, while $\Kr_\alpha(L)\le\frac{1}{\alpha}$ and $\sup_t\|E_\alpha(t^\alpha L)\|=M_\alpha(\theta)\le\Mst$, since $L$ is normal. Adding a
non-normal coupling changes the size of both sides of \Cref{thm:A}, but not their ratio by much; this is prototype (a), analysed
in \cref{prop:cascade} and computed in \cref{sec:exp-cascade}. Note that in this situation the numerical range of $L$ contains $\cos\theta\in\Sa$, so neither
\cref{thm:numrange} nor any of the spectral-set bounds discussed after it applies.

\subsection{Jordan blocks}

Let $Z_n\in\C^{n\times n}$ be the nilpotent shift, with ones on the superdiagonal, and $J_{n,\mu}:=-I+\mu Z_n$ with $\mu>0$,
acting on $\C^n$ with the Euclidean norm. For $\alpha=1$ these blocks are the standard source of the dimension dependence in
\eqref{eq:classical}: as $\mu\to\infty$, the ratio $\frac{\sup_t\|e^{tJ_{n,\mu}}\|}{\Kr_1(J_{n,\mu})}$ tends to $\frac{n^ne^{1-n}}{(n-1)!}\sim e\sqrt{\frac{n}{2\pi}}$. Put
\[
\Delta_n(\alpha):=\sup_{x\ge0}\frac{x^{n-1}|E_\alpha^{(n-1)}(-x)|}{(n-1)!},\qquad\varrho_n(\alpha):=\frac{n^n}{(n-1)^{n-1}}\,\Delta_n(\alpha).
\]

\begin{proposition}\label{prop:jordan}
Let $n\ge2$ and $0<\alpha\le1$.
\begin{enumerate}[label=\textup{(\alph*)}]
\item $\sup_{z\in\Sa}a_\alpha(z)|1+z|^{-n}=(n-1)^{n-1}n^{-n}$, and
$\lim_{\mu\to\infty}\frac{\sup_t\|E_\alpha(t^\alpha J_{n,\mu})\|}{\Kr_\alpha(J_{n,\mu})}=\varrho_n(\alpha)\ge1$.
\item $\varrho_n(1)=\frac{n^ne^{1-n}}{(n-1)!}$; in particular $\varrho_n(1)\sim e\sqrt{\frac{n}{2\pi}}$ as $n\to\infty$.
\item If $0<\alpha<1$, then $\varrho_n(\alpha)\le\vartheta_\alpha\big(\frac n{n-1}\big)^n\le4\vartheta_\alpha$, with $\vartheta_\alpha$ as in
\cref{lem:ML}(e).
\end{enumerate}
\end{proposition}

\begin{proof}
(a) For $z\in\Sa$ we have $a_\alpha(z)\le\RE z$ and $|1+z|\ge1+\RE z$, so $a_\alpha(z)|1+z|^{-n}\le x(1+x)^{-n}$ with $x=\RE z>0$. The
right-hand side is at most $(n-1)^{n-1}n^{-n}$, with equality at $x=\frac{1}{n-1}$, and $z=\frac{1}{n-1}$ attains this value. Let
$\tau:=t^\alpha$. Since $Z_n^n=0$ and $-\tau I$ commutes with $Z_n$, Taylor's formula gives
$E_\alpha(\tau J_{n,\mu})=\sum_{k=0}^{n-1}\frac{E_\alpha^{(k)}(-\tau)(\tau\mu)^kZ_n^k}{k!}$. For $\alpha<1$, differentiating
\cref{lem:ML}(d) under the integral sign gives $E_\alpha^{(k)}(-\tau)=\int_0^\infty u^ke^{-\tau u}m_\alpha(u)\,du\ge0$, hence
\begin{equation}\label{eq:jordan-coef}
\frac{\tau^kE_\alpha^{(k)}(-\tau)}{k!}=\int_0^\infty\frac{(\tau u)^ke^{-\tau u}}{k!}\,m_\alpha(u)\,du\le\sup_{s\ge0}\frac{s^ke^{-s}}{k!}\le1,
\end{equation}
which is trivial for $\alpha=1$. Consequently, uniformly in $\tau$,
\[
\Big\|\mu^{1-n}E_\alpha(\tau J_{n,\mu})-\frac{\tau^{n-1}E_\alpha^{(n-1)}(-\tau)}{(n-1)!}Z_n^{n-1}\Big\|\le\sum_{k<n-1}\mu^{k-n+1},
\]
and
$\mu^{1-n}\sup_t\|E_\alpha(t^\alpha J_{n,\mu})\|\to\Delta_n(\alpha)$ because $\|Z_n^{n-1}\|=1$. Similarly
$(z-J_{n,\mu})^{-1}=\sum_{k=0}^{n-1}\mu^kZ_n^k(z+1)^{-k-1}$, and since $a_\alpha(z)\le\RE z<|1+z|$ and $|1+z|\ge1$ on $\Sa$, the same argument
gives $\mu^{1-n}\Kr_\alpha(J_{n,\mu})\to(n-1)^{n-1}n^{-n}$. The ratio is at least $1$ by \cref{thm:kreiss-lower}.

(b) $\Delta_n(1)=\sup_x\frac{x^{n-1}e^{-x}}{(n-1)!}=\frac{(n-1)^{n-1}e^{1-n}}{(n-1)!}$, and Stirling's formula gives the asymptotics.

(c) Substituting $s=xu$ in \eqref{eq:jordan-coef},
\[
\frac{x^{n-1}E_\alpha^{(n-1)}(-x)}{(n-1)!}=\int_0^\infty\frac{s^{n-2}e^{-s}}{(n-1)!}\cdot\frac sx\,m_\alpha\Big(\frac sx\Big)ds\le\vartheta_\alpha\frac{(n-2)!}{(n-1)!}=\frac{\vartheta_\alpha}{n-1},
\]
so $\varrho_n(\alpha)\le\frac{\vartheta_\alpha n^n}{(n-1)^n}$, and $(\frac{n}{n-1})^n\le4$ for $n\ge2$.
\end{proof}

In the limit of strong non-normality, Jordan blocks thus produce only a bounded amplification ratio when $\alpha<1$, for instance
$\varrho_n(\frac{1}{2})\le\sqrt{\frac{2}{e\pi}}\,(\frac{n}{n-1})^n\to\sqrt{\frac{2e}{\pi}}\approx1.316$. The reason is the averaging in
\cref{lem:families}(b): the factor $\frac{s^{n-1}e^{-s}}{(n-1)!}$, whose sharp peak of height of order $n^{-\frac{1}{2}}$ produces the growth when
$\alpha=1$, is smeared out by the Mittag--Leffler distribution. The bound is not uniform in $\alpha$, since $\vartheta_\alpha\to\infty$ as
$\alpha\uparrow1$: the moments $\mathbb EY_\alpha^k=\frac{k!}{\Gamma(1+\alpha k)}$ tend to $1$, so $Y_\alpha\to1$ in distribution, which is incompatible with a
uniform bound on $u\,m_\alpha(u)$ near $u=1$. The values in \cref{fig:constants}(c), listed in Table~S2 of the Supplementary Material, show the transition: $\varrho_{16}$
equals $1.08$, $1.28$, $2.18$ and $4.32$ for $\alpha=0.3$, $0.5$, $0.8$ and $1$. These limits concern Jordan blocks in the regime of strong
non-normality only, so they do not exclude other nilpotent mechanisms, but they are consistent with the construction of
\Cref{thm:C}, where the growth is produced by oscillatory boundary modes; \cref{sec:exp-markov} shows the same taming in the
means of time-changed compartment chains.

\section{The prototypes revisited}\label{sec:experiments}

We return to the three prototypes of \Cref{sec:nutshell} and compare, for each of them, the propagator
$\sup_t\|E_\alpha(t^\alpha L)\|$ with the two sides of \Cref{thm:A} and with the bounds that are otherwise available. These are the
semigroup route (B1), $\sup_t\|E_\alpha(t^\alpha L)\|\le\sup_s\|e^{sL}\|$ by \cref{lem:families}(b), which is infinite whenever $\sigma(L)$ meets
$\HH$ or $L$ generates no semigroup; the numerical-range bounds (B2), $K\sup_{C_W}|E_\alpha|$ with $K=1+\sqrt2$
\citep{CrouzeixPalencia2017}, $K=K(\beta)$ \citep{Crouzeix2025} or $K=1$ (\cref{thm:numrange}), all infinite when $W(L)$ meets $\Sa$; and
crude deflation (B3), which estimates the slow block and the obliqueness of the splitting separately, as in \eqref{eq:crude}. The
Supplementary Material also reports (B4), \Cref{thm:A} with the dimension in place of the rank, and an eigenvector bound (B5) is used
where a complete eigenbasis is available.
Quantities attached to slow components are computed through the reductions of \cref{lem:struct,lem:slow-basic}, at a cost governed
by $m$, and Kreiss constants are maximised over polar grids in $\Sa$ with local refinement, so the reported values are accurate lower
estimates rather than certified upper bounds. The numerical methods, the complete tables and the code are described in the
Supplementary Material, which also reports a statistical study of $16200$ random operators of rank at most $32$ on spaces of
dimension up to $10^5$: there the ratio of the two sides of \Cref{thm:A} never exceeds $5.6$, nor $17\,\%$ of the corresponding
$\kappa_{\alpha,n}$, and is essentially independent of the dimension for $d\ge10^3$ (Section~S3).

\subsection{Prototype (a): the oscillator-driven subdiffusion}\label{sec:exp-cascade}

For the generator of \cref{prop:cascade} with $\alpha=\frac{1}{2}$, $\theta=\frac{\pi}{3}$ and the profiles $b_1=\pi^{-\frac{1}{2}}\mathbf 1_{(0,\pi)}$ and $b_2$ the
normalised parabola $x(\pi-x)$, all quantities are computed in the sine basis, with up to $10^4$ modes, at which the reported digits no
longer change. By \cref{prop:cascade} the semigroup bound (B1) and every numerical-range bound (B2) are infinite for every coupling
strength, so among the bounds that do not require an eigenbasis only the slow component gives information. For (B3) we evaluate \eqref{eq:crude} with the constant
$\kappa'_{\frac{1}{2},2}$, computing $\Kr_\alpha(L_P)$ in the norm of $H$ through the similarity $L_P=TL_1T^{-1}$ of the proof of \cref{prop:cascade}.
Since this operator has a complete eigenbasis, we also report the eigenvector bound (B5), $\operatorname{cond}(S)\max_{\lambda\in\sigma(L)}\sup_t|E_\alpha(t^\alpha\lambda)|$,
where $S$ maps the standard basis to an eigenbasis with normalised vectors.

\begin{figure}[t]
\centering
\includegraphics[width=\textwidth]{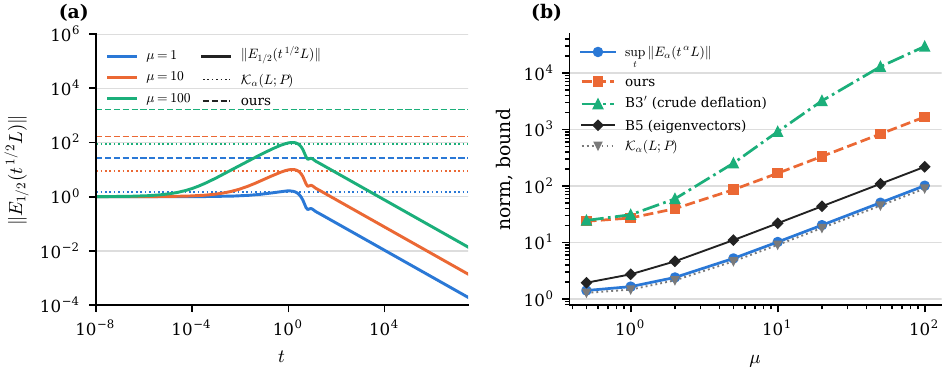}
\caption{Prototype (a) with $\alpha=\frac{1}{2}$ and $\theta=\frac{\pi}{3}$. (a) The envelope $t\mapsto\|E_{\frac{1}{2}}(t^{\frac{1}{2}}L)\|$ for three coupling strengths: initial
slope, transient peak near $t=1$ and algebraic decay of order $t^{-\frac{1}{2}}$; dotted lines mark $\Kr_{\frac{1}{2}}(L;P)$ and dashed lines the upper bound of
\Cref{thm:A}. (b) Propagator, deflated Kreiss constant and upper bounds as functions of $\mu$: ``ours'' is the upper bound of
\Cref{thm:A}, B3$'$ the crude deflation bound (B3) and B5 the eigenvector bound; the semigroup and numerical-range bounds are
infinite.}\label{fig:cascade}
\end{figure}

\Cref{fig:cascade} shows the result. The envelope rises with the initial slopes of \cref{prop:timescales}(a), peaks near $t=1$ and
decays like $t^{-\frac{1}{2}}$, as \cref{prop:timescales}(d) predicts, since $L$ satisfies its hypotheses with $\psi=\theta$. The four quantities $\Kr_\alpha(L)$,
$\Kr_\alpha(L;P)$, $C_Q=\|Q\|$ and $\sup_t\|E_\alpha(t^\alpha L)\|$ grow linearly in $\mu$, while $\Kr_\alpha(L_1)=1.229$ and the bound $1$ on the fast
part do not move; at $\mu=100$, for instance, $\|P\|=77.7$, $\Kr_\alpha(L)=67.3$, $\Kr_\alpha(L;P)=89.7$ and $\sup_t\|E_\alpha(t^\alpha L)\|=100.9$. The
propagator exceeds the deflated Kreiss constant by at most $13\,\%$ for every coupling strength, and by at most $12\,\%$ and $4\,\%$ for
$\alpha=0.3$ and $\alpha=0.8$, so in this model $\Kr_\alpha(L;P)$ is a faithful indicator of the transient although the semigroup grows
exponentially and the numerical range meets the sector of instability. The lower bound $\Kr_\alpha(L)$ is within a factor $1.5$ of the
propagator, and the upper bound of \Cref{thm:A}, $1662$ at $\mu=100$, is $16$ to $17$ times larger, a gap that reflects the worst-case
constant $\kappa'_{\frac{1}{2},2}=17.7$ rather than this operator; the quotients of \cref{thm:bernstein} for concrete rational functions stay
below $40\,\%$ of their bound (Table~S3). The crude bound (B3) is much larger: $\Kr_\alpha(L_P)$ grows to $21.4$ at
$\mu=100$, below the bound $\frac{y_0}{s_0}\Kr_\alpha(L_1)\approx32$ of \cref{prop:cascade}(b), and \eqref{eq:crude} multiplies it once more by $\|P\|$, so that
the obliqueness of the splitting is counted twice and (B3) exceeds our bound by a factor of $18$ at $\mu=100$. The deflated constant counts it
once, since $\Kr_\alpha(L;P)\le\|P\|\Kr_\alpha(L_1)$. The eigenvector bound (B5) is about twice the propagator and is the sharpest of all
bounds for this operator, but it needs a complete and well-conditioned eigenbasis of the whole generator, which is not available for
non-diagonalisable or nearly defective generators, whereas the slow-component bound uses only the slow block, the splitting and a
bound on the fast part.

\subsection{Prototype (b): the time-fractional Schr\"odinger equation}\label{sec:exp-schroedinger}

In this prototype the values follow from the theory alone. Let $A_V=-\frac{d^2}{dx^2}+V$ be a Schr\"odinger operator that is self-adjoint,
bounded below and unbounded above, with at least one positive spectral value, for instance the operator on $\R$ with the P\"oschl--Teller potential
$V(x)=-6\operatorname{sech}^2x$, whose bound states are $-4$ and $-1$, or its restriction to a large interval with Dirichlet conditions. For $L=e^{-i\frac{\alpha\pi}{2}}A_V$,
\cref{rem:normal-ray} and \cref{prop:normal} give $\sup_t\|E_\alpha(t^\alpha L)\|=\Mst$ and $\Kr_\alpha(L)=\frac{1}{\alpha}$ exactly, so the two sides of
\Cref{thm:A} differ by the factor $\alpha\Mst$, while $L$ generates no $C_0$-semigroup. The numerical range of $L$ lies on the line
through $0$ that contains the lower boundary ray of $\Sa$, so \cref{thm:numrange} applies with the half-plane bounded by this line and
returns the exact value $\Mst$, whereas for bounded discretisations the bound of \citet{CrouzeixPalencia2017} gives $(1+\sqrt2)\Mst$ and the
half-plane estimate 
of \citet{Crouzeix2025} does not apply.
The memory of the time-fractional
Schr\"odinger equation therefore costs an amplification of exactly $\Mst$, which is $2.34$ for $\alpha=\frac{1}{2}$, in a model for which
no semigroup bound exists. Finite-difference discretisations reproduce these values on every grid, while the norms $\|e^L\|$ of the
discretised semigroups grow without bound with the resolution (Table~S9).

\subsection{Prototype (c): time-changed Markov dynamics}\label{sec:exp-markov}

\emph{Compartment chains.} Let $\mathbf x_s$ be an Ornstein--Uhlenbeck process $d\mathbf x_s=-B\mathbf x_s\,ds+\Xi\,dW_s$ in $\R^d$ and let $\mathcal T_t$ be the
inverse of an independent standard $\alpha$-stable subordinator, the clock of a continuous-time random walk with heavy-tailed
waiting times. Since $\mathcal T_t$ has the law of $t^\alpha Y_\alpha$ \citep{Bingham1971,MeerschaertScheffler2004} and $\mathbb E\mathbf x_s=e^{-sB}\mathbf x_0$,
\cref{lem:families}(b) gives
\[
\mathbb E\big[\mathbf x_{\mathcal T_t}\big]=E_\alpha(-t^\alpha B)\,\mathbf x_0,
\]
so the mean of the time-changed process is governed by \Cref{thm:A}. We take the linear chain $B=I-\mu Z_d^\top$, with $Z_d$ the
nilpotent shift of \cref{sec:mechanisms}, in which each compartment decays at unit rate and feeds the next one with strength $\mu$;
then $-B=J_{d,\mu}^\top$ is a Jordan block, and the classical mean exhibits the large transient growth that is typical of
feed-forward networks. In this example the two Kreiss constants coincide: since $\mu\ge0$ and $(z+B)^{-1}=\sum_k\mu^k(Z_d^\top)^k(z+1)^{-k-1}$, its entries
are dominated in modulus by those of the nonnegative matrix $(\RE z+B)^{-1}$, and together with $a_\alpha(z)\le\RE z$ this gives
$\Kr_\alpha(-B)=\sup_{x>0}x\|(x+B)^{-1}\|=\Kr_1(-B)$ for every $\alpha$. The difference between the classical and the fractional transient is
therefore entirely in the propagators. For $\mu=2$ the classical ratio $\frac{\sup_s\|e^{-sB}\|}{\Kr_1(-B)}$ grows like $e\sqrt{\frac{d}{2\pi}}$, from
$1.88$ at $d=5$ to $15.3$ at $d=200$, while the fractional ratio at $d=200$ equals $1.09$, $1.31$ and $2.66$ for $\alpha=0.3$, $0.5$ and
$0.8$, just below the bounds $\vartheta_\alpha(\frac{d}{d-1})^d$ that \cref{prop:jordan}(c) gives in the limit of strong coupling
(Table~S10). The rank is $d$ here, so the example illustrates the mechanism of \cref{sec:mechanisms} rather than the
slow-component bound: the heavy-tailed clock averages out the nilpotent amplification, and up to $d=200$ the ratio stays below $2.7$.

\emph{Learned low-rank generators.} In data-driven practice the generator of a diffusion is replaced by a finite-rank estimate on an
infinite-dimensional function space, and the slow-component formulation then yields certificates whose constants depend on the rank
of the estimate and not on the number of features. For the transfer operator of a non-reversible Langevin diffusion in $\R^2$,
estimated by reduced-rank regression \citep{kostic2022koopman} of rank $r\le16$ on up to $4000$ random Fourier features, every
quantity of \cref{thm:slow} is evaluated at a cost governed by $r$. The estimated slow dynamics shows no transient growth beyond the
obliqueness of the estimator, $\Kr_\alpha(L_r;\Pi)=\sup_t\|E_\alpha(t^\alpha L_r)\Pi\|=\|\Pi\|$, so the lower bound of \Cref{thm:A} is attained,
and in one configuration a spurious eigenvalue of size $10^{-6}$ in $\Sa$ makes $\Kr_\alpha(L_r;\Pi)$ infinite and flags the model as
unusable for long-time forecasting under the fractional clock (Section~S7.2 of Supplementary material).

\subsection{Summary}

In the three prototypes new theory gives information where the classical tools give little or none. The deflated Kreiss constant
follows the propagator within $13\,\%$ in a model where every semigroup and numerical-range bound is infinite, the amplification
caused by memory in the time-fractional Schr\"odinger equation is exactly $\Mst$ although no semigroup exists, and nilpotent growth
that is unbounded in physical time is tamed by the fractional clock. Among the operators considered here, the ratio of the two sides of \Cref{thm:A} is large only in
the oscillatory construction of \cref{sec:optimal}, where it reaches $13.4$ at rank $32$, at which random operators stay below $5.6$,
and $138$ at rank $512$, in agreement with the linear upper bound of \Cref{thm:A} and the lower bound of \Cref{thm:C}.

\section{Conclusions and open problems}\label{sec:conclusions}

The fractional Kreiss constant is the resolvent quantity that controls the transient growth of Mittag--Leffler families. It bounds
the propagator from below for every generator of a bounded family, and from above, up to a factor linear in the rank, for every
generator with a finite slow component, on arbitrary Banach spaces and uniformly in the order $\alpha\in(0,1]$. We show that this rank dependence is genuine, 
and that in infinite dimensions a finite Kreiss constant does not guaranty neither boundedness nor existence of the family. 
The upper bounds rest on the Bernstein-type inequality for a class of non-rational functions, which can be of independent interest. 

The bigger picture is that the fractional transient also behaves differently from the classical one. Normal operators amplify, modes that are unstable in physical time are stable under the fractional clock, and for each fixed $\alpha<1$ the nilpotent amplification that drives the classical Kreiss matrix theorem is averaged out by the Mittag--Leffler distribution, while the growth in our optimality construction is oscillatory. The envelope has distinct short-time, intermediate and long-time regimes, governed by the numerical abscissa, by the
Kreiss constant and by the resolvent at the origin.

These results bear directly on the models that motivated them. For subdiffusion driven by a self-excited oscillator, the
deflated Kreiss constant tracked the transient within $13\,\%$ in all computed cases, although the semigroup grows exponentially and
the numerical range meets the sector of instability, which makes it a natural candidate for assessing how strongly an unstable
local mode can amplify a subdiffusive signal. For the time-fractional Schr\"odinger equation, the solution operators are bounded by exactly $\Mst$ although no
semigroup exists, which quantifies the price of memory in this quantum model. For time-changed Markov dynamics, the means of
compartment chains are far less amplified under the fractional clock than in physical time, and for finite-rank models learned from data the slow-component formulation yields certificates and diagnostics whose constants are set by the number of retained modes.

We conclude with an open question that remains for future work.  It concerns the sharpness of constants $\Jq_F$ in \Cref{thm:B} and the order of $\kappa^{\rm opt}_\alpha(n)$. For $\alpha=1$ the linear order in the dimension is sharp \citep{Kraaijevanger1994}; for $\alpha<1$ we know $n^\gamma\lesssim\kappa^{\rm opt}_\alpha(n)\lesssim n$ for every $\gamma<1$. Finally, we note that the proof of \Cref{thm:B} uses change $\omega\mapsto\omega^\alpha$ only when it maps vertical lines to curves whose tangent turns monotonically through an angle less than $\pi$. Hence, we believe that Kreiss-type bounds for resolvent families of Volterra equations with other completely positive kernels, as well as for the families $t^{\beta-1}E_{\alpha,\beta}(t^\alpha L)$ that appear in inhomogeneous problems, could follow as a natural extension of this work.

\section*{Data availability}
No external data were used. The numerical methods and the complete results are described in the Supplementary Material. The
code that produces them, together with the stored results, is available to the referees on request and will be released in a
public repository upon publication; every
table and figure of the paper and of the Supplementary Material is regenerated by a single script with fixed random seeds.

\section*{Declaration of competing interest}
The authors declare that they have no known competing financial interests or personal relationships that could have appeared to
influence the work reported in this paper.

\section*{Funding}  
This work was funded by the Provincial Secretariat for Higher Education and Scientific Research of Vojvodina (Research Grant 003870560 2025 09418 003 000 000 001 04 004).

\section*{Declaration of generative AI and AI-assisted technologies in the writing process}  
The problem of fractional Kreiss-type theorem for finite rank operators, conceptualization of novel Kreiss constant and developed mathematical approach to obtain the bound is human made, Claude Opus 5.5 has been used
to improve presentation, polish writing and identify issues in consistent notations and mathematical formalism. Additionally, Claude Code has been used to produce simple numerical demonstrations provided in the Supplementary material for reader's convenience.

\bibliographystyle{plainnat}
\bibliography{refs}


\appendix

\begin{center}
    \Large Appendix
\end{center}

The appendix complements the main text in three self-contained parts. \Cref{app:notation} collects the notation, grouped by topic. \Cref{app:aux} contains the proofs of the two results on sectorial generators, the
generation result \cref{lem:sectorial} and the long-time expansion of \cref{prop:timescales}(d); both follow the semigroup case
closely, with the Hankel contour adapted to the fractional order. \Cref{app:optimal} contains the harmonic-analysis ingredients
of the optimality construction of \cref{sec:optimal}: the $A_2$ property of the weight, the ray integral \eqref{eq:ray}, the
boundary behaviour \eqref{eq:bstar} of $E_\alpha$ and the Dirichlet-kernel estimate \eqref{eq:dirichlet}. Each part restates what it
proves, so that it can be read independently of the main text. The complete tables of numerical results, as well as the code used to produce them, are in the separate Supplementary Material.

\section{Notation}\label{app:notation}

\begin{center}\small
\renewcommand{\arraystretch}{1.15}
\begin{longtable}{@{}p{0.27\textwidth}p{0.66\textwidth}@{}}
\toprule
Symbol & Meaning \\
\midrule
\endhead
\bottomrule
\endfoot
\multicolumn{2}{@{}l}{\emph{General notation (\cref{sec:prelim})}}\\
$X$, $X^*$, $H$, $\B(X)$ & Banach space, its dual, Hilbert space, bounded operators\\
$L$, $D(L)$, $\rho(L)$, $\sigma(L)$ & closed, densely defined operator, its domain, resolvent set and spectrum\\
$\alpha\in(0,1]$ & order of the Caputo derivative $\partial_t^\alpha$\\
$\HH$, $\Sigma(\theta)$ & open right half-plane $\{\RE\omega>0\}$; open sector $\{z\ne0:|\Arg z|<\theta\}$\\
$\Sa$, $\La$ & open sector of instability $\Sigma(\frac{\alpha\pi}{2})$; closed stability region $\C\setminus\Sa$\\
$a_\alpha(z)$ & weight $|z|\cos(\frac{\Arg z}{\alpha})$ on $\Sa$; $a_\alpha(\omega^\alpha)=\RE(\omega)|\omega|^{\alpha-1}$\\
$n$, $m$, $d$ & rank of an operator, rank of a slow component, dimension\\
\multicolumn{2}{@{}l}{\emph{Mittag--Leffler function and families (\cref{sec:prelim})}}\\
$E_\alpha(z)$, $E_\alpha(t^\alpha L)$ & Mittag--Leffler function; Mittag--Leffler family generated by $L$ (\cref{def:family})\\
$q_\alpha(z)$ & $E_\alpha(z)-\alpha^{-1}\exp(z^{\frac{1}{\alpha}})$ (\cref{lem:ML})\\
$\Theta_\alpha$, $g_\alpha$ & one-sided $\alpha$-stable random variable and its density\\
$Y_\alpha$, $m_\alpha$, $\vartheta_\alpha$ & Mittag--Leffler random variable $\Theta_\alpha^{-\alpha}$, its density, $\sup_u u\,m_\alpha(u)$\\
$\Mst$, $M_\alpha(\theta)$ & $\sup_{\La}|E_\alpha|$; $\sup_{r\ge0}|E_\alpha(re^{i\theta})|$\\
$s_*$, $r_*$ & maximiser of $|E_\alpha|$ on the boundary ray, $s_*=r_*^{\frac{1}{\alpha}}$\\
$k_\beta(t)$ & Volterra kernel $\frac{t^{\beta-1}}{\Gamma(\beta)}$\\
$\psi$, $\psi'$, $\varphi$ & sectoriality angles and the Hankel angle $\varphi=\frac{\psi'}{\alpha}$ (\cref{lem:sectorial})\\
$s_{\psi'}(L)$, $\mathrm{Ha}_\varphi$, $c_\varphi$ & sectoriality constant, Hankel contour and its constant (\cref{lem:sectorial})\\
\multicolumn{2}{@{}l}{\emph{Kreiss constants and slow components (\cref{sec:prelim,sec:finite-rank,sec:slow})}}\\
$\Kr_\alpha(L)$, $\Kr_1(L)$ & fractional and classical Kreiss constants (\cref{def:kreiss})\\
$\Kr_\alpha^+(L)$ & alternative normalisation of \cref{rem:Kplus}\\
$\sigma_\eps(L)$, $a_{\alpha,\eps}(L)$ & $\eps$-pseudospectrum; supremum of $a_\alpha$ over $\sigma_\eps(L)\cap\Sa$ (\cref{lem:elementary})\\
$P$, $Q=I-P$, $L_P$, $L_Q$ & slow component, its complement, slow block, fast part (\cref{def:slow})\\
$\Kr_\alpha(L;P)$ & deflated Kreiss constant of a slow component (\cref{def:slow})\\
$C_Q$ & $\sup_t\|E_\alpha(t^\alpha L)Q\|$, size of the fast propagator\\
$U$, $V$, $\widehat L=VU$ & factors of a finite-rank operator $L=UV$ (\cref{lem:struct})\\
$\Pi=U\widehat L^{-1}V$ & slow projection of a finite-rank operator, onto $\Ran L$ along $\ker L$\\
$\widehat S(\omega)$, $\widehat g(\omega)$ & Laplace transforms of a family and of a matrix element $g(t)=\langle E_\alpha(t^\alpha L)x,x^*\rangle$\\
$\kappa_{\alpha,n}$, $\kappa'_{\alpha,n}$ & $\frac e2(3+4n(2-\alpha))$, $\frac e2(1+4n(2-\alpha))$\\
$\kappa^{\rm opt}_\alpha(n)$ & best constant in \cref{thm:finite-rank} on Hilbert spaces\\
$\mathbf n$ & inner unit normal of a boundary ray of $\Sa$\\
\multicolumn{2}{@{}l}{\emph{The Bernstein-type inequality (\cref{sec:bernstein})}}\\
$\RR_{\alpha,n}$ & nonzero rational functions $F=\frac{F_1}{F_0}$, coprime, $\deg F_1\le\deg F_0\le n$, no poles in $\Sa$\\
$h_F$, $D_F$, $\Jq_F(a)$ & $\omega^{-1}F(\omega^\alpha)$; $-F+\alpha zF'$; length-to-supremum quotient\\
$\ell_a$, $h^*_a$ & length of $y\mapsto h_F(a+iy)$ and its maximal modulus\\
$\Var_I(u)$, $N_u$, $N_u^{[c]}$ & total variation, counting function, counting function modulo $c$\\
$\mathcal C_\alpha$, $T$, $\Theta$ & image of $\{\RE\omega=1\}$ under $\omega\mapsto\omega^\alpha$, tangent angle, total turning\\
$\eta$, $\beta_\zeta$ & continuous lifts of $\arg h_F'$ and of $\arg((1+iy)^\alpha-\zeta)$\\
$G$, $G_N$, $G_D$ & numerator polynomial of $D_F$ and the reduced numerator and denominator (proof of \cref{thm:bernstein})\\
\multicolumn{2}{@{}l}{\emph{Hilbert spaces (\cref{sec:hilbert})}}\\
$W(L)$, $C_W$, $\nu(L)$ & numerical range, the closed convex cone it generates, numerical abscissa\\
$C_Q^W$ & closed convex cone generated by $W(L_Q)$ (\cref{cor:deflation})\\
$\Omega$, $\Omega_\phi$, $\chi$, $\delta$ & closed half-planes in \cref{sec:hilbert}, the angle and offset of $\Omega$\\
$\mathcal G$ & skew-adjoint generator of a unitary dilation\\
$G_\alpha(s)$ & $\sup_{y\in\R}|E_\alpha(-s+iy)|$\\
$K(\beta)$ & spectral constant of \citet{Crouzeix2025} for convex sets of aperture $2\beta$\\
\multicolumn{2}{@{}l}{\emph{Optimality (\cref{sec:optimal})}}\\
$\T$, $e_k$, $\mathcal P$ & circle, $e^{ikx}$, trigonometric polynomials\\
$w_\gamma$, $L^2_\gamma$, $\|\cdot\|_\gamma$ & weight $|\tan(\frac{x}{2})|^\gamma$, weighted space and its norm\\
$R_j$, $r_\gamma$, $S_N$, $V_N$ & Riesz projections, their uniform bound, Fourier partial sum projection and its range\\
$T_m$, $\Var_\Z(m)$ & Fourier multiplier with symbol $m=(m_k)$ and the variation of the symbol\\
$\mathsf R_\theta$ & ray $\{se^{i\theta}:s\ge0\}$\\
$\Phi_\alpha$, $b_\alpha$, $b^*_\alpha$ & $E_\alpha((is)^\alpha)$, its deviation from $\alpha^{-1}e^{is}$, $\sup|b_\alpha|+\Var(b_\alpha)$\\
$D_N$, $c_\gamma$, $Q_\gamma$ & Dirichlet kernel and the constants in \eqref{eq:dirichlet}\\
$\D^\alpha$, $B_N$, $\ell_N$ & Weyl fractional derivative; its truncation $N^{-\alpha}\D^\alpha S_N$ and the symbol of $B_N$\\
\multicolumn{2}{@{}l}{\emph{Mechanisms (\cref{sec:mechanisms})}}\\
$Z_n$, $J_{n,\mu}$ & nilpotent shift, Jordan block $-I+\mu Z_n$\\
$\Delta_n(\alpha)$, $\varrho_n(\alpha)$ & normalised peak of $x^{n-1}|E_\alpha^{(n-1)}(-x)|$, Jordan-block limit (\cref{prop:jordan})\\
\multicolumn{2}{@{}l}{\emph{Prototypes (\cref{sec:slow,sec:hilbert,sec:experiments})}}\\
$L_1$, $\lambda_\pm$, $e_\pm$, $\theta$ & oscillator matrix of prototype (a), its eigenvalues $e^{\pm i\theta}$ and eigenvectors\\
$A$, $\mathcal S$, $b_1$, $b_2$, $\mu$ & Dirichlet Laplacian, source operator, source profiles and coupling strength in prototype (a)\\
$Y=\mu Y_0$, $y_0$, $s_0$ & Sylvester solution, norm and smallest singular value of $Y_0$ (\cref{prop:cascade})\\
$A_V$ & Schr\"odinger operator $-\frac{d^2}{dx^2}+V$ in prototype (b)\\
$\mathbf x_s$, $B=I-\mu Z_d^\top$, $\mathcal T_t$ & Ornstein--Uhlenbeck process, its drift matrix (a compartment chain) and the inverse stable subordinator\\
$L_r$, $r$, $D$, $\tau$ & learned rank-$r$ generator, rank, number of features, lag\\
\end{longtable}
\end{center}

\section{Proofs of preliminary results
}\label{app:aux}

\subsection{Mittag--Leffler functions}

Here we provide the collection of useful facts on Mittag--Leffler functions, together with their proofs for completeness.

\begin{lemma}\label{lem:ML}
Let $0<\alpha<1$ and fix $\theta_0$ with $\frac{\alpha\pi}{2}<\theta_0<\alpha\pi$. There are constants $C,C'>0$ such that the
following hold.
\begin{enumerate}[label=\textup{(\alph*)}]
\item If $|z|\ge1$ and $|\Arg z|\le\theta_0$, then
$E_\alpha(z)=\alpha^{-1}\exp(z^{\frac{1}{\alpha}})-\frac{z^{-1}}{\Gamma(1-\alpha)}+r(z)$ with $|r(z)|\le C|z|^{-2}$.
\item If $|z|\ge1$ and $\theta_0\le|\Arg z|\le\pi$, then $|E_\alpha(z)+\frac{z^{-1}}{\Gamma(1-\alpha)}|\le C|z|^{-2}$.
\item $E_\alpha$ is bounded on $\La$; on every closed sector $\{\theta_1\le|\Arg z|\le\pi\}$ with $\theta_1>\frac{\alpha\pi}{2}$,
$E_\alpha^{(j)}(z)=O(|z|^{-1-j})$ as $|z|\to\infty$ for every $j\in\N_0$; and $q_\alpha(z):=E_\alpha(z)-\alpha^{-1}\exp(z^{\frac{1}{\alpha}})$
satisfies $|q_\alpha(z)|\le C'|z|^{-1}$ and $|q_\alpha'(z)|\le C'|z|^{-2}$ for $|\Arg z|=\frac{\alpha\pi}{2}$ and $|z|\ge2$.
\item $\int_0^\infty u^km_\alpha(u)\,du=\frac{k!}{\Gamma(1+\alpha k)}$ for $k\in\N_0$, and $\int_0^\infty e^{xu}m_\alpha(u)\,du=E_\alpha(x)<\infty$
for every $x\in\R$.
\item $\vartheta_\alpha:=\sup_{u>0}u\,m_\alpha(u)<\infty$, and $m_{\frac{1}{2}}(u)=\pi^{-\frac{1}{2}}e^{-\frac{u^2}{4}}$, so that $\vartheta_{\frac{1}{2}}=\sqrt{\frac{2}{e\pi}}$.
\end{enumerate}
\end{lemma}

\begin{proof}
(a) and (b) are the cases $\beta=1$, $p=1$ of the classical expansions of the two-parameter function
$E_{\alpha,\beta}$ \citep[Section~1.2]{Podlubny1999}, see also \citet[Chapter~4]{Gorenflo2020}.

(c) For $z$ with $\frac{\alpha\pi}{2}\le|\Arg z|\le\theta_0$ we have $\frac{|\Arg z|}{\alpha}\in[\frac{\pi}{2},\pi)$, and therefore
$|\exp(z^{\frac{1}{\alpha}})|=\exp(|z|^{\frac{1}{\alpha}}\cos(\frac{\Arg z}{\alpha}))\le1$; together with (a), (b) and continuity on
$\{|z|\le1\}$ this gives boundedness on $\La$. On $\{\theta_1\le|\Arg z|\le\pi\}$ the exponential term decays faster
than any power, so $E_\alpha(z)=O(|z|^{-1})$ there, and Cauchy's estimates on discs of radius proportional to $|z|$
inside a slightly larger sector give the bounds on the derivatives. By (a), $|q_\alpha(z)|\le C'|z|^{-1}$ for $|z|\ge1$ and
$|\Arg z|\le\theta_0$, where $q_\alpha$ is holomorphic since $\theta_0<\pi$. Let $|\Arg z|=\frac{\alpha\pi}{2}$, $|z|\ge2$ and
$\varrho:=\frac{|z|\sin(\theta_0-\frac{\alpha\pi}{2})}{2}$. Every $\zeta$ with $|\zeta-z|\le\varrho$ satisfies $|\Arg\zeta|\le\theta_0$ and $|\zeta|\ge\frac{|z|}{2}\ge1$, so
Cauchy's estimate gives
\[
|q_\alpha'(z)|\le\varrho^{-1}\max_{|\zeta-z|=\varrho}|q_\alpha(\zeta)|=O(|z|^{-2}),
\]
and it remains to enlarge $C'$.

(d) The first identity is the moment formula $\mathbb E\Theta_\alpha^{-\alpha k}=\frac{\Gamma(1+k)}{\Gamma(1+\alpha k)}$
\citep[Section~XIII.8]{Feller1971}, see also \citet{Pollard1948}. Since $\sum_k\frac{x^k\,\mathbb EY_\alpha^k}{k!}=E_\alpha(x)$, the second
follows for $x\ge0$ by monotone convergence and for $x<0$ by dominated convergence.

(e) With $s=u^{-\frac{1}{\alpha}}$ we have $u\,m_\alpha(u)=\alpha^{-1}s\,g_\alpha(s)$, which is bounded because $g_\alpha$ is bounded
and continuous with $g_\alpha(s)=O(s^{-1-\alpha})$ as $s\to\infty$ \citep[Section~XVII.6]{Feller1971}. For $x\ge0$ one checks
after some algebra that $\pi^{-\frac{1}{2}}\int_0^\infty e^{-xu-\frac{u^2}{4}}\,du=e^{x^2}\operatorname{erfc}(x)=E_{\frac{1}{2}}(-x)$, and a distribution on
$(0,\infty)$ is determined by its Laplace transform.
\end{proof}

\subsection{Sectorial generators}

Here we provide the proof of  \cref{lem:families}.

\begin{proof}
(a) Norm continuity and the bound are immediate from the power series. For $\RE\omega>0$ we have
$\int_0^\infty e^{-\omega t}t^{\alpha k}dt=\Gamma(1+\alpha k)\,\omega^{-1-\alpha k}$, and
$\sum_k\|L\|^k(\RE\omega)^{-1-\alpha k}<\infty$ when $\|L\|<(\RE\omega)^\alpha$. Termwise integration gives
$\sum_kL^k\omega^{-1-\alpha k}=\omega^{-1}(I-\omega^{-\alpha}L)^{-1}=\omega^{\alpha-1}(\omega^\alpha-L)^{-1}$.

(b) By \cref{lem:ML}(d),
\[
\int_0^\infty m_\alpha(u)\|e^{t^\alpha uL}\|\,du\le ME_\alpha(\varpi t^\alpha)<\infty ,
\]
and $E_\alpha(\varpi t^\alpha)=O(e^{\varpi^{\frac{1}{\alpha}}t})$ by \cref{lem:ML}(a), so the integrability condition of \cref{def:family} holds
for every $w>\varpi^{\frac{1}{\alpha}}$.
Hence $S(t)x:=\int_0^\infty m_\alpha(u)e^{t^\alpha uL}x\,du$ converges absolutely, $S(0)=I$,
$\|S(t)\|\le\sup_s\|e^{sL}\|$, and $S$ is strongly continuous by dominated convergence. For $\lambda>0$ and
$s>0$ put $\psi_\lambda(s):=\int_0^\infty e^{-\lambda t}t^{-\alpha}m_\alpha(st^{-\alpha})\,dt$. For $\mu\ge0$,
Tonelli's theorem, the substitution $s=t^\alpha u$ and \cref{lem:ML}(d) give
\[
\int_0^\infty e^{-\mu s}\psi_\lambda(s)\,ds=\int_0^\infty e^{-\lambda t}E_\alpha(-\mu t^\alpha)\,dt
=\frac{\lambda^{\alpha-1}}{\lambda^\alpha+\mu}=\int_0^\infty e^{-\mu s}\lambda^{\alpha-1}e^{-\lambda^\alpha s}\,ds.
\]
The second equality holds for $\lambda>\mu^{\frac{1}{\alpha}}$ by (a) on $X=\C$, and then for every $\lambda>0$ by
analytic continuation, since both sides are holomorphic in $\lambda\in\HH$. By uniqueness of the
Laplace transform, $\psi_\lambda(s)=\lambda^{\alpha-1}e^{-\lambda^\alpha s}$ for almost every $s>0$.
Consequently, for $\lambda^\alpha>\varpi$, Fubini's theorem and the same substitution give
\[
\int_0^\infty e^{-\lambda t}S(t)x\,dt=\int_0^\infty\psi_\lambda(s)\,e^{sL}x\,ds
=\lambda^{\alpha-1}\int_0^\infty e^{-\lambda^\alpha s}e^{sL}x\,ds=\lambda^{\alpha-1}(\lambda^\alpha-L)^{-1}x.
\]

(c) By the Lumer--Phillips theorem $\|e^{sL}\|\le1$ \citep{EngelNagel2000}; for $\alpha<1$ apply (b).
\end{proof}



\begin{proof}[Proof of \cref{lem:sectorial}]
If $|\Arg\omega|\le\varphi$ and $\omega\ne0$, then $\omega^\alpha$ lies in the closure of $\Sigma(\psi')$, which is contained in
$\Sigma(\psi)\cup\{0\}$. By continuity of the resolvent,
\begin{equation}\label{eq:sect-bound}
\|\omega^{\alpha-1}(\omega^\alpha-L)^{-1}\|\le s_{\psi'}(L)\,|\omega|^{-1}\qquad(0<|\omega|,\ |\Arg\omega|\le\varphi).
\end{equation}
Let $S(t)$ denote the right-hand side of \eqref{eq:hankel}. On the rays $|e^{\omega t}|=e^{-t\varrho|\cos\varphi|}$, so the integral converges
absolutely for $t>0$ and is norm continuous in $t>0$. The integrand is holomorphic on $\{\omega\ne0:|\Arg\omega|<\varphi+\eps\}$ for small
$\eps>0$, because $\alpha(\varphi+\eps)<\psi$. By Cauchy's theorem the unit arc may be replaced by the arc of radius $\frac{1}{t}$, and the
substitution $\omega=\frac{s}{t}$ together with \eqref{eq:sect-bound} gives $\|S(t)\|\le c_\varphi s_{\psi'}(L)$ for $t>0$.

\emph{Strong continuity at $0$.} Let $x\in D(L)$. Since
$\omega^{\alpha-1}(\omega^\alpha-L)^{-1}x=\omega^{-1}x+\omega^{-1}(\omega^\alpha-L)^{-1}Lx$ and
$(2\pi i)^{-1}\int_{\mathrm{Ha}_\varphi}e^{\omega t}\omega^{-1}d\omega=1$ by Hankel's formula, we have
\[
S(t)x-x=\frac1{2\pi i}\int_{\mathrm{Ha}_\varphi}e^{\omega t}\omega^{-1}(\omega^\alpha-L)^{-1}Lx\,d\omega .
\] For $0<t\le1$ the
integrand is bounded by $e\,s_{\psi'}(L)\|Lx\|\,|\omega|^{-1-\alpha}$, which is integrable on $\mathrm{Ha}_\varphi$, so by dominated convergence
$S(t)x-x$ converges as $t\downarrow0$ to $(2\pi i)^{-1}\int_{\mathrm{Ha}_\varphi}\omega^{-1}(\omega^\alpha-L)^{-1}Lx\,d\omega$. This integral vanishes:
close the contour by the arcs $\{Re^{i\vartheta}:|\vartheta|\le\varphi\}$, on which the integrand is $O(R^{-1-\alpha})$, and apply Cauchy's theorem in
$\{1<|\omega|<R,\ |\Arg\omega|<\varphi\}$. Since $D(L)$ is dense and $S$ is bounded on $(0,1]$, $S(t)x\to x$ for every $x$; we set $S(0):=I$.

\emph{Laplace transform.} Let $\lambda>1$. On $\mathrm{Ha}_\varphi$ we have $\RE\omega\le1<\lambda$ and $\lambda-\RE\omega\ge\lambda+\varrho|\cos\varphi|$ on the rays, so
Fubini's theorem gives
$\int_0^\infty e^{-\lambda t}S(t)\,dt=(2\pi i)^{-1}\int_{\mathrm{Ha}_\varphi}(\lambda-\omega)^{-1}\omega^{\alpha-1}(\omega^\alpha-L)^{-1}\,d\omega$. Truncate $\mathrm{Ha}_\varphi$ at
$|\omega|=R>\lambda$ and close it by the arc $\{Re^{i\vartheta}:|\vartheta|\le\varphi\}$, traversed clockwise, on which the integrand is $O(R^{-2})$. The
resulting closed curve bounds $\{1<|\omega|<R,\ |\Arg\omega|<\varphi\}$ with negative orientation, and the only singularity inside is the simple pole
at $\omega=\lambda$, with residue $-\lambda^{\alpha-1}(\lambda^\alpha-L)^{-1}$. Letting $R\to\infty$ gives \eqref{eq:laplace}.
\end{proof}

\subsection{Time scales in Hilbert spaces}\label{app:timescales}

The envelope $t\mapsto\|E_\alpha(t^\alpha L)\|$ has a short-time, an intermediate and a long-time regime, each governed by a
different resolvent quantity. Let
$\nu(L):=\sup\RE W(L)\in(-\infty,\infty]$ be the numerical abscissa of an operator on $H$.

\setcounter{theorem}{2}
\begin{proposition}[Time scales]\label{prop:timescales}
Let $0<\alpha\le1$ and let $L$ be closed and densely defined on $H$, generating a Mittag--Leffler family.
\begin{enumerate}[label=\textup{(\alph*)}]
\item For every $v\in D(L)$ with $\|v\|=1$, as $t\downarrow0$,
$\|E_\alpha(t^\alpha L)v\|=1+\frac{\RE\langle Lv,v\rangle\,t^\alpha}{\Gamma(1+\alpha)}+o(t^\alpha)$.
\item If $L\in\B(H)$, then $\|E_\alpha(t^\alpha L)\|=1+\frac{\nu(L)\,t^\alpha}{\Gamma(1+\alpha)}+O(t^{2\alpha})$.
\item $\|E_\alpha(t^\alpha L)\|\le1$ for all $t\ge0$ if and only if $L$ is m-dissipative.
\item Let $0<\alpha<1$, and suppose that $L$ satisfies the hypothesis of \cref{lem:sectorial} for some $\psi\in(\frac{\alpha\pi}{2},\pi]$ and that
$0\in\rho(L)$. Then for every $q\in\N$, in operator norm,
\[
E_\alpha(t^\alpha L)=-\sum_{k=1}^{q}\frac{t^{-\alpha k}L^{-k}}{\Gamma(1-\alpha k)}+O\big(t^{-\alpha(q+1)}\big)\qquad(t\to\infty),
\]
and in particular $\|E_\alpha(t^\alpha L)\|=\frac{\|L^{-1}\|\,t^{-\alpha}}{\Gamma(1-\alpha)}+O(t^{-2\alpha})$.
\end{enumerate}
\end{proposition}

\begin{proof}
(a) Let $c_t:=\frac{t^\alpha}{\Gamma(1+\alpha)}=\int_0^tk_\alpha(t-s)\,ds$. By \eqref{eq:integrated},
$\|E_\alpha(t^\alpha L)v-v-c_t Lv\|$ is bounded by $c_t\sup_{0\le s\le t}\|E_\alpha(s^\alpha L)Lv-Lv\|=o(c_t)$, and
$\|v+c_t Lv\|^2=1+2c_t\RE\langle Lv,v\rangle+c_t^2\|Lv\|^2$.

(b) $E_\alpha(t^\alpha L)=I+c_t L+O(t^{2\alpha})$, and taking the supremum over unit vectors in the identity of (a) gives
$\|I+c_t L\|=1+c_t\nu(L)+O(c_t^2)$.

(c) If $L$ is m-dissipative, apply \cref{lem:families}(c). Conversely, (a) gives $\RE\langle Lv,v\rangle\le0$, and \cref{thm:kreiss-lower}
gives $(0,\infty)\subset\Sa\subset\rho(L)$.

(d) Fix $\psi'$ and $\varphi=\frac{\psi'}{\alpha}$ as in \cref{lem:sectorial}, so that \eqref{eq:hankel} holds. The closed set
$\overline{\Sigma(\psi')}=\{z:z=0\text{ or }|\Arg z|\le\psi'\}$ lies in $\rho(L)$, since $0\in\rho(L)$ and $\overline{\Sigma(\psi')}\setminus\{0\}\subset\Sigma(\psi)$.
Hence the quantity $c_L:=\sup\{\|(z-L)^{-1}\|:z\in\overline{\Sigma(\psi')}\}$ is finite: for $|z|\ge1$ the resolvent is bounded by
$\frac{s_{\psi''}(L)}{|z|}$ with $\psi'<\psi''<\psi$, and on the compact remainder it is continuous. For $z\in\rho(L)$ we have $(z-L)^{-1}=-L^{-1}+zL^{-1}(z-L)^{-1}$, and by
induction $(z-L)^{-1}=-\sum_{k=0}^{q-1}z^kL^{-k-1}+z^qL^{-q}(z-L)^{-1}$. Insert this with $z=\omega^\alpha$ into \eqref{eq:hankel}. Hankel's formula
$(2\pi i)^{-1}\int_{\mathrm{Ha}_\varphi}e^{\omega t}\omega^{-\beta}\,d\omega=\frac{t^{\beta-1}}{\Gamma(\beta)}$ with $\beta=1-\alpha(k+1)$ gives the terms of the expansion.
In the remainder substitute $\omega=\frac{s}{t}$ and return to $\mathrm{Ha}_\varphi$ by Cauchy's theorem, as in the proof of \cref{lem:sectorial}; this bounds
its norm by $t^{-\alpha(q+1)}(2\pi)^{-1}c_L\|L^{-q}\|\int_{\mathrm{Ha}_\varphi}|e^{s}|\,|s|^{\alpha(q+1)-1}\,|ds|$.
\end{proof}

Thus, as for semigroups, a Mittag--Leffler family is contractive exactly when its generator is m-dissipative, independently of
$\alpha$, and the initial slope in the variable $t^\alpha$ is $\frac{\nu(L)}{\Gamma(1+\alpha)}$. In physical time the envelope has an infinite slope at
$t=0$ whenever $\alpha<1$ and $\nu(L)>0$. The hypothesis of (d) holds for bounded $L$ with $\sigma(L)\cap\overline{\Sa}=\emptyset$, for normal $L$
with $\sigma(L)\subset\{|\Arg z|\ge\psi\}\setminus\{0\}$, and for $L$ on $H$ with $0\in\rho(L)$, $W(L)\subset\{|\Arg z|\ge\psi\}$ and
$\rho(L)\cap\Sa\ne\emptyset$, by the spectral theorem and the argument of \cref{thm:numrange}(a).

Together, \cref{prop:timescales} and \Cref{thm:A} describe the regimes of the envelope: the initial slope is
governed by the numerical abscissa, the transient maximum by the fractional Kreiss constant, and the algebraic decay by the
resolvent at the origin, $\|L^{-1}\|$, rather than by the spectral abscissa. For $\alpha=1$ all coefficients $\frac{1}{\Gamma(1-k)}$ vanish,
consistently with exponential decay. The envelopes in \cref{fig:cascade}(a) display all of them.
\section{Ingredients of the optimality construction}\label{app:optimal}

Throughout this appendix $0<\gamma<1$, $w:=w_\gamma$ and $\|\cdot\|:=\|\cdot\|_\gamma$. Recall that $w\in A_2(\T)$ if
$\sup_I\big(|I|^{-1}\int_Iw\big)\big(|I|^{-1}\int_Iw^{-1}\big)<\infty$, the supremum being taken over all arcs $I\subset\T$.

We recall \cref{lem:multiplier}(a): $w_\gamma\in A_2(\T)$, the Riesz projections $R_j$ are uniformly bounded on $L^2_\gamma$ with bound
$r_\gamma<\infty$, and $\|S_N\|\le2r_\gamma$.

\begin{proof}[Proof of \cref{lem:multiplier}(a)]
Since $w$ and $w^{-1}$ are integrable, the $A_2$ quotient is bounded over arcs of length at least $\frac{\pi}{2}$. An arc of length less than
$\frac{\pi}{2}$ is contained in $\{|x|\le\frac{3\pi}{4}\}$ or in $\{|x|\ge\frac{\pi}{4}\}$. As $u\le\tan u\le3u$ on $[0,\frac{3\pi}{8}]$, on the first set $w(x)$ is comparable
to $|x|^\gamma$, and on the second, with $y:=\pi-|x|$, to $y^{-\gamma}$. Power weights $|t|^\beta$ with $|\beta|<1$ belong to $A_2(\R)$
\citep[Example~7.1.7]{Grafakos2014}, and $A_2$ bounds are stable under multiplication by functions bounded above and below. Hence
$w\in A_2(\T)$. By the theorem of \citet{HuntMuckenhouptWheeden1973}, the conjugate-function operator $f\mapsto\tilde f$,
$\tilde e_k=-i\operatorname{sgn}(k)e_k$, is bounded on $L^2_\gamma$. For $f\in\mathcal P$,
\[
R_0f=\tfrac12\big(f+i\tilde f\big)+\tfrac12\hat f(0)e_0,\qquad
|\hat f(0)|\,\|e_0\|\le\frac{1}{2\pi}\Big(\int_{-\pi}^{\pi}w^{-1}\Big)^{\frac{1}{2}}\Big(\int_{-\pi}^{\pi}w\Big)^{\frac{1}{2}}\|f\|,
\]
the latter by the Cauchy--Schwarz inequality, so $R_0$ is bounded. Finally, $R_j=e_jR_0e_{-j}$ and multiplication by $e_{\pm j}$ is an
isometry of $L^2_\gamma$, so $r_\gamma<\infty$, and $\|S_N\|\le\|R_{-N}\|+\|R_{N+1}\|\le2r_\gamma$.
\end{proof}

We recall \eqref{eq:ray}: for a ray $\mathsf R_\theta=\{se^{i\theta}:s\ge0\}$ and $z\notin\mathsf R_\theta$,
$\int_0^\infty|z-se^{i\theta}|^{-2}ds\le\frac{\pi}{\dist(z,\mathsf R_\theta)}$.

\begin{proof}[Proof of \eqref{eq:ray}]
Write $ze^{-i\theta}=u+iv$; the integral equals $\int_0^\infty\big((s-u)^2+v^2\big)^{-1}ds$. If $v=0$, then $u<0$ and the integral is
$\frac{1}{|u|}=\frac{1}{\dist(z,\mathsf R_\theta)}$. If $v\ne0$, it equals $|v|^{-1}\big(\frac\pi2+\arctan(\frac{u}{|v|})\big)$. For $u\ge0$ we have
$\dist(z,\mathsf R_\theta)=|v|$, and the bound is clear. For $u<0$ we have $\dist(z,\mathsf R_\theta)=|z|$, and with
$\phi:=\arctan(\frac{|v|}{|u|})\in(0,\frac{\pi}{2})$ the integral equals $\frac{\phi}{|z|\sin\phi}\le\frac{\pi}{2|z|}$.
\end{proof}

We recall \eqref{eq:bstar}: $b_\alpha(s)=E_\alpha((is)^\alpha)-\alpha^{-1}e^{is}$ is continuous on $\R$, bounded and of bounded
variation, and $b_1\equiv0$.

\begin{proof}[Proof of \eqref{eq:bstar}]
For $\alpha=1$ we have $\Phi_1(s)=e^{is}$. Let $0<\alpha<1$. Since $E_\alpha$ has real Taylor coefficients and $(-is)^\alpha=\overline{(is)^\alpha}$ for
$s>0$, we have $b_\alpha(-s)=\overline{b_\alpha(s)}$, so it suffices to work on $[0,\infty)$. Put $\sigma_\alpha:=2^{\frac{1}{\alpha}}$. For $s\ge\sigma_\alpha$ the point
$z=(is)^\alpha=s^\alpha e^{i\frac{\alpha\pi}{2}}$ satisfies $|z|\ge2$, $|\Arg z|=\frac{\alpha\pi}{2}$ and $z^{\frac{1}{\alpha}}=is$. Hence $b_\alpha(s)=q_\alpha(z)$, and
\cref{lem:ML}(c) gives
\[
|b_\alpha(s)|\le C's^{-\alpha},\qquad|b_\alpha'(s)|=|q_\alpha'(z)|\,\alpha s^{\alpha-1}\le\alpha C's^{-1-\alpha}\qquad(s\ge\sigma_\alpha).
\]
On $[0,\sigma_\alpha]$ we have $\Phi_\alpha(s)=\sum_{k\ge0}\frac{e^{i\frac{\alpha k\pi}{2}}s^{\alpha k}}{\Gamma(1+\alpha k)}$, so
$\int_0^{\sigma_\alpha}|\Phi_\alpha'|\le\sum_{k\ge1}\frac{\sigma_\alpha^{\alpha k}}{\Gamma(1+\alpha k)}=E_\alpha(2)-1$, while $s\mapsto\alpha^{-1}e^{is}$ has variation $\frac{\sigma_\alpha}{\alpha}$ there.
Thus $b_\alpha$ is continuous and bounded on $[0,\infty)$, and
$\Var_{[0,\infty)}(b_\alpha)\le E_\alpha(2)-1+\frac{\sigma_\alpha}{\alpha}+\alpha C'\int_{\sigma_\alpha}^\infty s^{-1-\alpha}\,ds<\infty$.
\end{proof}

We recall \eqref{eq:dirichlet}: $\|D_N(\cdot+\pi)\|_\gamma\ge c_\gamma N^\gamma\|D_N\|_\gamma$ for all $N\in\N$, with $c_\gamma$ and $Q_\gamma$ as stated there.

\begin{proof}[Proof of \eqref{eq:dirichlet}]
Since $\sin(\frac{y}{2})\ge\frac{y}{\pi}$ on $[0,\pi]$, we have $|D_N(y)|\le\min\big(2N+1,\frac{\pi}{|y|}\big)$ for $0<|y|\le\pi$, and $|D_N(y)|\le\sqrt2$ for
$\frac{\pi}{2}\le|y|\le\pi$. Since $\tan u\le2u$ on $[0,\frac{\pi}{4}]$, we have $w(y)\le|y|^\gamma$ for $|y|\le\frac{\pi}{2}$. Splitting at $|y|=\frac{\pi}{2}$, and the integral
over $[0,\frac{\pi}{2}]$ at $y_N:=\frac{1}{2N+1}$,
\[
\|D_N\|^2\le2\int_0^{y_N}(2N+1)^2y^\gamma\,dy+2\int_{y_N}^{\frac{\pi}{2}}\pi^2y^{\gamma-2}\,dy+4\int_{\frac{\pi}{2}}^{\pi}\tan^\gamma\frac y2\,dy\le Q_\gamma(2N+1)^{1-\gamma}.
\]
Next, the substitution $y=x+\pi$ and $|\tan(\frac{y-\pi}{2})|=|\cot(\frac{y}{2})|$ give
\[
\|D_N(\cdot+\pi)\|^2=\int_{-\pi}^{\pi}|D_N(y)|^2\,|\cot(\tfrac{y}{2})|^\gamma\,dy .
\]
For $0<|y|\le y_N$ we have $(N+\tfrac12)|y|\le\tfrac12$, hence $D_N(y)\ge\frac2\pi(2N+1)\ge N$, using $\sin u\ge\frac{2u}{\pi}$ on $[0,\frac{\pi}{2}]$ and
$\sin u\le u$, and moreover $|\cot(\frac{y}{2})|\ge\frac{1}{|y|}$. Therefore
$\|D_N(\cdot+\pi)\|^2\ge2N^2\int_0^{y_N}y^{-\gamma}\,dy=\frac{2N^2(2N+1)^{\gamma-1}}{1-\gamma}$. Dividing the two estimates and using $2N+1\le3N$ and
$2\gamma-2<0$,
\[
\frac{\|D_N(\cdot+\pi)\|^2}{\|D_N\|^2}\ge\frac{2N^2(2N+1)^{2\gamma-2}}{(1-\gamma)Q_\gamma}\ge\frac{2\cdot3^{2\gamma-2}}{(1-\gamma)Q_\gamma}\,N^{2\gamma}=c_\gamma^2N^{2\gamma}.\qedhere
\]
\end{proof}
\end{document}